\documentclass[11pt]{article}

\usepackage[T1]{fontenc}
\usepackage[latin1]{inputenc}
\usepackage{amsfonts,amsmath,amssymb,amsthm,color,enumerate,graphicx}
\usepackage{dsfont}
\usepackage[colorlinks,urlcolor=blue]{hyperref}
\usepackage{bbm}
\usepackage{fullpage}
\usepackage[numbers]{natbib}
\newcommand{\abs}[1]{\lvert #1 \rvert}
\newcommand{\norm}[1]{\lVert #1 \rVert}
\newcommand{\set}[1]{\ensuremath{\{#1\}}}
\newcommand{\cutnorm}[1]{\norm{#1}_{\square}}
\newcommand{\dcut}{d_{\square}}

\newcommand{\Kcut}{\square_{K}}
\newcommand{\Kcutnorm}[1]{\norm{#1}_{\Kcut}}
\newcommand{\dKcut}{d_{\Kcut}}
\newcommand{\deltaKcut}{\delta_{\Kcut}}
\newcommand{\NdKcut}[2][\delta]{\widehat{N}_{\Kcut}(n, #1; #2)}
\newcommand{\NdeltaKcut}[2][\delta]{N_{\Kcut}(n, #1; #2)}
\newcommand{\kcut}{\square_{k}}

\newcommand{\dkcut}{d_{\kcut}}
\newcommand{\deltakcut}{\delta_{\kcut}}

\newcommand{\limitsx}[1]{\widehat{#1}}
\newcommand{\maxlimits}[1]{\widehat{#1}^*}
\newcommand{\oi}{[0,1]}

\newcommand{\sequence}[2][n]{(#2_{#1})_{{#1} = 1}^{\infty}}
\newcommand{\tend}{\longrightarrow}
\newcommand{\pto}{\overset{\mathrm{p}}{\tend}}

\newcommand{\average}[1]{\overline{#1}}
\renewcommand{\vector}[1]{\mathbf{#1}}

\newcommand{\BorelMeasures}{\mathcal{B}}

\newcommand{\family}{\mathcal{F}}
\newcommand{\GenSet}{\mathcal{F}}

\newcommand{\labeled}{\mathcal{L}}
\newcommand{\PP}{\mathcal{P}}
\newcommand{\measures}{\mathcal{P}}
\newcommand{\partition}{\mathcal{Q}}
\newcommand{\QQ}{\mathcal{Q}}
\newcommand{\SSS}{\mathcal{S}}
\newcommand{\graphs}{\mathcal{U}}
\newcommand{\WW}{\mathcal{W}}

\newcommand{\N}{\mathbb{N}}
\newcommand{\Prob}{\mathbb{P}}
\newcommand{\R}{\mathbb{R}}

\newcommand{\Branch}{\mathbf{B}}
\newcommand{\Graphseq}{\mathbf{G}}
\newcommand{\Grid}{\mathbf{Grid}}
\newcommand{\Complete}{\mathbf{K}}
\newcommand{\Path}{\mathbf{P}}
\newcommand{\Hypercube}{\mathbf{Q}}

\newcommand{\WWk}{\WW_k}

\DeclareMathOperator{\ex}{ex}
\DeclareMathOperator{\supp}{supp}

\DeclareMathOperator{\Binom}{Binom}
\DeclareMathOperator{\Ent}{Ent}
\DeclareMathOperator{\Forb}{Forb}

\numberwithin{equation}{section}

\newtheorem{theorem}{Theorem}[section]
\newtheorem{lemma}[theorem]{Lemma}
\newtheorem{corollary}[theorem]{Corollary}

\newtheorem{claim}[theorem]{Claim}
\newtheorem{proposition}[theorem]{Proposition}

\newtheorem{observation}[theorem]{Observation}

\theoremstyle{definition}
\newtheorem{definition}[theorem]{Definition}
\newtheorem{example}[theorem]{Example}

\newtheorem{problem}[theorem]{Problem}

\theoremstyle{remark}
\newtheorem{remark}[theorem]{Remark}

\title{Multicolour containers and the entropy of decorated graph limits}

\author{Victor Falgas-Ravry \\
\small Ume{\aa} Universitet\\
\small \tt victor.falgas-ravry@umu.se\\
\and
Kelly O'Connell \\
\small Vanderbilt University\\ 
\small \tt kelly.m.oconnell@vanderbilt.edu
\and
Johanna Str\"omberg\\
\small Uppsala Universitet\\ 
\small \tt johanna.stromberg@math.uu.se
\and 
Andrew Uzzell\\
\small University of Nebraska--Lincoln\\ 
\small \tt andrew.uzzell@unl.edu}

\begin{document}
\maketitle
	\begin{abstract}
		In recent breakthrough results, Saxton--Thomason and Balogh--Morris--Samotij have developed powerful theories of hypergraph containers. These theories have led to a large number of new results on transference, and on counting and characterising typical graphs in hereditary properties. In a different direction, Hatami--Janson--Szegedy proved results on the entropy of graph limits which enable us to count and characterise graphs in dense hereditary properties.

\noindent In this paper, we make a threefold contribution to this area of research:
\begin{enumerate}[(i)]
\item We generalise results of Saxton--Thomason to obtain container theorems for general, dense hereditary properties of multicoloured graphs. Our main tool is the adoption of an entropy-based framework. As corollaries, we obtain general counting, characterization and transference results. We further give a streamlined extension of our results to cover a great variety of combinatorial structures: directed graphs, oriented graphs, tournaments, multipartite graphs, multi-graphs, hypercubes and hypergraphs.
\item We generalise the results of Hatami--Janson--Szegedy on the entropy of graph limits to the setting of decorated graph limits. In particular we define a cut norm for decorated graph limits and prove compactness of the space of decorated graph limits under that norm.
\item We explore a weak equivalence between the container and graph limit approaches to counting and characterising graphs in hereditary properties. In one direction, we show how our multicolour containers may be used to fully recover decorated versions of the results of Hatami--Janson--Szegedy. In the other direction, we show that our decorated extensions of Hatami--Janson--Szegedy's results on graph limits imply
counting and characterization
applications.
\end{enumerate}
Finally, we raise the problem of determining the possible structure of entropy maximisers in a multicoloured setting, and discuss the contrasts between the container and the graph limit approaches to counting.

Similar container results were recently obtained independently by Terry.

	\end{abstract}

{
  \hypersetup{linkcolor=black}
  \tableofcontents
}
	\section{Introduction}
	\subsection{Notation and basic definitions}
	Given a natural number~$r$, we write $A^{(r)}$ for the collection of all subsets of $A$ of size~$r$. We denote the powerset of~$A$ by $\{0,1\}^{A}$.
An \emph{$r$-uniform hypergraph}, or \emph{$r$-graph}, is a pair~$G=(V,E)$, where $V=V(G)$ is a set of \emph{vertices} and $E=E(G)\subseteq V^{(r)}$ is a set of \emph{$r$-edges}. We shall usually write `graph' for `$2$-graph' and, when there is no risk of confusion, `edge' for `$r$-edge'. We denote by $e(G):=\vert E(G)\vert$ the \emph{size} of~$G$ and by $v(G):=\vert V(G)\vert$ its \emph{order}.

	A \emph{subgraph} of an $r$-graph~$G$ is an $r$-graph~$H$ with $V(H)\subseteq V(G)$ and $E(H)\subseteq E(G)$. Given a set of vertices~$A\subseteq V(G)$, the subgraph of~$G$ \emph{induced} by $A$ is $G[A]:=(A, E(G)\cap A^{(r)})$. A set of vertices~$A$ is \emph{independent} in $G$ if the subgraph it induces contains no edges. The \emph{degree} of a set~$A\subseteq V(G)$ of size at most~$r-1$ is 
	\[\deg(A):=\vert \{f\in E(G)\, : \, A \subseteq f\}\vert. \] 
Finally an \emph{isomorphism} between $r$-graphs $G_1$ and~$G_2$ is a bijection~$\phi: \ V(G_1)\rightarrow V(G_2)$ which sends edges to edges and non-edges to non-edges.

Let $[n]:=\{1,2,\ldots, n\}$. A property $\mathcal{P}$ of (labelled) $r$-graphs is a sequence $(\mathcal{P}_n)_{n \in \N}$, where $\mathcal{P}_n$ is a collection of $r$-graphs on the \emph{labelled} vertex set $[n]$.  (Hereafter, we shall not distinguish between a property~$\PP$ and the class of $r$-graphs with $\PP$.)  An $r$-graph property is \emph{symmetric} if it is closed under relabelling of the vertices, i.e.\ under permutations of the vertex set $[n]$.  An $r$-graph property is \emph{monotone (decreasing)} if for every $r$-graph~$G\in \mathcal{P}$, every subgraph~$H$ of~$G$ is isomorphic to an element of~$\mathcal{P}$. A symmetric $r$-graph property is \emph{hereditary} if for every $r$-graph~$G\in \mathcal{P}$ every \emph{induced} subgraph~$H$ of~$G$ is isomorphic to an element of $\mathcal{P}$. Note that every monotone property is hereditary, but that the converse is not true. For example, the property of not containing a $4$-cycle as an induced subgraph is hereditary but not monotone.

In order to encode certain combinatorial objects of interest, such as directed graphs, we will consider a weaker notion of symmetry for hereditary graph properties.
\begin{definition}[Order-hereditary]
Let $m$,~$n\in\N$ with $m\leq n$. An \emph{order-preserving} map from $[m]$ to $[n]$ is a function $\phi: \ [m]\rightarrow [n]$ such that $\phi(i)\leq\phi(j)$ whenever $i\leq j$. Given graphs $G_1$ on $[m]$ and~$G_2$ on $[n]$, we say that $G_2$ contains $G_1$ as an \emph{order-isomorphic subgraph} if there is an order-preserving isomorphism from $G_1$ to an $m$-vertex subgraph~$H$ of~$G_2$. We further say that $G_2$ contains $G_1$ as an \emph{induced} order-isomorphic subgraph if the $m$-vertex subgraph~$H$ in question is an induced subgraph of~$G_2$.

A graph property $\mathcal{P}$ is said to be \emph{order-hereditary} if for every $G\in \mathcal{P}_n$ and every order-preserving injection $\phi: \ [m]\rightarrow [n]$, the graph~$G'=([m], \{f: \ \phi(f)\in E(G)\})$ is a member of $\mathcal{P}_m$. 	
\end{definition}
Clearly, every symmetric hereditary property is order-hereditary, but the converse is not true. As an example, consider the property~$\mathcal{P}$ of not containing an increasing path of length $2$, that is, the collection of graphs on $[n]$ ($n\in\N$) not containing vertices $i<j<k$ such that $ij$ and $jk$ are both edges. This is order-hereditary, but not symmetric --- and, as we shall see in Section~\ref{section: examples}, is much larger than the symmetric monotone property of not containing a path of length~$2$.

Finally, we shall use standard Landau notation throughout this paper, which we recall here. Given functions $f$,~$g:\N\to\mathbb{R}$, we have $f=O(g)$ if there exists a constant $C>0$ such that $\limsup_{n\to\infty}{f(n)}/{g(n)}\leq C$. If $\lim_{n\to\infty}{f(n)}/{g(n)}=0$, then we write $f=o(g)$. We write $f=\Omega (g)$ and $f=\omega(g)$ to denote $g=O(f)$ and $g=o(f)$ respectively. If we have both $f=O(g)$ and $f=\Omega(g)$, we say that $f$ and $g$ are of the same order and denote this by $f=\theta(g)$. We shall sometimes use $f\ll g$ and $f \gg g$ as alternatives to $f=o(g)$ and $f=\omega(g)$, respectively. Finally, we say that a sequence of events~$A_n$ occurs \emph{with high probability (whp)} if $\lim_{n\rightarrow \infty} \mathbb{P}(A_n)=1$.

	\subsection{Background: hereditary graph properties and their speeds}\label{section: background hereditary props}
The problem of counting and characterising graphs in a given symmetric hereditary property $\mathcal{P}$ has a long and distinguished history. The \emph{speed} $n\mapsto \vert \mathcal{P}_n \vert$ of a graph property was introduced in 1976 by Erd{\H o}s, Kleitman and Rothschild~\cite{ErdosKleitmanRothschild76}. Together with the structural properties of a `typical' element of $\mathcal{P}_n$, it has received extensive attention from the research community.

Early work focussed on the case where $\mathcal{P}=\Forb(F)$, the monotone decreasing property of not containing a fixed graph $F$ as a subgraph. We refer to the graphs in $\Forb(F)$ as \emph{$F$-free} graphs. The \emph{Tur\'an number} of~$F$, denoted by $\ex(n, F)$, is the maximum number of edges in an $F$-free graph on $n$ vertices. Clearly, any subgraph of an $F$-free graph is also $F$-free. This gives the following lower bound on the number of $F$-free graphs on $n$ labelled vertices:
\[\Forb(F)_n \geq 2^{\ex(n, F)}.\]
Erd{\H o}s, Kleitman and R{\"o}dl~\cite{ErdosKleitmanRothschild76} showed that if $F=K_t$, the complete graph on $t$ vertices, then the exponent in this lower bound is asymptotically tight:
\[\Forb(K_t)_n \leq 2^{\bigl(1+o(1)\bigr)\ex(n, K_t)}.\]
Their work was generalised by Erd{\H o}s, Frankl and R{\"o}dl~\cite{ErdosFranklRodl86} to the case of arbitrary forbidden subgraphs~$F$ and by Pr\"omel and Steger~\cite{PS92}, who considered the property $\Forb^{*}(F)$ of not containing $F$ as an induced subgraph.  Finally, Alekseev~\cite{Alekseev93} and Bollob\'as--Thomason~\cite{BT} independently determined the asymptotics of the logarithm of the speed for any symmetric hereditary property in terms of its \emph{colouring number}, which we now define.
\begin{definition}\label{definition: colouring number}
For each $r\in \N$ and $\vector{v}\in \{0,1\}^r$, let  $\mathcal{H}(r,\vector{v})$ be the collection of all graphs $G$ such that $V(G)$ may be partitioned into $r$ disjoint sets $A_1$, \dots,~$A_r$ such that for each $i$, $G[A_i]$ is an empty graph if $v_i=0$ and a complete graph if $v_i=1$. The \emph{colouring number} $\chi_c(\mathcal{P})$ of a symmetric hereditary property is defined to be
\[\chi_c(\mathcal{P}):=\sup\bigl\{r\in \N \, : \, \mathcal{H}(r,\vector{v})\subseteq \mathcal{P} \textrm{ for some }\vector{v} \in \{0,1\}^r\bigr\}.\]
\end{definition}
\begin{theorem}[Alekseev--Bollob\'as--Thomason Theorem]\label{theorem: alekseevbollobasthomason}
If\/ $\mathcal{P}$ is a symmetric hereditary property of graphs with $\chi_c(\mathcal{P})=r$, then
\[\lim_{n\rightarrow \infty} \dfrac{\log_2\vert\mathcal{P}_n\vert}{\binom{n}{2}}= 1-\dfrac{1}{r}.\]
\end{theorem}
Subsequently, the rate of convergence of $\log_2 \vert \mathcal{P}_n\vert/\binom{n}{2}$ and the structure of typical graphs were investigated by Balogh, Bollob\'as and Simonovits~\cite{BBS04, BBS09} for symmetric monotone properties, and by Alon, Balogh, Bollob\'as and Morris~\cite{ABBM11} for symmetric hereditary properties.

There has also been interest in the speed of monotone properties in other discrete structures. Kohayakawa, Nagle and R\"odl~\cite{KohayakawaNagleRodl03}, Ishigami~\cite{Ishigami07}, Dotson and Nagle~\cite{DotsonNagle09} and Nagle, R\"odl and Schacht~\cite{NagleRodlSchacht06} investigated the speed of hypergraph properties, while in a series of papers Balogh, Bollob\'as and Morris~\cite{BaloghBollobasMorris06, BaloghBollobasMorris07b, BaloghBollobasMorris07} studied the speed of properties of ordered graphs, oriented graphs and tournaments. Many of these results relied on the use of graph and hypergraph regularity lemmas. See the survey of Bollob\'as~\cite{Bollobas07} for an overview of the state of the area \emph{before} the breakthroughs discussed in the next subsection.
\subsection{Background: transference and containers}
Recently, there has been great interest in \emph{transference} theorems, in which central results of extremal combinatorics are shown to also hold in `sparse random' settings.   These results are motivated by, \emph{inter alia}, the celebrated Green--Tao theorem on arithmetic progressions in the primes~\cite{GreenTao08} and the K{\L}R conjecture of Kohayawa, {\L}uczak and R\"odl~\cite{KLR} and its applications. (Very roughly, the K{\L}R conjecture says that, given a graph~$H$ and $p = p(n) \in \oi$ large enough, with high probability, every subgraph of an Erd{\H{o}}s--R{\'{e}}nyi random graph~$G(n, p)$ has approximately the `right' number of copies of~$H$.  See~\cite{ConlonGowersSamotijSchacht14} for a discussion of the conjecture and its applications.)

In major breakthroughs a little over five years ago, Conlon and Gowers~\cite{ConlonGowers10} and independently Friedgut, R\"odl and Schacht~\cite{FriedgutRodlSchacht10} and Schacht~\cite{Schacht2009} proved very general transference results, which in particular settled many cases of the K{\L}R conjecture. Their work was soon followed by another dramatic breakthrough: Balogh, Morris and Samotij~\cite{BaloghMorrisSamotij15} and independently Saxton and Thomason~\cite{SaxtonThomason15}, building on work of Kleitman--Winston~\cite{KleitmanWinston82} and of Sapozhenko~\cite{Sapozhenko87,Sapozhenko01} for graphs, developed powerful theories of hypergraph containers.

These container theories essentially say that hereditary properties can be `covered' by `small' families of `containers', which are themselves `almost in the property'. We discuss containers with more precision and details in Section~\ref{section: containers}. As an application of their theories, Balogh--Morris--Samotij and Saxton--Thomason gave both new proofs of known counting/characterization results and many new counting/characterization results for hereditary properties, and in addition a spate of transference results.
In particular Balogh--Morris--Samotij and Saxton--Thomason settled the K{\L}R conjecture in full generality. 

We refer the reader to the excellent ICM survey of Conlon~\cite{Conlon14} for an in-depth discussion of the recent groundbreaking progress made by researchers in the area.
	\subsection{Background: entropy and graph limits}
A parallel but separate development at the intersection of extremal combinatorics and discrete probability has been the rise of theories of limit objects for sequences of discrete structures. The first to appear was the theory of exchangeable random variables, originating in the work of de~Finetti in the 1930s and further developed in the 1980s by Aldous, Hoover and Kallenberg amongst others, see the monograph of Aldous~\cite{Aldous} on the subject. After the turn of the century, two more approaches to limit objects from a more combinatorial perspective garnered attention. First of all, the study of left convergence/dense graph limits was initiated by Borgs, Chayes, Lov\'asz, S\'os and Vesztergombi~\cite{BCLSV08} and by Lov\'asz and Szegedy~\cite{LS06}.  (For a thorough development of the accompanying theory, see the monograph of Lov\'asz~\cite{LovaszBook}.)  In a different direction, Razborov~\cite{Razborov07} developed flag algebras with a view to applications to extremal combinatorics; see also~\cite{Razborov13} for an introduction to Razborov's theory and its ramifications. We refer the reader to Austin~\cite{Austin08} for an exposition and thorough analysis of the links between these three limit object theories.

In this paper, we focus on the theory of graph limits.  We shall give precise definitions later, but for now, it is enough to say that certain sequences of graphs are defined to be \emph{convergent}. If $(G_n)_{n=1}^{\infty}$ is a convergent graph sequence, then its limit can be represented by a \emph{graphon}, i.e., a symmetric, measurable function~$W : \oi^2 \to \oi$.  A deep result of Lov\'asz and Szegedy~\cite{LS:analyst} says that the set of graphons forms a compact topological space with respect to a certain `cut metric'.

Recently, Hatami, Janson and Szegedy~\cite{HJS} defined and studied the entropy of a graphon.  They used this notion to recover Theorem~\ref{theorem: alekseevbollobasthomason} and to describe the typical structure of a graph in a hereditary property.
The Hatami--Janson--Szegedy notion of entropy can be viewed as a graphon analogue of the classical notion of the entropy of a discrete random variable, which first appeared in Shannon's foundational paper~\cite{Shannon48}. Using entropy to count objects is an old and celebrated technique in discrete probability --- see for example Galvin~\cite{Galvin14} for an exposition of the applications of entropy to counting.

	\subsection{Contributions of this paper}
Our first contribution in this paper is to prove very general multicolour hypergraph container statements. Amongst other structures of interest, our results cover directed graphs, oriented graphs, tournaments, multipartite graphs, square grids and both edge- and vertex-subgraphs of hypercubes. We use our results to obtain general counting, characterization and transference results for hereditary properties of the aforementioned structures. As we restrict ourselves to the study of `dense' properties, our container statements and their corollaries are (we believe) simple and easy to apply (albeit weaker than the full strength of the Balogh--Morris--Samotij and Saxton--Thomason container theorems), which we hope may be useful to other researchers.

Our main tool in this part of the paper is a container theorem of Saxton--Thomason for linear hypergraphs together with the adoption of an entropy-based framework. We should like to emphasize here the intellectual debt this paper owes to the pioneering work of Balogh--Morris--Samotij and Saxton--Thomason: our work relies on theirs in a crucial way, and some of our ideas exist already in their papers in an embryonic form, which we explore further. The usefulness of our exploration is vindicated by the fact that some of the applications of containers to other discrete structures which we treat are new, and were not well understood by the mathematical community at the time of writing. For example, finding a container theorem for digraphs was a problem raised by K\"uhn, Osthus, Townsend and Zhao~\cite{KuhnOsthusTownsendZhao14}, which we resolve in the present paper.

Our second main contribution is to relate container theorem to the work of Hatami--Janson--Szegedy on the entropy of graph limits.  Given a set~$K$, a \emph{$K$-decorated graph of order~$n$} is a labelling of~$E(K_n)$ with elements of~$K$.  (An ordinary graph may be viewed as a $\{0, 1\}$-decorated graph with edges labelled~$1$ and non-edges labelled~$0$.)  We use our multicolour container theorems to obtain generalizations of Hatami--Janson--Szegedy's results to the setting of decorated graph limits. In the other direction, we obtain a second proof of these generalizations by working directly in the world of decorated graphons (which requires us to construct a cut metric for decorated graphons and show compactness of the space under that metric, amongst other things). We then show how these analytic results can be used to recover many of the main combinatorial applications of containers, namely counting and characterization
for hereditary properties of multicoloured graphs, and contrast the container and entropy of graph limit approaches to counting. This is part of an attempt to build links between the rich and currently quite distinct theories of graph limits and hypergraph containers.

We note that there is a significant overlap between the container, enumeration, and stability results presented here and those obtained independently by Terry~\cite{Terry16}, although the emphases of our papers are rather different.
		
	\subsection{Structure of the paper}
Section~\ref{section: containers} gathers together our main results on multicolour containers. Section~\ref{subsection: key definitions templates and entropy} contains our key definitions of templates and entropy.  In Section~\ref{subsection: containers proof}, we state and prove our first multicolour container theorem (Theorem~\ref{theorem: multi-colour container}), and in Section~\ref{subsection: entropy density, supersaturation} we introduce entropy density and prove a supersaturation result that is key to several of our applications.  In Section~\ref{subsection: approximation of arbitrary hereditary properties}, we use these tools to prove container theorems for general hereditary properties (Corollary~\ref{corollary: containers for arbitrary hereditary properties}) and prove a general counting result (Corollary~\ref{corollary: speed of arbitrary hereditary properties}).  Finally in Sections \ref{subsection: stability} and~\ref{subsection: transference} we obtain general characterization and transference results (Theorems \ref{theorem: strong stability and containers} and~\ref{theorem: transference}, respectively).  As indicated earlier, the results of Sections \ref{subsection: containers proof}--\ref{subsection: stability} are very similar to those proved by Terry~\cite{Terry16}.  In particular, our Terry's Theorems 2, 3, 6 and~7 correspond to our Proposition~\ref{proposition: entropy density}, Corollary~\ref{corollary: speed of arbitrary hereditary properties}, Theorem~\ref{theorem: multi-colour container} and Lemma~\ref{lemma: supersaturation}, respectively, while her Theorem~5 is very similar to our Theorem~\ref{theorem: strong stability and containers}.  Furthermore, Terry's results hold for uniform hypergraphs, and so do ours, as shown in Section~\ref{subsection: hypergraphs}.

In Section~\ref{section: other structures}, we extend our main results to a number of other discrete structures. Section~\ref{subsection: oriented} describes how our theorems apply to oriented and directed graphs; as mentioned earlier, this addresses an issue raised in~\cite{KuhnOsthusTownsendZhao14}. In Section~\ref{subsection: other host graphs} we extend our main results to cover colourings of sequences of graphs (rather than sequences of complete graphs). Oour results in that subsection cover a multitude of examples including grid graphs, multipartite graphs and hypercube graphs. In Section~\ref{subsection: hypergraphs} we extend our results to a general hypergraph setting, which allows us amongst other things to prove in Section~\ref{subsection: vertex colouring hypercube} results on vertex-colourings of hypercubes.

Section~\ref{section: examples} is dedicated to applications of our results to a variety of examples (graphs, digraphs, multigraphs, multicoloured graphs and hypercubes).  In particular, we give a new, short proof of the Alekseev--Bollob\'as--Thomason theorem and prove counting and characterization results for hereditary properties of directed graphs.

In the next part of our paper, we turn to graph limits.  In Section~\ref{section: cut metric for graphons}, we define a cut norm for decorated graphons and prove compactness of the space of decorated graphons under that norm. This continues a program of Lov\'asz and Szegedy~\cite{LS:decorated}.  Section~\ref{section: containers to graphon entropy} shows how we may use our container and compactness results to obtain generalizations of results of Hatami, Janson and Szegedy~\cite{HJS} to decorated graphons.   Finally, in Section~\ref{section: decorated generalization of HJS}, we use the results of Section~\ref{section: cut metric for graphons} to give a second proof of the generalizations
of the results of Hatami--Janson--Szegedy on the entropy of graph limits.

We end this long paper in Section~\ref{section: concluding remarks} with an open problem on the possible structure of entropy maximizers in the multicolour/decorated setting and a discussion of the differences between the container and the graph limit approaches we have explored.
	\section{Multi-colour containers}\label{section: containers}
\subsection{Key definitions: templates and entropy}\label{subsection: key definitions templates and entropy}
Let $K_n$ denote the complete graph $([n], [n]^{(2)})$. We study $k$-colourings of (the edges of) $K_n$, that is to say, we work with the set of colouring functions~$c: \ E(K_n)\rightarrow [k]$. Denote by $[k]^{K_n}$ the set of all such colourings. Note that each colour~$i$ induces a graph~$c_i$ on $[n]$, $c_i=([n], c^{-1}(i))$. An ordinary graph~$G$ may be viewed as a $2$-colouring of~$E(K_n)$, with $G=c_1$ and $\overline{G}=c_2$. An oriented graph~$\vec{G}$ may be viewed as a $3$-colouring of~$E(K_n)$, such that each edge~$ij$ with $i<j$ is coloured $1$ if $\vec{ij}\in D$, $2$ if $\vec{ji} \in D$ and $3$ otherwise.

Our notions of subgraph and isomorphism carry over to the $k$-colouring setting in the natural way: two $k$-colourings $c$ and $c'$ of~$K_n$ are \emph{isomorphic} if there is a bijection $\phi:\ [n] \rightarrow [n]$ such that $\phi$ is an isomorphism from $c_i$ to $c'_i$ for each colour $i\in [k]$. Given $m\leq n$ and $k$-colourings $c$,~$c'$ of~$K_m$~and~$K_n$ respectively, we say that $c$ is an \emph{(order-preserving) subcolouring} of $c'$ if there is an order-preserving injection~$\phi: \ [m]\rightarrow [n]$ such that $c(ij)=c'(\phi(i)\phi(j))$ for all $1\leq i<j\leq m$. Finally, given an $m$-set~$A=\{a_1, a_2, \ldots, a_m\}\subseteq [n]$ with $a_1<a_2<\cdots <a_m$ and a $k$-colouring~$c$ of~$K_n$, the \emph{(order-preserving) restriction} of~$c$ to~$A$ is the $k$-colouring~$c_{\vert A}$ of~$K_m$ defined by $c_{\vert A}(ij)=c(a_ia_j)$. Thus $c'$ is a subcolouring of~$c$ if and only if there exists a set~$A$ such that $c'=c_{\vert A}$.

Our main object of study in this section will be \emph{order-hereditary} properties of~$[k]^{K_n}$.
\begin{definition}[Order-hereditary property]\label{definition: order hereditary}
An \emph{order-hereditary property} of $k$-colourings is a sequence~$\mathcal{P}= \left(\mathcal{P}_n\right)_{n\in \N}$, such that:
\begin{enumerate}[(i)]
	\item $\mathcal{P}_n$ is a family of $k$-colourings of~$K_n$,
	\item for every $c\in \mathcal{P}_n$ and every $A\subseteq [n]$, $c_{\vert A}\in \mathcal{P}_{\vert A\vert}$.
\end{enumerate} 		
\end{definition}
A key tool in extending container theory to $k$-coloured graphs will be the following notion of a \emph{template}:
\begin{definition}[Template]\label{definition: template}
A \emph{template} for a $k$-colouring of~$K_n$ is a function $t: \ E(K_n)\rightarrow	\{0,1\}^{[k]}$, associating to each edge $f$ of~$K_n$ a non-empty list of colours $t(f)\subseteq [k]$; we refer to this set $t(f)$ as the \emph{palette} available at $f$.


Given a template~$t$, we say that a $k$-colouring~$c$ of~$K_n$ \emph{realises} $t$ if $c(f)\in t(f)$ for every edge~$f\in E(K_n)$.  We write $\langle t \rangle$ for the collection of realisations of~$t$.
\end{definition}
 In other words, a template~$t$ gives for each edge of~$K_n$ a palette of permitted colours, and $\langle t \rangle$ is the set of $k$-colourings of~$K_n$ that respect the template. We observe that a $k$-colouring of~$K_n$ may itself be regarded as a template, albeit with only~one colour allowed at each edge. We extend our notion of subcolouring to templates in the natural way.
 \begin{definition}[Subtemplate]\label{definition: subtemplate}
 The \emph{(order-preserving) restriction} of a $k$-colouring template~$t$ for $K_n$ to an $m$-set~$A=\{a_1,a_2,\ldots, a_m\}\subseteq [n]$ with $1\leq a_1<a_2<\cdots a_m\leq m$ is the $k$-colouring template for $K_m$ $t_{\vert A}$ defined by $t_{\vert A}(ij)=t(a_ia_j)$.

 Given $m\leq n$ and  $k$-colouring templates $t$,~$t'$ for $K_m$,~$K_n$ respectively, we say that $t$ is a \emph{subtemplate} of~$t'$, which we denote by $t\leq t'$, if there exists an $m$-set~$A\subseteq [m]$ such that $t(f)\subseteq t'_{\vert A}(f)$ for each $f\in E(K_m)$. Furthermore $t$ is an \emph{induced} subtemplate of~$t'$ if $t=t'_{\vert A}$ for some $m$-set~$A\subseteq [n]$.
\end{definition}
Our notion of subtemplates can be viewed as the template analogue of the notion of an order-preserving subgraph for graphs on a linearly ordered vertex set.

 Templates enable us to generalise the notion of containment to the $k$-coloured setting.
 \begin{definition}[Container family]
	Given a family of $k$-colourings $\mathcal{F}$ of $E(K_n)$, a \emph{container family} for $\mathcal{F}$ is a collection $\mathcal{C}=\{t_1, t_2, \ldots, t_m\}$ of $k$-colouring templates such that $\mathcal{F}\subseteq \bigcup_i \langle t_i \rangle$.  (In other words, every colouring in $\mathcal{F}$ is a realisation of some template in $\mathcal{C}$.)
\end{definition}
Next we introduce the key notion of the \emph{entropy} of a template.
\begin{definition}\label{def: entropy}
The \emph{entropy} of a $k$-colouring template $t$ is 
\[\Ent(t):= \log_k\prod_{e\in E(K_n)} \vert t(e)\vert.\]
\end{definition}
Observe that for any template~$t$, $0\leq \Ent(t)\leq \binom{n}{2}$, and that the number of distinct realisations of~$t$ is exactly~$\vert\langle t \rangle\vert= k^{\Ent(t)}$. There is a direct correspondence between our notion of entropy and that of Shannon entropy in discrete probability: given a template $t$, we can define a \emph{$t$-random} colouring~$\mathbf{c}_t$, by choosing for each $f\in E(K_n)$ a colour~$\mathbf{c}_t(f)$ uniformly at random from $t(f)$. The entropy of~$t$ as defined above is exactly the $k$-ary Shannon entropy of the discrete random variable~$\mathbf{c}$. Finally, observe that zero-entropy templates correspond to $k$-colourings of~$K_n$ and that if $t$ is a subtemplate of~$t'$ then $\Ent(t)\leq \Ent(t')$.

\subsection{Containers}\label{subsection: containers proof}

Let $N\in \N$ be fixed and let $\mathcal{F}=\{c_1, c_2, \ldots, c_m\}$ be a nonempty collection of $k$-colourings of~$E(K_N)$. Let $\Forb(\mathcal{F})$ be the collection of all $k$-colourings~$c$ of~$K_n$, $n\in \N$, such that $c_i\not\leq c$ for $i: 1\leq i\leq m$. 
More succinctly, $\Forb(\mathcal{F})$ is the collection of all $k$-colourings avoiding $\mathcal{F}$, which clearly is an order-hereditary property of $k$-colourings.

\begin{theorem}\label{theorem: multi-colour container}
Let $N\in \N$ be fixed and let $\mathcal{F}$ be a nonempty collection of $k$-colourings of~$E(K_N)$.  For any $\varepsilon>0$, there exist constants $C_0$,~$n_0>0$, depending only on $(\varepsilon, k, N)$,  such that for any $n\geq n_0$ there exists a collection~$\mathcal{T}$ of $k$-colouring templates for $K_n$ satisfying:
\begin{enumerate}[(i)]
\item $\mathcal{T}$ is a container family for $\left(\Forb(\mathcal{F})\right)_n$;
\item for each template $t\in\mathcal{T}$, there are at most $\varepsilon \binom{n}{N}$ sets $A$ such that $c_i\leq t_{\vert A}$ for some $c_i\in \mathcal{F}$;
\item $\log_k \vert \mathcal{T}\vert\leq C_0n^{-1/(2\binom{N}{2}-1)} \binom{n}{2}$.
\end{enumerate}
\end{theorem}
In other words, the theorem says that we can find a \emph{small} (property (iii)) collection of templates, that together \emph{cover} $\Forb(\mathcal{F})_n$ (property (i)), and whose realisations are \emph{close} to lying in $\Forb(\mathcal{F})_n$ (property (ii)).

We shall deduce Theorem~\ref{theorem: multi-colour container} from a hypergraph container theorem of Saxton and Thomason. Say that an $r$-graph~$H$ is \emph{linear} if each pair of distinct $r$-edges of~$H$ meets in at most~$1$ vertex. Saxton and Thomason proved the following:
\begin{theorem}[Saxton--Thomason (Theorem 1.2 in~\cite{SaxtonThomason16})]\label{theorem: saxton thomason simple containers}
	Let $0<\delta<1$. Then there exists $d_0=d_0(r, \delta)$ such that if $G$ is a linear $r$-graph of average degree~$d\geq d_0$ then there exists a collection~$\mathcal{C}$ of subsets of~$V(G)$ satisfying:
	\begin{enumerate}
		\item if $I\subseteq V(G)$ is an independent set, then there exists $C\in \mathcal{C}$ with $I\subseteq C$;
		\item $e(G[C])< \delta e(G)$ for every $C\in \mathcal{C}$;
		\item  $\vert \mathcal{C} \vert \leq 2^{\beta v(G)}$, where $\beta=(1/d)^{1/(2r-1)}$.
	\end{enumerate}  	
\end{theorem}

In the proof of Theorem~\ref{theorem: multi-colour container} and elsewhere, we shall use the following standard Chernoff bound: if $X\sim\Binom(n,p)$, then for any $\delta\in[0,1]$,
\begin{equation}\label{equation: Chernoff} \mathbb{P}\bigl(\vert X-np\vert \geq \delta np\bigr) \leq 2e^{-\frac{\delta^2np}{4}}.\end{equation}

\begin{proof}[Proof of Theorem~\ref{theorem: multi-colour container}]
If $N=2$, then $\mathcal{F}$ just gives us a list of forbidden colours, say $F\subseteq [k]$. Then $\left(\Forb(\mathcal{F})\right)_n$ is exactly the collection of all realisations of the template~$t$ assigning to each edge~$e$ of~$K_n$ the collection~$[k]\setminus F$ of colours not forbidden by $\mathcal{F}$. Thus in this case our result trivially holds, and we may therefore assume $N\ge3$.

We define a hypergraph~$H$ from $\mathcal{F}$ and $K_n$ as follows. Set $r=\binom{N}{2}$. We let the vertex set of $H$ consist of~$k$ disjoint copies of~$E(K_n)$, one for each of our $k$ colours: $V(H)=E(K_n)\times [k]$; this key idea, allowing us to apply Theorem~\ref{theorem: saxton thomason simple containers}, first appeared (as far as we know) in a $2$-colour form in the seminal paper of Saxton and Thomason~\cite{SaxtonThomason15}. For every $N$-set~$A\subseteq [n]$ and every $k$-colouring~$c$ of the edges of the (order-preserving) copy of~$K_N$ induced by $A$ with $c\in \mathcal{F}$, we add to $H$ an $r$-edge~$e_{c,A}$, where
\[e_{c, A}=\bigl\{\bigl(e, c(e)\bigr) \, : \, e\in A^{(2)} \bigr\}.\]
This gives us an $r$-graph~$H$. Let us give bounds on its average degree. 
Since $\mathcal{F}$ is nonempty, for every $N$-set $A\subseteq [n]$, there are at least~$1$ and at most~$k^{\binom{N}{2}}$ colourings~$c$ of $A^{(2)}$ which are order-isomorphic to an element of~$\mathcal{F}$. Thus 
\begin{equation}\label{equation: bounds on e(H)}
\frac{n^N}{N^N} \leq  \binom{n}{N} \leq e(H)\leq k^{\binom{n}{2}}\binom{n}{N}\leq k^{\binom{N}{2}}\left(\frac{en}{N}\right)^N .
\end{equation}
Thus $e(H)$ is of order~$n^N$ and the average degree in $H$ is of order~$n^{N-2}$, which tends to infinity as $n\rightarrow \infty$. We are almost in a position to apply Theorem~\ref{theorem: saxton thomason simple containers}, with one significant caveat: the hypergraph~$H$ we have defined is in no way linear. Following~\cite{SaxtonThomason16}, we circumvent this difficulty by considering a \emph{random sparsification} of~$H$.

Let $\varepsilon_1 \in(0,1)$ be a constant to be specified later and let
\begin{equation}\label{eq: sparsification probability}
p=\varepsilon_1 \big/\left(24k^{2\binom{N}{2}-3}\binom{N}{3}\binom{n-3}{N-3}\right).
\end{equation}
We shall keep each $r$-edge of~$H$ independently with probability~$p$, and delete it otherwise, to obtain a random subgraph~$H'$ of~$H$. Standard probabilistic estimates will then show that with positive probability the $r$-graph~$H'$ is almost linear, has large average degree and respects the density of~$H$. More precisely, we show:
\begin{lemma}\label{lemma: random sparsification}
Let $p$ be as in~\eqref{eq: sparsification probability}, let $H'$ be the random subgraph of~$H$ defined above and consider the following events:
\begin{itemize}
	\item $F_1$ is the event that $e(H')\geq \frac{pe(H)}{2}\geq \frac{p\binom{n}{N}}{2}$;
	\item $F_2$ is the event that $H'$ has at most $3p^2k^{2\binom{N}{2}-3}\binom{n}{N}\binom{N}{3}\binom{n-3}{N-3}= \frac{\varepsilon_1}{8}p\binom{n}{N}$ pairs of edges~$(f,f')$ with $\vert f\cap f'\vert \geq 2$;
	\item $F_3$ is the event that for all $S\subseteq V(H)$ with $e(H[S])\geq \varepsilon_1 e(H)$, we have $e(H'[S])\geq \frac{\varepsilon_1}{2} e(H')$.
\end{itemize}
There exists $n_1=n_1(\varepsilon_1, k, N)\in\N$ such that for all $n\geq n_1$, $F_1 \cap F_2 \cap F_3$ occurs with strictly positive probability.
\end{lemma}
\begin{proof}
By (\ref{equation: bounds on e(H)}), we have $\mathbb{E}e(H')=pe(H)\geq p\binom{n}{N}$. Applying the Chernoff bound (\ref{equation: Chernoff}) with $\delta=1/2$, we get that the probability that $F_1$ fails in $H'$ is at most
\begin{equation*}\label{equation: bound on prob 1. fails}
\mathbb{P}\biggl(e(H')\leq \frac{n^3}{2kN^N\log n}\biggr)\leq 2e^{-\frac{p\binom{n}{N}}{16}}= e^{-\Omega(n^3)}.
\end{equation*}

Next consider the pairs of $r$-edges $(f,f')$ in $H$ with $\vert f\cap f'\vert \geq 2$, which we refer to hereafter as \emph{overlapping pairs}.  Let $Y_H$ denote the number of overlapping pairs in $H$ and define $Y_{H'}$ similarly. Note that $Y_H$ is certainly bounded above by the number of ways of choosing an $N$-set~$A$, a $3$-set~$B$ from $A$ and an $(N-3)$-set~$A'$ from $[n]\setminus B$ (thereby making an overlapping pair of $N$-sets~$(A, A'\cup B)$) and assigning an arbitrary $k$-colouring to the edges in $A^{(2)}\cup {A'}^{(2)}$. Thus,
\begin{equation*} Y_H \leq \binom{n}{N}\binom{N}{3} \binom{n-3}{N-3}k^{2\binom{N}{2}-3}
\end{equation*}
and
\begin{equation*}
\mathbb{E}(Y_{H'}) = p^2 Y_H \leq p^2\binom{n}{N}\binom{N}{3} \binom{n-3}{N-3}k^{2\binom{N}{2}-3} =\frac{\varepsilon_1}{24} p\binom{n}{N}.
\end{equation*} 
Applying Markov's inequality, we have that with probability at least~$\frac{2}{3}$, $Y_{H'} \leq \frac{\varepsilon_1}{8} p\binom{n}{N}$ and $F_2$ holds.

Finally, consider a set $S\subseteq V(H)$ with $e(H[S])\geq \varepsilon_1 e(H)$. Applying the Chernoff bound~\eqref{equation: Chernoff} and our lower bound~\eqref{equation: bounds on e(H)} on $e(H)$, we get
\begin{equation}\label{eq: too few induced edges}
\mathbb{P}\left(e\bigl(H'[S]\bigr)\leq\frac{1}{\sqrt {2}}\mathbb{E}e\bigl(H'[S]\bigr)\right) \leq 2e^{-\frac{\mathbb{E} e(H'[s])}{8}}=2e^{-\frac{pe(H[s])}{8}} \leq 2e^{-\frac{p\varepsilon_1 e(H)}{8}}= e^{-\Omega(n^3)}.
\end{equation}
Moreover, by \eqref{equation: Chernoff} and~\eqref{equation: bounds on e(H)} again,
\begin{equation}\label{eq: too many induced edges}
\mathbb{P}\left(e(H')\ge\sqrt{2}\mathbb{E} e(H')\right)\leq 2e^{-\frac{(\sqrt{2}-1)^2p e(H)}{4}}= e^{-\Omega(n^3)}.
\end{equation}
Say that a set~$S\subseteq V(H)$ is \emph{bad} if $e(H[S])\geq \varepsilon_1 e(H)$ and $e(H'[S])\leq \frac{\varepsilon_1}{2} e(H')$.   By \eqref{eq: too few induced edges},~\eqref{eq: too many induced edges} and the union bound, the probability that $F_3$ fails, i.e., that there exists some bad $S\subseteq V(H)$, is at most
\begin{equation*}
\mathbb{P}\bigl(\exists \textrm{ bad }S\bigr) \leq \mathbb{P}\bigl(e(H')\geq \sqrt{2}\mathbb{E} e(H')\bigr)+\sum_S \mathbb{P}\biggl(e(H(S)\leq\frac{1}{\sqrt {2}}\mathbb{E}e\bigl(H'[S]\bigr)\biggr)\leq 2^{k\binom{n}{2}} e^{-\Omega(n^3)}=e^{-\Omega(n^3)}.
\end{equation*}

Therefore with probability at least~$2/3-o(1)$ the events $F_1$, $F_2$ and~$F_3$ all occur, and in particular they must occur simultaneously with strictly positive probability for all $n\geq n_1=n_1(\varepsilon_1, k, N)$. 
\end{proof}
By Lemma~\ref{lemma: random sparsification}, for any $\varepsilon>0$ and $\varepsilon_1=k^{-\binom{N}{2}}\varepsilon$ fixed and any $n\geq n_1(\varepsilon_1)$, there exists a sparsification~$H'$ of~$H$ for which the events $F_1$, $F_2$ and~$F_3$ from the lemma all hold.  Deleting one $r$-edge from each overlapping pair in $H'$, we obtain a linear $r$-graph~$H''$ with average degree~$d$ satisfying
\begin{equation}\label{eq: linear hypergraph average degree}
d=\frac{e(H'')}{v(H'')} \geq \frac{1}{k\binom{n}{2}}\left(e(H') - Y_{H'}\right)
\geq \frac{1}{k\binom{n}{2}}\biggl(\frac{1}{2} -\frac{\varepsilon_1}{4}\biggr)p\binom{n}{N}=\Omega(n).
\end{equation}
We are now in a position to apply the container theorem for linear $r$-graphs, Theorem~\ref{theorem: saxton thomason simple containers}, to $H''$. Let $\delta=\delta(\varepsilon_1)$ satisfy $0<\delta <\varepsilon_1/4$ and let $d_0=d_0(\delta, r)$ be the constant in Theorem~\ref{theorem: saxton thomason simple containers}. For $n\geq n_2(k, N, \delta)$, we have $d\geq d_0$. Thus there exists a collection~$\mathcal{C}$ of subsets of~$V(H'')=V(H)$ satisfying conclusions 1.--3.\ of Theorem~\ref{theorem: saxton thomason simple containers}.

For each $C\in \mathcal{C}$, we obtain a template~$t=t(C)$ for a \emph{partial} $k$-colouring of~$K_n$, with the palette for each edge~$e$ given by $t(e)=\{i\in[k]: \ (e,i)\in C\}$ (note that some edges may have an empty palette). Let $\mathcal{T}$ be the collection of \emph{proper} templates obtained in this manner, that, is the collection of~$t=t(C)$ with $C\in \mathcal{C}$ and with each edge~$e\in E(K_n)$ having a nonempty palette~$t(e)$. We claim that the template family~$\mathcal{T}$ satisfies the conclusions (i)--(iii) of Theorem~\ref{theorem: multi-colour container} that we are trying to establish.

Indeed, by definition of~$H$, any colouring $c\in \mathcal{P}_n$ gives rise to an independent set in the $r$-graph~$H$ and hence its subgraph~$H''$, namely $I=\{(e,i): \ c(e)=i\}$. Thus there exist $C\in \mathcal{C}$ with $I\subseteq C$, giving rise to a proper template~$t \in \mathcal{T}$ with $c\in\langle t \rangle$. Conclusion~(i) is therefore satisfied by $\mathcal{T}$.

Further for each $C\in \mathcal{C}$, conclusion 2.\ of Theorem~\ref{theorem: saxton thomason simple containers} and the event~$F_2$ together imply
\begin{equation*}
e\bigl(H'[C]\bigr)\leq e\bigl(H''[C]\bigr) + \left(e(H')-e(H'')\right)\leq \delta e(H'') + \frac{\varepsilon_1}{4}e(H')< \frac{\varepsilon_1}{2}e(H').
\end{equation*}
Together with the fact that $F_3$ holds, this implies $e(H[C])<\varepsilon_1 e(H)$, which by (\ref{equation: bounds on e(H)}) is at most $\varepsilon_1 k^{\binom{N}{2}} \binom{n}{N} =\varepsilon \binom{n}{N}$. In particular, by construction of~$H$, we have that for each $t=t(C)\in \mathcal{T}$ there are at most~$\varepsilon \binom{n}{N}$ pairs~$(A, c_i)$ of $N$-sets~$A$ and forbidden colourings~$c_i\in \mathcal{F}$ with $c_i\leq t_{\vert A}$. This establishes (ii).

Finally by property 3.\ of Theorem~\ref{theorem: saxton thomason simple containers} and our bound on the average degree~$d$ in $H''$, inequality~\eqref{eq: linear hypergraph average degree}, we have
\begin{equation*}
\vert \mathcal{T}\vert \leq \vert \mathcal{C}\vert \leq 2^{\beta(d) k\binom{n}{2}}=k^{\Omega( n^{-1/(2r-1)})\binom{n}{2}},
\end{equation*}
so that there exist constants $C_0$,~$n_3>0$ such that for all $n\geq n_3$~sufficiently large, $\log_k\vert \mathcal{T}\vert \leq C_0n^{-1/(2\binom{N}{2}-1) \binom{n}{2}}$ and (iii) is satisfied. This establishes the statement of Theorem~\ref{theorem: multi-colour container} for $n\geq n_0=\max(n_1, n_2, n_3)$.
\end{proof}


\subsection{Extremal entropy and supersaturation}\label{subsection: entropy density, supersaturation}
In this section we obtain the two key ingredients needed in virtually all applications of containers, namely the existence of the limiting `density' of a property and a supersaturation result.

\begin{definition}
Let $\mathcal{P}$ be an order-hereditary property of $k$-colourings with $\mathcal{P}_n \neq \emptyset$ for every $n\in \N$. For every $n\in \N$, we define the \emph{extremal entropy} of $\mathcal{P}$ to be 
\[\ex(n, \mathcal{P})=\max\left\{\Ent(t)\, : \,  t \textrm{ is a $k$-colouring template for }K_n \textrm{ with } \langle t \rangle\subseteq \mathcal{P}_n\right\}.\]
\end{definition}
Note that this definition generalises the concept of the Tur\'an number: if $k = 2$, $F$ is a graph and $\mathcal{P} = \Forb(F)$, then $\ex(n, \mathcal{P}) = \ex(n, F)$.
\begin{proposition}\label{proposition: entropy density}
If\/ $\mathcal{P}$ is an order-hereditary property of $k$-colourings with\/ $\mathcal{P}_n \neq \emptyset$ for every $n\in \N$, then the sequence $\bigl(\ex(n, \mathcal{P})/\binom{n}{2}\bigr)_{n\in \N}$ tends to a limit $\pi(\mathcal{P})\in[0,1]$ as $n\rightarrow \infty$.
\end{proposition}
\begin{proof}
This is similar to the classical proof of the existence of the Tur\'an density. As observed after Definition~\ref{def: entropy}, $0\leq \Ent(t)\leq \binom{n}{2}$ for any $k$-colouring template $t$ of~$K_n$, so that $\ex(n, \mathcal{P})/\binom{n}{2}\in[0,1]$. It is therefore enough to show that $\left(\ex(n, \mathcal{P})/\binom{n}{2}\right)_{n\in \N}$ is nonincreasing.

Let $t$ be any $k$-colouring template for $K_{n+1}$ with $\langle t \rangle\subseteq \mathcal{P}_{n+1}$. For any $n$-subset $A\subseteq [n+1]$, the restriction $t_{\vert A}$ is a $k$-colouring template for $K_n$. Since $\mathcal{P}$ is order-hereditary, $\langle t \rangle\subseteq \mathcal{P}_{n+1}$ implies $\langle t_{\vert A} \rangle \subseteq \mathcal{P}_n$. By averaging over all choices of $A$, we have:
\begin{align*}
\frac{\Ent(t)}{\binom{n+1}{2}} &=\frac{1}{\binom{n+1}{2}} \log_k \left(\prod_{e\in [n+1]^{(2)}} \bigl\lvert t(e)\bigr\rvert\right)\\
&= \frac{1}{\binom{n+1}{2}}\log_k\left(\prod_{A\in [n+1]^{(n)}} \prod_{e\in A^{(2)}} \bigl\lvert t_{\vert A}(e)\bigr\rvert \right)^{1/n-1}\\
& =\frac{1}{n+1} \frac{1}{\binom{n}{2}}\sum_{A\in [n+1]^{(n)}} \Ent(t_{\vert A})\\
&\leq \frac{1}{n+1} \frac{1}{\binom{n}{2}}(n+1) \ex(n, \mathcal{P}).
\end{align*}
Thus $\ex(n+1, \mathcal{P})/\binom{n+1}{2}\leq \ex(n, \mathcal{P})/\binom{n}{2}$ as required and we are done.
\end{proof}
We call the limit $\pi(\mathcal{P})$ the \emph{entropy density} of $\mathcal{P}$. 
 Observe that the entropy density gives a lower bound on the \emph{speed} $\vert\mathcal{P}_n\vert$ of the property $\mathcal{P}$: for all $n\in\N$,
\begin{equation}\label{eq: entropy bound on P_n}
k^{\pi(\mathcal{P})\binom{n}{2}}\leq k^{\ex(n, \mathcal{P})}\leq \vert \mathcal{P}_n\vert.
\end{equation}
We shall show that the exponent in this lower bound is asymptotically tight.


\begin{lemma}[Supersaturation]\label{lemma: supersaturation}
Let $N \in \N$ be fixed and let $\mathcal{F}=\{c_1, c_2, \ldots, c_m\}$ be a nonempty collection of $k$-colourings of~$K_N$. Set\/ $\mathcal{P}=\Forb(\mathcal{F})$.  For every $\varepsilon$ with $0<\varepsilon<1$, there exist constants $n_0\in\N$ and~$C_{0}>0$ such that for any $k$-colouring template~$t$ for $K_n$ with $n\geq n_0$ and $\Ent(t)> \left(\pi(\mathcal{P})+\varepsilon\right)\binom{n}{2}$, there are at least~$C_{0} \varepsilon\binom{n}{N}$ pairs~$(A, c_i)$ with $A\in [n]^{(N)}$ and $c_i \in \mathcal{F}$ with $c_i\leq t_{\vert A}$.
\end{lemma}
\begin{proof}
We use a probabilistic bootstrapping technique. Given a $k$-colouring template~$t$ of~$K_m$, let $B(t)$ denote the number of pairs~$(A, c_i)$ with $A$ an $N$-set and $c_i \in \mathcal{F}$ such that $c_i\leq t_{\vert A}$.
Since every extra choice we are given above the extremal entropy~$\ex(m, \mathcal{P})$ must give rise to a new such pair, and since increasing the size of~$t(e)$ by $1$ for some edge~$e$ increases $\Ent(t)$ by at most $\log_k 2$, we have that
\begin{equation}\label{eq: bound on bad pairs}
B(t) \geq\frac{1}{\log_k 2} \bigl(\Ent(t)- \ex(m, \mathcal{P})\bigr).
\end{equation}

Now fix $\varepsilon >0$. By the monotonicity established in the proof of Proposition~\ref{proposition: entropy density}, there exists $n_1$ such that for all $n\geq n_1$ we have $\ex(n, \mathcal{P})\leq \pi(\mathcal{P})\binom{n}{2} +\frac{\varepsilon}{3}\binom{n}{2}$. Let $n_2=\max(n_1, \frac{1}{2}\log \frac{6}{\varepsilon})$ and $n_0=16n_2$. 

Let $t$ be a $k$-colouring template of~$K_n$, for some $n\geq n_0$. Suppose $\Ent(t)\geq \pi(\mathcal{P})\binom{n}{2}+\varepsilon \binom{n}{2}$. Let $p=\frac{8n_2}{n}$ and let $X$ be a random subset of $V(K_n)$ obtained by retaining each $v\in V(K_n)$ independently with probability~$p$ and casting it out otherwise. Denote by $t_{\vert X}$ the random $k$-colouring subtemplate of~$t$ induced by $X$.

Let $\mathcal{A}$ be the event that $\vert X\vert < n_1$. Since $\mathbb{E}\vert X\vert=np= 8n_2$, a standard Chernoff bound gives 
\begin{equation}\label{eq: chernoff bound on A}
\mathbb{P}(\mathcal{A})=\mathbb{P}\left(\vert X\vert < \frac{1}{8}\mathbb{E}\vert X\vert\right)\leq e^{-\frac{(7/8)^2\mathbb{E}\vert X\vert}{2}}<e^{-2n_2},
\end{equation}
which by our choice of $n_2$ is at most~$\varepsilon/6$. Now if $\mathcal{A}^c$ occurs, we may use~(\ref{eq: bound on bad pairs}) to bound $B(t_{\vert X})$ as follows: 
\begin{align}\label{eq:bound on induced bad pairs}
B\bigl(t_{\vert X}\bigr)&\geq \frac{1}{\log_k (2)} \left(\Ent\bigl(t_{\vert X}\bigr)- \ex\bigl(\vert X\vert, \mathcal{P}\bigr)\right) \nonumber\\
&\geq \mathbbm{1}_{\mathcal{A}^c} \frac{1}{\log_k (2)} \left(\Ent\bigl(t_{\vert X}\bigr)- \left(\pi+\frac{\varepsilon}{3}\right)\binom{\vert X\vert }{2}\right) \nonumber\\
&\geq \frac{1}{\log_k (2)}\left\{\sum_{e=\{xy\}\in E(K_n)} \mathbbm{1}_{\mathcal{A}^c}\mathbbm{1}_{x\in X}\mathbbm{1}_{y\in X}\log_k \bigl(\vert t(e)\vert\bigr)   -\left(\pi(\mathcal{P})+\frac{\varepsilon}{3}\right)\mathbbm{1}_{x\in X}\mathbbm{1}_{y\in X} \right\}.
\end{align}

Now the events $\mathcal{A}^c$, $\{x\in X\}$ and~$\{y\in X\}$ are all increasing events, and so by the Harris--Kleitman inequality they are positively correlated. Also for $x\neq y$, $\{x\in X\}$ and $\{y\in X\}$ are independent events, each occurring with probability~$p$.  It follows that
\begin{equation}\label{eq:induced entropy}
\sum_{e=\{xy\}\in E(K_n)} \mathbbm{1}_{\mathcal{A}^c}\mathbbm{1}_{x\in X}\mathbbm{1}_{y\in X}\log_k \bigl(\vert t(e)\vert\bigr) = \sum_{e=\{xy\}\in E(K_n)} \mathbb{P}(\mathcal{A}^c)p^2 \log_k \bigl(\vert t(e)\vert\bigr) = p^2 \bigl(1-\mathbb{P}(\mathcal{A})\bigr)\Ent(t)
\end{equation}
and
\begin{equation}\label{eq:size of X}
\sum_{e=\{xy\}\in E(K_n)} \mathbbm{1}_{x\in X}\mathbbm{1}_{y\in X} = p^2\binom{n}{2}.
\end{equation}
Taking the expectation of~\eqref{eq:bound on induced bad pairs} and applying \eqref{eq:induced entropy} and~\eqref{eq:size of X}, we thus have:
\begin{align*}
p^N B(t)&=\sum_{(A,c_i) \in B_t} \mathbb{P}(A\subseteq X)\\
&=\mathbb{E}B\bigl(t_{\vert X}\bigr)\\
&\geq\frac{1}{\log_k (2)} p^2\left(\bigl(1-\mathbb{P}(\mathcal{A})\bigr)\Ent(t) -\left(\pi(\mathcal{P})+\frac{\varepsilon}{3}\right) \binom{n}{2} \right)\\
&>\frac{1}{\log_k (2)} p^2\binom{n}{2}\left(\biggl(1-\frac{\varepsilon}{6}\biggr)\left(\pi(\mathcal{P})+\varepsilon\right)- \left(\pi(\mathcal{P})+\frac{\varepsilon}{3}\right)  \right)\\
&>\frac{1}{\log_k(2)}\frac{\varepsilon}{3} p^2\binom{n}{2},
\end{align*} 
where in the penultimate equality we used~\eqref{eq: chernoff bound on A}.
Dividing through by $p^N$, we deduce that
\begin{equation*}
B(t) > \frac{1}{\log_k(2)}\frac{\varepsilon}{3}p^{2-N}\binom{n}{2}=\frac{1}{\log_k (2)}\frac{\varepsilon}{3} \frac{n^{N-2}}{(8n_2)^{N-2}}\binom{n}{2}> \left(\frac{N!}{6(8n_2)^N\log_k (2)}\right)\varepsilon\binom{n}{N}.
\end{equation*}
This proves the lemma with $C_{0}=\frac{N!}{6(8n_2)^N\log_k 2}$ and $n_0=16n_2$. 
\end{proof}
\subsection{Speed of order-hereditary properties}\label{subsection: approximation of arbitrary hereditary properties}
In this section, we relate the speed of an order-hereditary property to its extremal entropy density and obtain container and counting theorems for arbitrary order-hereditary properties (i.e., properties defined by a possibly \emph{infinite} set of forbidden colourings).
\begin{theorem}\label{theorem: counting result for Forb(F), F finite hereditary families}
Let $N \in \N$ be fixed and let $\mathcal{F}=\{c_1, c_2, \ldots, c_m\}$ be a nonempty collection of $k$-colourings of~$E(K_N)$. Set\/ $\mathcal{P}=\Forb(\mathcal{F})$.  For all $\varepsilon>0$, there exists $n_0$ such that for all $n\geq n_0$ we have
\[k^{\pi(\mathcal{P})\binom{n}{2}}\leq \vert \mathcal{P}_n \vert \leq k^{\pi(\mathcal{P})\binom{n}{2}+\varepsilon \binom{n}{2}}.\]
\end{theorem}
\begin{proof}
Inequality~(\ref{eq: entropy bound on P_n}) already established the lower bound on the speed $n \mapsto \vert \mathcal{P}_n\vert$. 
For the upper bound, we first apply our multicolour container result. By Theorem~\ref{theorem: multi-colour container} applied to $\mathcal{P}$, for any $\eta>0$ there exists $n_1\in \N$ such that for all $n\geq n_1$, there exists a collection~$\mathcal{C}_n$ of at most~$k^{\eta\binom{n}{2}}$ templates such that $\mathcal{C}_n$ is a container family for $\mathcal{P}$ and there are at most~$\eta\binom{n}{N}$ pairs~$(A, c_i)$ with $A\in [n]^{(N)}$ and $c_i \in \mathcal{F}$ such that there exists some realisation~$c$ of a template in $\mathcal{C}_n$ with $c_{\vert A}$ isomorphic to $c_i$.

Provided we pick $\eta>0$~sufficiently small (less than~$C_{\varepsilon/2}$), we deduce from our supersaturation result, Lemma~\ref{lemma: supersaturation}, that there exists $n_2\in \N$ such that for all $n\geq n_2$, if $t$ is a $k$-colouring template for $E(K_n)$ for which there are fewer than~$\eta\binom{n}{N}$ pairs~$(A, c_i)$ with $A\in [n]^{(N)}$ and $c_i \in \mathcal{F}$ such that there exists $c\in\langle t \rangle$ with $c_{\vert A}$ isomorphic to~$c_i$,
then $\Ent(t)\leq \pi(\mathcal{P})\binom{n}{2}+\frac{\varepsilon}{2}\binom{n}{2}$.

Thus choosing $0<\eta<\min\left(\frac{\varepsilon}{2}, C_{\varepsilon/2}\right)$ and $n_0\geq \max(n_1, n_2)$, we have that for $n\geq n_0$ every template $t\in  \mathcal{C}_n$ has entropy at most $\pi(\mathcal{P})\binom{n}{2}+\frac{\varepsilon}{2}\binom{n}{2}$, whence we may at last bound above the speed of $\mathcal{P}$: for $n\geq n_0$,
\begin{equation*}
\vert \mathcal{P}_n\vert \leq \vert \mathcal{C}_n\vert k^{\max_{t\in \mathcal{C}_n} \Ent(t)}\leq k^{\eta\binom{n}{2} + \pi(\mathcal{P})\binom{n}{2} + \frac{\varepsilon}{2} \binom{n}{2}}\leq k^{ \pi(\mathcal{P})\binom{n}{2} + \varepsilon \binom{n}{2}}. \qedhere
\end{equation*}
\end{proof}

\begin{theorem}[Approximation of arbitrary order-hereditary properties]\label{theorem: approximation of general hereditary properties}
Let\/ $\mathcal{P}$ be an order-hereditary property with\/ $\mathcal{P}_n\neq\emptyset$ for every $n\in \N$ and let $\varepsilon>0$ be fixed. There exist constants $N$ and~$n_0\in\N$ and a nonempty symmetric family~$\mathcal{F}=\{c_1, c_2, \ldots, c_m\}$ of $k$-colourings of $E(K_N)$ such that for all $n\geq n_0$, we have
\begin{enumerate}[(i)]
\item $\mathcal{P}_n\subseteq \Forb(\mathcal{F})_n$, and 
\item $\vert\Forb(\mathcal{F})_n \vert \leq \lvert \mathcal{P}_n\rvert k^{\varepsilon \binom{n}{2}}$.
\end{enumerate}
\end{theorem}

\begin{proof}
For every $n\in\N$, let $\mathcal{F}_n$ denote the collection of $k$-colourings of $E(K_n)$ which are not in $\mathcal{P}_n$ (and thus, as $\mathcal{P}$ is order-hereditary, do not appear as induced subcolourings of any elements of $\mathcal{P}_{n'}$ with $n'\geq n$). Set $Q^n=\Forb(\bigcup_{m\leq n}\mathcal{F}_m)$ to be the order-hereditary property of $k$-colourings which avoids exactly the same $k$-colourings on at most $n$ vertices as $\mathcal{F}$. Note that by construction we have $\left(\mathcal{Q}^n\right)_m =\mathcal{P}_m$ for every $m\leq n$ and $\left(\mathcal{Q}^n\right)_m \supseteq \mathcal{P}_m$ for every $m>n$. We thus have a chain of inclusions
\[\mathcal{Q}^1\supseteq \mathcal{Q}^2\supseteq \cdots \supseteq \mathcal{Q}^n \supseteq \cdots \supseteq \mathcal{P}.\] 
Also, the sequence of entropy densities $\left(\pi(\mathcal{Q}^n)\right)_{n\in \N}$ is nonincreasing and bounded below by $\pi(\mathcal{P})$. We claim that $\lim_{n\rightarrow \infty} \pi(\mathcal{Q}^n)=\pi(\mathcal{P})$. Indeed, suppose this was not the case. Then there exists $\eta>0$ such that $\pi(\mathcal{Q}^n)>\pi(\mathcal{P})+\eta$ for all $n\in\N$. Since, as shown in the proof of Proposition~\ref{proposition: entropy density}, the sequence $\left(\ex(m, \mathcal{Q}^n)/\binom{m}{2}\right)_{m\in \N}$ is nonincreasing, there exists for every $n\in \N$ a $k$-colouring template~$t^n$ for $E(K_n)$ such that 
\[\langle t^n \rangle\subseteq \left(Q^n\right)_n= \mathcal{P}_n\] and 
\[\pi(\mathcal{P})+\eta<\Ent(t^n)/\binom{n}{2}\leq \ex(n, \mathcal{P})/\binom{n}{2},\]
contradicting Proposition~\ref{proposition: entropy density}.

Thus we must have $\lim_{n\rightarrow \infty} \pi(\mathcal{Q}^n)=\pi(\mathcal{P})$, as claimed. In particular there must exist some $N \in \N$ for which $\pi(\mathcal{Q}^N)<\pi(\mathcal{P})+\varepsilon/2$. 

Now for $n\geq n_1$, using the monotonicity and definition of $\ex(n, \mathcal{P})/\binom{n}{2}$ and $\pi(\mathcal{P})$ one last time, we have:
\begin{equation*}
k^{\pi(\mathcal{P})\binom{n}{2}}\leq k^{\ex(n, \mathcal{P})}\leq \vert \mathcal{P}_n\vert. 
\end{equation*}
On the other hand, by Corollary~\ref{theorem: counting result for Forb(F), F finite hereditary families} (applied to the property $\mathcal{Q}^N$ with parameter $\varepsilon/2$) there exists $n_2\in \N$ such that for all $n\geq n_2$ we have:
\begin{equation*}
\bigl\lvert (\mathcal{Q}^N)_n\bigr\rvert\leq k^{\pi(\mathcal{Q}^N)\binom{n}{2} +\frac{\varepsilon}{2} \binom{n}{2}}< k^{\pi(\mathcal{P})\binom{n}{2} +\varepsilon\binom{n}{2}}\leq \vert \mathcal{P}_n\vert k^{\varepsilon{\binom{n}{2}}}. 
\end{equation*}
Setting $n_0=\max(N, n_1, n_2)$ and observing that for $n\geq n_0$ we have $\left(\mathcal{Q}^N\right)_n = \Forb(\mathcal{F}_N)_n$ we see that the triple $(N, n_0, \mathcal{F}_N)$ satisfies the conclusion of the theorem. 
\end{proof}

\begin{corollary}\label{corollary: containers for arbitrary hereditary properties}
Let\/ $\mathcal{P}$ be an order-hereditary property with\/ $\mathcal{P}_n\neq\emptyset$ for every $n\in \N$ and let $\varepsilon>0$,~$m\in\N$ be fixed. There exists $n_0$ such that for any $n\geq n_0$ there exists a collection~$\mathcal{T}$ of templates of $k$-colourings of~$K_n$ satisfying:
\begin{enumerate}[(i)]
\item $\mathcal{T}$ is a container family for\/ $\mathcal{P}_n$;
\item every template~$t$ in\/ $\mathcal{T}$ has entropy at most $\left(\pi(\mathcal{P})+\varepsilon \right)\binom{n}{2}$;
\item in every realisation~$c$ of a template in\/ $\mathcal{T}$ there are at most~$\varepsilon \binom{n}{m}$ subsets~$A$ of~$[n]$ of order~$m$ such that the restriction of~$c$ to~$A$ is not order-isomorphic to an element of\/ $\mathcal{P}_m$;
\item $\vert \mathcal{T}\vert\leq k^{\varepsilon \binom{n}{2}}$.
\end{enumerate}
\end{corollary}
\begin{proof}
We let $\mathcal{Q}^n$ be defined as in the proof of Theorem~\ref{theorem: approximation of general hereditary properties}. As was shown there, we have $\lim_{n\to\infty}\pi(\mathcal{Q}^n)=\pi(\mathcal{P})$.  Therefore, for some $N$, we have $\pi(\mathcal{Q}^N)<\pi(\mathcal{P})+\frac{\varepsilon}{2}$. Without loss of generality we may take $N>m$.

Now let $C_{\varepsilon/2}$ be as given by Lemma~\ref{lemma: supersaturation} with $N$ as above and $\mathcal{F}=\bigcup_{m\le N} \mathcal{F}_m$ as in Theorem~\ref{theorem: approximation of general hereditary properties}, so that $\mathcal{Q}^N=\Forb(\mathcal{F})$. Fix $\delta<\min\{\varepsilon, C_{\varepsilon/2}\}$. 

Apply Theorem \ref{theorem: multi-colour container} with $\mathcal{F}$ as given and using $\delta$ in place of $\varepsilon$ in the statement of the theorem.  For large enough~$n$, consider the container family $\mathcal{C}_n$ for $\Forb(\mathcal{F})$ given by the theorem. Since $\mathcal{P}_n \subseteq (\mathcal{Q}^N)_n=\Forb(\mathcal{F})_n$ and $\mathcal{C}_n$ is a container family for $\Forb(\mathcal{F})_n$, $\mathcal{C}_n$ must also be a container family for $\mathcal{P}_n$.

We also know from Theorem \ref{theorem: multi-colour container} that for every template~$t\in\mathcal{C}_n$, there are at most~$\delta{\binom{n}{N}}$ pairs~$(A, c_i)$ with $A\in[n]^{(N)}$ and $c_i\in\mathcal{F}$ such that there is some realisation~$c$ of~$t$ whose restriction~$c_{\vert A}$ is isomorphic to $c_i$. Since $\delta<C_{\varepsilon/2}$, as long as $n$ is large enough Lemma~\ref{lemma: supersaturation} then implies that for every $t\in\mathcal{C}_n$, we have 
\begin{equation*}
\text{Ent}(t) \le\biggl(\pi\bigl(\mathcal{Q}^N\bigr)+\frac{\varepsilon}{2}\biggr){\binom{n}{2}}<\biggl(\pi(\mathcal{P})+\frac{\varepsilon}{2}+\frac{\varepsilon}{2}\biggr){\binom{n}{2}}
=\pi(\mathcal{P}){\binom{n}{2}}+\varepsilon{\binom{n}{2}}
\end{equation*}
which proves (ii).

To prove part (iii), consider a realisation~$c$ of a template~$t\in\mathcal{C}_n$. Let $B$ be a subset of~$[n]$ of size~$m$ such that the restriction of~$c$ to $B$ is isomorphic to an element of~$\mathcal{F}$. Now since $\Forb(\mathcal{F})$ is hereditary, and $m<N$, every superset of~$B$ of size~$N$ is also isomorphic to an element of~$\mathcal{F}$. However, from Theorem \ref{theorem: multi-colour container}, there are at most~$\varepsilon{\binom{n}{N}}$ subsets~$A$ of $[n]$ of order~$N$ such that the restriction of~$c$ to $A$ is isomorphic to an element of~$\mathcal{F}$. This gives us the following:
\begin{equation*}
\bigl\lvert \{B\in [n] \, : \, c|_B\in \mathcal{F}\} \bigr\rvert {\binom{n-m}{N-m}} \le \varepsilon {\binom{n} {N}}{\binom{N}{m}}.
\end{equation*}

Therefore there are at most~$\varepsilon {\binom{n}{N}}{\binom{N}{m}}/\binom{n-N}{N-M}=\varepsilon {\binom{n}{m}}$ subsets~$B$ of~$[n]$ of size $m$ such that the restriction of~$c$ to $B$ is isomorphic to an element of~$\mathcal{F}$. Now since $m<N$, $\mathcal{P}_m=\Forb(\mathcal{F})_m$ and so we have proved part~(iii).

Since $\delta<\varepsilon$, part (iv) of the corollary follows immediately from part (iii) of Theorem \ref{theorem: multi-colour container}.
\end{proof}

\begin{corollary}\label{corollary: speed of arbitrary hereditary properties}
Let\/ $\mathcal{P}$ be an order-hereditary property of $k$-colourings with\/ $\mathcal{P}_n\neq\emptyset$ for every $n\in \N$ and let $\varepsilon>0$ be fixed. There exists $n_0\in \N$ such that for all $n\geq n_0$,
\[k^{\pi(\mathcal{P})\binom{n}{2}}\leq \vert \mathcal{P}_n\vert \leq k^{\pi(\mathcal{P})\binom{n}{2}+\varepsilon \binom{n}{2}}.\]
\end{corollary}
\begin{proof}
The lower bound is given by inequality~(\ref{eq: entropy bound on P_n}). For the upper bound, we apply Corollary~\ref{corollary: containers for arbitrary hereditary properties} (with parameter~$\varepsilon/2$) to obtain for all $n\geq n_0$ a family of templates~$\mathcal{C}_n$ satisfying properties (i), (ii) and~(iv). Thus for $n\geq n_0$,
\[\vert\mathcal{P}_n\vert \leq \vert\mathcal{C}_n\vert k^{\max_{t\in \mathcal{C}_n} \Ent(t)} \leq k^{\frac{\varepsilon}{2}\binom{n}{2}} k^{\pi(\mathcal{P})\binom{n}{2}+\frac{\varepsilon}{2}\binom{n}{2}}= k^{\pi(\mathcal{P})\binom{n}{2}+\varepsilon\binom{n}{2}}. \qedhere\]
\end{proof}

\subsection{Stability and characterization of typical colourings}\label{subsection: stability}

\begin{definition}
	A ($k$-colouring) \emph{template sequence} is a sequence $\mathbf{t}=(t_n)_{n\in \N}$, where $t_n$ is a ($k$-colouring) template for $K_n$. Given a family $\mathcal{S}$ of $k$-colouring template sequences, we denote by $\langle \mathcal{S} \rangle$ the sequence of realisations from $\mathcal{S}$, i.e., $\langle \mathcal{S} \rangle=\{c: \ c\in \langle t \rangle \textrm{ for some }t\in \mathcal{S}\}$.
\end{definition}

\begin{definition}[Edit distance]\label{definition: edit distance}
	The \emph{edit distance}~$\rho(s,t)$ between two $k$-colouring templates $s$,~$t$ of~$K_n$ is the number of edges~$e\in E(K_n)$ on which $s(e)\neq t(e)$. We also define the \emph{edit distance}~$\rho(c,t)$ between a $k$-colouring~$c$ and a $k$-colouring template~$t$ to be the number of edges~$e\in E(K_n)$ on which $c\notin t(e)$.

	Finally, the \emph{edit distance} $\rho(\mathcal{S}, t)$ between a $k$-colouring template~$t$ and a family~$\mathcal{S}$ of $k$-colouring template sequences is 
	\[\rho(\mathcal{S}, t):=\min_{s\in \mathcal{S}}\rho(s,t).\]  
	Similarly, we define $\rho(\langle \mathcal{S} \rangle, c)$ for a $k$-colouring $c$ to be $\min_{c'\in\langle \mathcal{S} \rangle} \rho(c', c)$. 
\end{definition}

\begin{definition}[Strong stability]\label{definition: strong stability family}
Let $\mathcal{P}$ be a hereditary property of $k$-colourings of~$K_n$. A family~$\mathcal{S}$ of $k$-colouring template sequences is a \emph{strong stability family} for $\mathcal{P}$ if for all $\varepsilon>0$, there exist $\delta>0$ and $m$,~$n_0\in \N$ such that for all $n\geq n_0$, every $k$-colouring template~$t$ for $K_n$ satisfying
	\begin{enumerate}[(i)]
		\item (almost extremality) $\Ent(t)\geq (\pi(\mathcal{P})-\delta)\binom{n}{2}$;
		\item (almost locally  in $\mathcal{P}$) there are at most~$\delta \binom{n}{m}$ pairs~$(A, c)$ where $A$ is an $m$-subset of~$V(K_n)$ and $c$ is a realisation of~$t_{\vert_A}$ such that $c\notin \mathcal{P}_m$
	\end{enumerate}
	must lie within edit distance at most $\varepsilon \binom{n}{2}$ of a template $t_n$ drawn from a template sequence $\mathbf{t} \in \mathcal{S}$.
\end{definition}

\begin{theorem}\label{theorem: strong stability and containers}
Let\/ $\mathcal{P}$ be an order-hereditary property of $k$-colourings and suppose $\mathcal{S}$ is a strong stability family for\/ $\mathcal{P}$.
For all $\varepsilon>0$, there exists $n_0\in \N$ such that for all $n\geq n_0$ there are at most~$\varepsilon \vert \mathcal{P}_n\vert$ colourings~$c\in\mathcal{P}_n$ with $\rho(c, \langle \mathcal{S} \rangle)>\varepsilon \binom{n}{2}$.
\end{theorem}

\begin{proof}

Let $\varepsilon>0$, and let $\delta$ be as given by Definition~\ref{definition: strong stability family} for $\mathcal{S}$. Apply Corollary \ref{corollary: containers for arbitrary hereditary properties} to $\mathcal{P}$ with $\varepsilon_1<\delta$ to get (for large enough~$n$) a container family $\mathcal{C}_n$ for $\mathcal{P}_n$.

Now remove from $\mathcal{C}_n$ any templates $t$ with $\text{Ent}(t)<(\pi(\mathcal{P})-\delta){\binom{n}{2}}$, to get $\mathcal{C}_n'\subseteq\mathcal{C}_n$.  By part~(iii) of Corollary~\ref{corollary: containers for arbitrary hereditary properties} and Corollary~\ref{theorem: counting result for Forb(F), F finite hereditary families}, the number of elements of $\mathcal{P}_n$ which are not realisable from a template $t\in \mathcal{C}_n'$ is then at most
\[\vert \mathcal{C}_n\vert k^{(\pi(\mathcal{P})-\delta){\binom{n}{2}}} \le k^{(\pi(\mathcal{P})+\varepsilon_1-\delta){\binom{n}{2}}} \leq\vert \mathcal{P}_n\vert k^{(\varepsilon_1-\delta){\binom{n}{2}}}.\]
Since $\varepsilon_1<\delta$, for large enough~$n$ the right hand side is less than~$\varepsilon \vert \mathcal{P}_n\vert$.

Now let $c$ be a member of $\mathcal{P}_n$ which \textit{is} realisable from a template $t\in\mathcal{C}_n'$. Since $t\in\mathcal{C}_n'$, we have $\text{Ent}(t)\ge (\pi(\mathcal{P})-\delta){\binom{n}{2}}$. Also since $\varepsilon_1<\delta$, Corollary \ref{corollary: containers for arbitrary hereditary properties} implies that $t$ satisfies condition~(ii) of Definition \ref{definition: strong stability family} and so there is a template~$s_n\in \mathcal{S}$ such that $\rho(t,s_n)<\varepsilon{\binom{n}{2}}$. Since $c$ realises $t$, this readily implies that $\rho(c,s_n)<\varepsilon{\binom{n}{2}}$ and so $\rho(c,\mathcal{S})<\varepsilon{\binom{n}{2}}$.
\end{proof}

\subsection{Transference}\label{subsection: transference}
\begin{definition}[Multicolour monotonicity]
An order-hereditary property $\mathcal{P}$ of $k$-colourings is \emph{monotone} with respect to colour $i\in[k]$ if whenever $c$ is a $k$-colouring of~$K_n$ which lies in $\mathcal{P}$ and $e$ is any edge of~$K_n$, the colouring~$\tilde{c}$ obtained from $c$ by changing the colour of~$e$ to $i$ also lies in $\mathcal{P}$.
\end{definition}
\begin{definition}[Meet of two templates]
Given two $k$-colouring templates $t$,~$t'$ of~$K_n$ which have at least~one realisation in common, we denote by $t\wedge t'$ the template with $(t\wedge t')(e) =t(e)\cap t'(e)$ for each $e\in E(K_n)$; we refer to $t\wedge t'$ as the \emph{meet} of $t$~and~$t'$.  More generally, given a set~$\mathcal{S}$ of $k$-colouring templates of~$K_n$ and a $k$-colouring template~$t'$ of~$K_n$, we denote by $\mathcal{S}\wedge t'$ the collection~$\{t\wedge t': \ t\in\mathcal{S}\}$. 
\end{definition}

\begin{definition}[Complete, random and constant templates]
Let $T_n$ denote the \emph{complete $k$-colouring template} for $K_n$, that is, the unique template allowing all $k$ colours on all edges. Given a fixed colour $i\in[k]$ and $p\in[0,1]$, we define the \emph{$p$-random template} $T_{n,p}=T_{n,p}(i)$ to be the random template for a $k$-colouring of~$K_n$ obtained by letting
\[T_{n,p}(e)=\left\{\begin{array}{ll}
[k] & \textrm{with probability $p$}\\
\{i\} & \textrm{otherwise,}\end{array} \right.\]
independently for each edge~$e\in E(K_n)$. Finally, let $E_n=E_n(i)$ denote the \emph{$i$-monotone} template with $E_n(e)=\{i\}$ for each $e\in E(K_n)$.
\end{definition}
The $p$-random template~$T_{n,p}$ is our $k$-colouring analogue of the celebrated Erd{\H o}s--R\'enyi binomial random graph~$G_{n,p}$, while the zero-entropy template $E_n$ is a $k$-colouring analogue of the empty graph. Just as extremal theorems for the graph~$K_n$ can be reproved in the sparse random setting of~$G_{n,p}$, so also extremal entropy results for $i$-monotone properties in $T_n$ can be transferred to the $T_{n,p}(i)$ setting.

Of course, there are other very natural notions of random templates.  For example, for each edge~$e \in E(K_n)$, one could choose the palette available at $e$ by including colours independently with probability~$p$, independently of all other edges.  Another possibility is discussed in Remark~\ref{remark: partial order}.

\begin{definition}[Relative entropy]
Let $\mathcal{P}$ be a property of $k$-colourings of complete graphs that is monotone with respect to colour $i$. Let $t$ be a template for a $k$-colouring of~$K_n$ which contains the $i$-monotone template $E_n(i)$. We define the \emph{extremal entropy of $\mathcal{P}$ relative to $t$} to be:
\[\ex(t,\mathcal{P}):=\max\bigl\{\Ent(t') \, : \, v(t')=n,\ t'\leq t, \ \langle t' \rangle \subseteq \mathcal{P}_n\bigr\}.\]
\end{definition}
Note that this notion of relative entropy extends the notion of extremal entropy introduced in Definition~\ref{def: entropy}: $\ex(n, \mathcal{P})=\ex(T_n,\mathcal{P})$. Our next theorem states that for $p$ not too small, with high probability the extremal entropy of $\mathcal{P}$ relative to a $p$-random template $T_{n,p}$ is $p\ex(n, \mathcal{P})+o(pn^2)$. 
\begin{theorem}[Transference]\label{theorem: transference}
Let\/ $\mathcal{P}$ be an order-hereditary, $i$-monotone property of $k$-colourings of complete graphs defined by forbidden colourings on at most $N$ vertices. Let $p=p(n)\gg n^{-1/(2\binom{N}{2}-1)}$ and let $T$ denote an instance of the $p$-random template $T_{n,p}(i)$. For any fixed $\varepsilon>0$, with high probability
\begin{equation}\label{eq:TransferenceEntropy}
p\left(\ex(n,\mathcal{P})-\varepsilon n^2 \right) \leq \ex(T,\mathcal{P})\leq p\left(\ex(n,\mathcal{P})+2\varepsilon n^2 \right).
\end{equation}
\end{theorem}
\begin{proof}
The result is trivial for $N=1$, so assume $N\geq 2$. Let $\varepsilon>0$ be fixed. Apply Theorem~\ref{theorem: multi-colour container} to find for all $n\geq n_0=n_0(\varepsilon)$ a collection~$\mathcal{C}_n$ of at most $k^{C_0n^{-1/(2\binom{N}{2}-1)} \binom{n}{2}}$ $k$-colouring templates covering $\mathcal{P}_n$, each having entropy at most $\ex(n, \mathcal{P})+\varepsilon \binom{n}{2}$.

Now let $T$ be an instance of the random $k$-colouring template $T_{n,p}$. By the container property, every $k$-colouring template in $\mathcal{P}_n\wedge T$ is a subtemplate of a template in $\mathcal{C}_n\wedge T$. Let us therefore estimate the maximal entropy attained in that family. For each $t \in \mathcal{C}_n$, we have 
\[\Ent(t\wedge T)=\sum_{e\in E(K_n)} \log_k \bigl\lvert t\wedge T(e)\bigr\rvert=\sum_{e: \ T(e)=[k]} \log_k\bigl\lvert t(e)\bigr\rvert.\]
By the Chernoff bound (\ref{equation: Chernoff}), for all $t \in \mathcal{C}_n$, we have
\begin{equation*}\mathbb{P}\biggl(\Ent(t\wedge T)\geq p\ex(n, \mathcal{P})+2p\varepsilon \binom{n}{2}\biggr) \leq \exp\biggl(-\varepsilon \frac{pn^2}{8}\biggr). 
\end{equation*}
In particular, with probability at least 
\begin{equation*}
1- \lvert \mathcal{C}_n\wedge T\rvert e^{-\varepsilon \frac{pn^2}{8}}\geq 1- \left(k^{n^{-1/(2\binom{N}{2}-1)}}e^{-\frac{\varepsilon  p}{4}}\right)^{n^2/2},
\end{equation*}
the maximum entropy in $\mathcal{C}_n\wedge T$ is at most $p\left(\ex(n,\mathcal{P})+\varepsilon n^2 \right)$. For $\varepsilon p > 8(\log k)  n^{-1/(2\binom{N}{2}-1)}$, this probability is
\[
1-\exp\Biggr(-\Omega\Biggl(n^{2-\frac{1}{2\binom{N}{2}-1}} \Biggr)\Biggr)=1-o(1).
\]
In particular, for $p \gg n^{-1/(2\binom{N}{2}-1)}$ and any fixed $\varepsilon$, whp the extremal entropy of~$\mathcal{P}$ relative to an instance~$T$ of $T_{n,p}$ will be at most~$p\ex(n, \mathcal{P})+ 2\varepsilon pn^2$. This establishes the upper bound in~\eqref{eq:TransferenceEntropy}.

For the lower bound, consider a maximum entropy template $t_ {\star}$ for $\mathcal{P}_n$. Then using the Chernoff bound (\ref{equation: Chernoff}) again, we have that
\[\mathbb{P}\Bigl(\Ent(t_{\star}\wedge T)\geq p\bigl(\ex(n,\mathcal{P})-\varepsilon n^2\bigr)\Bigr) =1-\exp \bigl(-\Omega\bigl(\varepsilon p n^2\bigr)\bigr)=1-o(1),\]
 provided $p\gg n^{-2}$ and $\varepsilon$ is fixed. Together with the argument above, this establishes that for
\[p\gg \max\left(n^{-1/(2\binom{N}{2}-1)}, n^{-2}\right)= n^{-1/(2\binom{N}{2}-1}\]
and any fixed $\varepsilon>0$, whp 
\[p\ex(n,\mathcal{P})-p\varepsilon n^2\leq \ex(T,\mathcal{P})\leq p\ex(n,\mathcal{P})+2p\varepsilon n^2,  \]
as required.
\end{proof}

\begin{remark}
	The bound on $p$ required in Theorem~\ref{theorem: transference} is not best possible in general (see~\cite{Conlon14}). 
	This bound can be improved by using the more powerful container theorems of~\cite{BaloghMorrisSamotij15} and~\cite{SaxtonThomason15} rather than the simple hypergraph container theorem, Theorem~\ref{theorem: saxton thomason simple containers}. However, as this is not the focus of this paper, we do not pursue such improvements further here.
\end{remark}

We can extend Theorem~\ref{theorem: transference} to cover general order-hereditary properties.
\begin{corollary}\label{corollary: transference for general properties}
Let\/ $\mathcal{P}$ be an order-hereditary, $i$-monotone property of $k$-colourings of complete graphs. Let $p=p(n)$ be a sequence of probabilities satisfying $\log (1/p)=o(\log n)$, and let $T$ denote an instance of the $p$-random template $T_{n,p}(i)$. For any fixed $\varepsilon>0$, with high probability
\[p\left(\ex(n,\mathcal{P})-\varepsilon n^2\right)\leq \ex(T,\mathcal{P})\leq p\left(\ex(n,\mathcal{P})+3\varepsilon n^2\right).\]
\end{corollary}
\begin{proof}
Let $\varepsilon>0$ be fixed. As in Theorem~\ref{theorem: approximation of general hereditary properties}, approximate $\mathcal{P}$ from above by some property~$\mathcal{Q}$ defined by forbidden colourings on at most $N$ vertices and satisfying $\mathcal{P}\subseteq\mathcal{Q}$ and $\pi(\mathcal{Q}) \leq \pi(\mathcal{P}) +\varepsilon$. For $n$~sufficiently large, $\ex(n,\mathcal{Q})\leq \pi(\mathcal{Q})\binom{n}{2}+\varepsilon \binom{n}{2}\leq \pi(\mathcal{P})\binom{n}{2}+\varepsilon n^2$. Applying Theorem~\ref{theorem: transference} to $\mathcal{Q}$, and noting that our condition on $p$ ensures $p\gg n^{-(1/(2\binom{N}{2}-1))}$ for all $N \in \N$, we obtain the desired result. 
\end{proof}
	The container approach allows us to transfer stability results as well as extremal ones to the sparse random setting. It is straightforward to adapt the proof of 
	Theorem~\ref{theorem: strong stability and containers} to prove the following:
\begin{theorem}\label{theorem: sparse random stability}
Let\/ $\mathcal{P}$ be an order-hereditary, $i$-monotone property of $k$-colourings of complete graphs. Suppose\/ $\mathcal{S}$ is a strong stability family for\/ $\mathcal{P}$. Then for all $p=p(n)$ with $\log (1/p)=o(\log n)$ the subtemplates of $T_{n,p}$ with extremal entropy for\/ $\mathcal{P}$ are close to\/ $\mathcal{S}\wedge T_{n,p}$.

More precisely: whp, all subtemplates~$t$ of~$T_{n,p}$ with $\Ent(t)\geq (1+o(1))p \ex(n, \mathcal{P})$ satisfy $\rho(t, \mathcal{S}\wedge T_{n,p})=o(n^2)$.
\end{theorem}
\begin{proof}
First apply Theorem~\ref{theorem: strong stability and containers} to the near-extremal templates for $k$-colourings of~$K_n$ to show they are close to templates from $\mathcal{S}$. Then proceed to the proof of Theorem~\ref{theorem: transference} and Corollary~\ref{corollary: transference for general properties} as before.
\end{proof}

\begin{remark}\label{remark: partial order}
To conclude this section we remark that in the multicolour setting, other, more sophisticated notions of transference are possible. Explicitly, equip the set of colours~$[k]$ with a partial order~$\leq$.
Call a $k$-colouring property~$\mathcal{P}$ \emph{$\leq$-monotone} if whenever $i\leq j$ and $c\in \mathcal{P}_n$ has $c(e)=j$ for some $e\in E(K_n)$, then the colouring obtained from $c$ by changing $e$'s colour to~$i$ also lies in $\mathcal{P}$. Define a probability distribution~$\mu$ on the collection of all palettes containing an $\leq$-minimal element. From this distribution we obtain a new, different notion of a sparse random template: let $T_n(\mu)$ be the random template obtained by assigning each edge~$e\in E(K_n)$ a $\mu$-random subset of~$[k]$ independently of all other edges. As $\mathcal{P}$ is $\leq$-monotone, we know at least~one realisation of~$T_n(\mu)$ must lie in $\mathcal{P}$, and we can ask about the extremal entropy~$\ex(T_n(\mu), \mathcal{P})$ as before. This appears more complicated than the setting we explored in this section, but it would be a natural direction for further work.
\end{remark}




\section{Other discrete structures}\label{section: other structures}
Our container results so far allow us to compute the speed of (dense) order-hereditary properties of $k$-colourings of~$K_n$, as well as to characterise typical colourings and transfer extremal entropy results to the sparse random setting. However, as we show in this section, the container theory of Saxton--Thomason and Balogh--Morris--Samotij is robust enough to cover $k$-colourings of many other interesting discrete structures: in essence, all we need to apply a container theorem is a sensible notion of substructure. In this respect, container theory is reminiscent of the theory of flag algebras developed by Razborov~\cite{Razborov07}.

In the following subsections, we outline how our $k$-colouring extensions of the container theorems of Saxton--Thomason can be applied to tournaments, oriented graphs, directed graphs, multipartite graphs, graph sequences, edge- and vertex-subsets of the hypercubes and hypergraphs.

\subsection{Tournaments, oriented graphs and directed graphs}\label{subsection: oriented}
Tournaments, oriented graphs and directed graphs are central objects of study in discrete mathematics and computer science, with a number of applications both to other branches of mathematics and to real-world problems. As we show below, our framework of order-hereditary properties for $k$-colourings of~$K_n$ allow us to cover these structures with our container, supersaturation, counting, transference and characterization theorems. This is perhaps the most interesting application of our work, as containers had not been successfully applied to the directed setting before (see Section~\ref{subsection: digraph examples} for a discussion, or the remark after Corollary~3.4 in K\"uhn, Osthus, Townsend and Zhao~\cite{KuhnOsthusTownsendZhao14}).

Formally, a \emph{directed graph}, or \emph{digraph}, is a pair $D=(V,E)$, where $V=V(D)$ is a set of vertices and $E=E(D)\subset V\times V$ is a collection of \emph{ordered} pairs from $V$. By convention, we write $\vec{ij}$ to denote $(i,j)\in E$. Note that we could have both $\vec{ij}\in E(D)$ and $\vec{ji}\in E(D)$, in which case we say that $ij$ is a \emph{double edge} of $D$.

An \emph{oriented graph}, or \emph{orgraph}, is a digraph $\vec{G}$ in which for each pair $ij\in V(\vec{G})$ at most one of $\vec{ij}$ and $\vec{ji}$ lies in $E(\vec{G})$. A \emph{tournament} $\vec{T}$ is a digraph in which for each pair $ij\in V(\vec{T})$ exactly one $\vec{ij}$ and $\vec{ji}$ lies in $E(\vec{T})$ --- or, more helpfully a tournament can be viewed as an orientation of the edges of the complete graph.

A \emph{monotone (decreasing)} property of digraphs/orgraphs is a property of digraphs/orgraphs which is closed with respect to taking subgraphs (i.e.\ closed under the deletion of vertices and oriented edges). A \emph{hereditary} property of digraphs/orgraphs/tournaments is a property of digraphs/orgraphs/tournaments which is closed with respect to taking induced subgraphs.

\begin{observation}[Key observation]\label{observation: encoding of digraphs}
Tournaments, oriented graphs and directed graphs on the labelled vertex set $[n]$ vertices can be encoded as $2$--, $3$-- and $4$-colourings of~$K_n$. Moreover, under this encoding, hereditary properties of tournaments, oriented graphs and directed graphs correspond to order-hereditary properties of $2$--, $3$-- and $4$-colourings of~$K_n$.
\end{observation}	
\begin{proof}
	Given a directed graph $D$ on $[n]$, we define a colouring $c$ of $E(K_n)$ by setting for each pair~$ij\in [n]^{(2)}$ with $i<j$
	\[c(ij):=\left\{\begin{array}{ll}1 & \textrm{if neither of $\vec{ij}$,~$\vec{ji}$ lies in $E(D)$,}\\
	2 & \textrm{if $\vec{ij}\in E(D)$, $\vec{ji}\notin E(D)$,}\\
	3 & \textrm{if $\vec{ij}\notin E(D)$, $\vec{ji}\in E(D)$,}\\
	4 & \textrm{if both of $\vec{ij}$,~$\vec{ji}$ lie in $E(D)$.}
	\end{array}
	 \right.\]
	Observation~\ref{observation: encoding of digraphs} is immediate from this colouring and our definition of order-hereditary properties. Tournaments correspond to colouring by colours~$\{2,3\}$, oriented graphs to colourings by colours~$\{1,2,3\}$ and digraphs to colourings with the full palette of possible colours~$\{1,2,3,4\}$.
\end{proof}	
\begin{remark}
Monotone properties of digraphs/orgraphs correspond to order-hereditary properties of $4$--/$3$-colourings of~$K_n$ which are monotone with respect to colour~$1$.	
\end{remark}
\begin{corollary}
If\/ $\mathcal{P}$ is a hereditary property of digraphs/orgraphs/tournaments defined by forbidden configurations on at most~$N$ vertices and let $k=4/3/2$, then the conclusions of Theorems \ref{theorem: multi-colour container}, \ref{theorem: counting result for Forb(F), F finite hereditary families}, \ref{theorem: strong stability and containers}, and~\ref{theorem: transference} hold for\/ $\mathcal{P}$.\qed
\end{corollary}
\begin{corollary}
If\/ $\mathcal{P}$ is a hereditary property of digraphs/orgraphs/tournaments and let $k=4/3/2$, then the conclusions of Corollary~\ref{corollary: containers for arbitrary hereditary properties} and Theorem~\ref{theorem: strong stability and containers} hold for\/ $\mathcal{P}$.\qed
\end{corollary}

In particular, we have general counting, stability and transference results for hereditary properties of digraphs, orgraphs and tournaments. As mentioned earlier, this overcomes an obstruction to the extension of containers to the digraph setting.


\subsection{Other host graphs: grids, multipartite graphs and hypercubes}\label{subsection: other host graphs}

Our results thus far concern $k$-colourings of the complete graph~$K_n$. However, the container theory of Balogh--Morris--Samotij and Saxton--Thomason is more than robust enough to cover the case of $k$-colouring of other graphs, in particular $q$-partite graphs and hypercube graphs (defined below). To tackle such cases, we need to define notions of a template and of extremal entropy \emph{relative} to a graph.

\begin{definition}[Template and entropy relative to a graph]\label{definition: template/entropy relative to G}
	Let $G$ be a graph. A \emph{template} for a $k$-colouring of~$G$ is a function~$t$, associating to each edge~$e$ of~$G$ a non-empty list of colours~$t(e)\subseteq [k]$. The set of all such templates is denoted by $[k]^{E(G)}$.

	 Given a template $t\in [k]^{E(G)}$, we write $\langle t \rangle$ for the collection of \emph{realisations} of $t$, that is, the collection of $k$-colourings~$c$ of~$E(G)$ such that $c(e)\in t(e)$ for every edge~$e\in E(G)$. The \emph{entropy} of a $k$-colouring template~$t$ of~$G$ is 
	\[\Ent(t):= \sum_{e\in E(G)}\log_k\vert t(e)\vert.\]
\end{definition}
Observe that $0\leq \Ent(t)\leq e(G)$ and $\vert \langle t \rangle\vert =k^{\Ent(t)}$.
\begin{definition}[Extremal entropy relative to a graph sequence]\label{definition: extremal entropy/relative to G}
	Let $\Graphseq=(G_n)_{n \in \N}$ be a sequence of graphs on linearly ordered vertex sets. A \emph{$k$-colouring property of $\Graphseq$} is a sequence $\mathcal{P}=(\mathcal{P}_n)_{n\in \N}$, where $\mathcal{P}_n$ is a collection of $k$-colourings of $G_n$. The \emph{extremal entropy} of $\mathcal{P}$ \emph{relative to $\Graphseq$} is
	\[\ex(\Graphseq, \mathcal{P})=\ex(G_n, \mathcal{P}_n):=\max\left\{\Ent(t) \,: \, t\in [k]^{E(G_n)}, \ \langle t \rangle\subseteq \mathcal{P}_n\right\}.\]
\end{definition}
Thus the work in the previous section was concerned with $k$-colouring properties of~$\Complete=(K_n)_{n\in \N}$. However, many other natural graph sequences have been studied from an extremal and/or counting perspective. Examples of such sequences include:
\begin{itemize}
	\item $\Path=(P_n)_{n\in \N}$, the sequence of paths on $[n]$, $P_n=([n], \{i(i+1): \ 1\leq i \leq n-1\})$;
	\item $\Grid= (P_n\times P_n)_{n \in \N}$, the sequence of $n\times n$ grids~$P_n\times P_n$ obtained by taking the Cartesian product of $P_n$ with itself, or, more generally for $(a,b)\in \N^2$ the sequence of rectangular grids~$\Grid(a,b)=(P_{an}\times P_{bn})_{n\in \N}$;
	\item $\Branch_b=(B_{b,n})_{n \in \N}$, the sequence of $b$-branching trees with $n$ generations from a single root;
	\item $\Complete_q=(K_q(n))_{n \in \N}$, the sequence of complete balanced $q$-partite graphs on $qn$ vertices;
	\item $\Hypercube=(Q_n)_{n\in \N}$, the sequence of $n$-dimensional discrete hypercube graphs~$Q_n=(\{0,1\}^n, \{\vector{x}\vector{y}: \ \vector{x}_i =\vector{y}_i \textrm{ for all but exactly~one index }i\}$.
\end{itemize}
Outside of extremal combinatorics, the sequences $\Hypercube$ and~$\Branch_2$ are of central importance in theoretical computer science and discrete probability (they represent $n$-bit sequences and binary search trees respectively), while the sequence~$\Grid$ has been extensively studied in the context of percolation theory, in particular with respect to crossing probabilities.

Each of the graph sequences above comes equipped with a natural notion of `substructure' --- subpaths of a path, subgrids of a grid, subtrees of a branching tree, subgraphs of a $q$-partite graph, subcubes of a hypercube --- which in turn  gives rise to a notion of an (order-) hereditary property. Our next theorem says that we can find a generalization of Theorem~\ref{theorem: multi-colour container} for any graph sequence~$\Graphseq$ with a notion of substructure which is `sufficiently rich'.

\begin{definition}[Embedding]
	Let $\Graphseq$ be a graph sequence and let $n \geq N$. An \emph{order-preserving embedding} of~$G_N$ into $G_n$ is an injection $\phi: \ V(G_N)\rightarrow V(G_n)$ such that $\phi$ preserves edges and the linear order on the vertices: if $x\leq y$ in $G_N$ then $\phi(x)\leq \phi(y)$ in $G_n$, and if $xy\in E(G_N)$ then $\phi(x)\phi(y)\in E(G_n)$. We denote by $\binom{G_n}{G_N}$ the number of order-preserving embeddings of~$G_N$ into $G_n$.	
\end{definition}
\begin{definition}[Intersecting embedding]
Let $N_1$,~$N_2 \leq n$.  An \emph{$i$-intersecting (order-preserving) embedding} of~$(G_{N_1}, G_{N_2})$ into $G_n$ is a function $\phi: \ V(G_{N_1})\sqcup V(G_{N_2})\rightarrow V(G_n)$ such that:
\begin{enumerate}[(i)]
	\item the restriction of $\phi$ to either of $V(G_{N_1})$ or $V(G_{N_2})$ is an order-preserving embedding;
	\item $\vert \phi(E(G_{N_1}))\cap \phi(E(G_{N_2}))\vert=i$.
\end{enumerate}
We denote by $I_i( (G_{N_1}, G_{N_2}), G_n)$ the number of $i$-intersecting embeddings of of $(G_{N_1}, G_{N_2})$ into $G_n$, and set
\[I(N, n):=\sum_{1<i<e(G_N)} I_i\bigl( (G_{N}, G_{N}), G_n\bigr).\]		
\end{definition}
\begin{definition}[Good graph sequence]\label{definition: good graph sequence}
A graph sequence $\Graphseq$ is \emph{good} if all of the following hold:
\begin{enumerate}[(i)]
	\item $e(G_n)\rightarrow \infty$ (`the graphs in the sequence become large');
	\item for all $N \in \N$ with $N\geq3$, $\binom{G_n}{G_N}\rightarrow \infty$ (`the sequence has many embeddings of~$G_N$');
	\item for all $N \in \N$ with $N\geq 3$, $e(G_n)I(N,n)\big/ \binom{G_n}{G_N}^2 \rightarrow 0$ as $n\rightarrow \infty$ (`most pairs of embeddings of~$G_N$ are edge-disjoint').
\end{enumerate}	
\end{definition}
Roughly speaking, a graph sequence is good if it is sufficiently rich in embeddings --- there must be many ways of embeddings $G_N$ into $G_n$ relative to the number of edges. In practice, (i) and (ii) are obvious, and only checking (iii) will require a calculation. It is easy to verify that the sequences $\Hypercube$, $\Complete_q$ and~$\Grid$ defined above are good, but that the sequences $\Path$ and~$\Branch_b$, for instance, have too few embeddings and so fail to be good. Our next result says that we have multicolour container theorems for good graph sequences. Furthermore, as we shall show in Section~\ref{subsection: nonexample containers fail for paths}, some form of the `goodness' assumption is necessary --- the conclusion of Theorem~\ref{theorem: general multicolour container result for graph sequences} fails for the sequence~$\Path$, for instance.

Given a collection $\mathcal{F}$ of forbidden $k$-colourings of $G_N$, denote by $\Forb_{\Graphseq}(\mathcal{P})$ the order-hereditary property of $k$-colourings of $\Graphseq$ of not containing an embedding of a colouring in $\mathcal{F}$.
\begin{theorem}\label{theorem: general multicolour container result for graph sequences}
	Let $\Graphseq$ be a good graph sequence, and let $k$,~$N \in \N$. Let $\mathcal{F}$ be a collection of forbidden $k$-colourings of~$G_N$ and let\/ $\mathcal{P}=\Forb_{\Graphseq}(\mathcal{P})$.  For any $\varepsilon>0$, there exist constants $C_0$,~$n_0>0$ (depending on $\varepsilon$, $k$, $N$ and~$\Graphseq$) such that for any $n\geq n_0$ there exists a collection\/  $\mathcal{T}$ of $k$-colouring templates for $G_n$ satisfying:
	\begin{enumerate}[(i)]
		\item $\mathcal{T}$ is a container family for\/ $\mathcal{P}_n$;
		\item for each template $t\in\mathcal{T}$, there are at most $\varepsilon \binom{G_n}{G_N}$ embeddings $A$ of $G_N$ such that $c\leq t_{\vert A}$ for some $c\in \mathcal{F}$;
	\item $\vert \mathcal{T}\vert\leq \exp \left( C_0(\log k){e(G_n)}^{1 -\frac{1}{2e(G_N)-1}}\right)$.
\end{enumerate}	
\end{theorem} 
\begin{proof}
Let $\mathcal{P}=\Forb_{\Graphseq}(\mathcal{F})$ be a property defined by a forbidden family $\mathcal{F}$ of $k$-colourings of $G_N$. The proof of our general container result shall closely follow that of Theorem~\ref{theorem: multi-colour container}.

First we modify the construction of the hypergraph~$H=H(\mathcal{F}, n)$ in the proof of Theorem~\ref{theorem: multi-colour container} as follows:
\begin{itemize}
	\item we set $r=e(G_N)$ (rather than $\binom{N}{2}$);
	\item we let $V(H)= E(G_n)\times [k]$ (rather than $E(K_n)\times [k]$);
\item for every order-preserving embedding $\phi: \ G_N\rightarrow G_n$, and every colouring $c\in \mathcal{F}$, we add to $E(H)$ the $r$-edge
	\[f_{c, \phi}=\bigl\{\bigl(\phi(e), c(e)\bigr) \, : \, e\in E(G_N)\bigr\}.\] 
\end{itemize}
As before, we can easily bound $e(H)$:
\begin{equation*}\label{equation: bound on e(H), general graph sequence}
\binom{G_n}{G_N}\leq e(H) \leq k^{e(G_N)} \binom{G_n}{G_N}.
\end{equation*}
The last thing we have to check is that a $\Graphseq$-analogue of our sparsification lemma, Lemma~\ref{lemma: random sparsification}, holds. Our technical assumptions on the graph sequence $\Graphseq$ (its `goodness') are exactly what is needed for the proof to go through as before.

Fix $\varepsilon_1 \in (0,1)$ and let
\begin{equation}\label{eq: general sparsification probability}
p=\varepsilon_1\binom{G_n}{G_N} \big/\left(12k^{2e(G_N)-2} I(N,n)\right).
\end{equation}
We shall keep each $r$-edge of~$H$ independently with probability~$p$, and delete it otherwise, to obtain a random subgraph~$H'$ of~$H$.
\begin{lemma}\label{lemma: sparsification lemma, general graph sequence}
Let $p$ be as in~\eqref{eq: general sparsification probability} and let $H'$ be the random subgraph of~$H$ defined above.  Consider the following events:
\begin{itemize}
	\item $F_1$ is the event that \[e(H')\geq \frac{p\dbinom{G_n}{G_N}}{2}=\frac{\varepsilon_1}{12k^{2e(G_N)-2}} \frac{\dbinom{G_n}{G_N}^2}{I(N,n)};\] 
	\item $F_2$ is the event that $H'$ has at most $3p^2k^{2e(G_N)-2}I(N,n)=\frac{\varepsilon_1}{4}p\binom{G_n}{G_N}$ pairs of $r$-edges~$(f,f')$ with $\vert f\cap f'\vert \geq 2$;
	\item $F_3$ is the event that for all $S\subseteq V(H)$ with $e(H[S])\geq \varepsilon_1 e(H)$, we have $e(H'[S])\geq \frac{\varepsilon_1}{2} e(H')$.
\end{itemize}
There exists $n_1\in\N$ such that for all $n\geq n_1$, $F_1 \cap F_2 \cap F_3$ occurs with strictly positive probability.
\end{lemma}
\begin{proof}
Simply follow the proof of Lemma~\ref{lemma: random sparsification}. Property~(ii) of a good graph sequence ensures $F_1$ holds with probability~$1-\exp(-\Omega(p\binom{G_n}{G_N}))=1-o(1)$. By Markov's inequality, $F_2$ holds with probability at least~$2/3$. Finally property~(iii) of a good graph sequence ensures $F_3$ holds with probability at least~$1-2^{ke(G_n)}\exp(-\Omega(p\binom{G_n}{G_N}))=1-o(1)$, so that $F_1$, $F_2$ and~$F_3$ hold simultaneously with probability at least~$2/3-o(1)$ which is strictly positive for $n$~sufficiently large.
\end{proof}
Having obtained the sparsification lemma, we apply the container theorem of Saxton--Thomason for linear $r$-graphs (Theorem~\ref{theorem: saxton thomason simple containers}) and finish the proof in exactly the same way as in Theorem~\ref{theorem: multi-colour container}.
\end{proof}
Theorem~\ref{theorem: general multicolour container result for graph sequences} gives us container theorems for hereditary properties of a wide variety of graph sequences $\Graphseq$. To obtain the standard applications of containers, we need two more ingredients, namely (a) the existence of the entropy density function for $\Graphseq$ (i.e.\ an analogue of Proposition~\ref{proposition: entropy density}) and (b) a supersaturation theorem for $\Graphseq$ (i.e.\ an analogue of Lemma~\ref{lemma: supersaturation}).

These ingredients are obtained on a more ad hoc basis than the general container theorem, Theorem~\ref{theorem: general multicolour container result for graph sequences} --- the proofs have to be tailored to $\Graphseq$ to a greater extent ---  though in many cases the same arguments as those we used in Section~\ref{subsection: entropy density, supersaturation} will work with only trivial modifications. Provided we can obtain them, we have as an immediate corollary of our container theorem the following theorem:
\begin{theorem}\label{theorem: counting for graph sequences}
		Let $\Graphseq$ be a good graph sequence and let $k$,~$N \in \N$. Let $\mathcal{F}$ be a collection of forbidden $k$-colourings of~$G_N$ and let\/ $\mathcal{P}$ be the order-hereditary property of $k$-colourings of~$\Graphseq$ of not containing an embedding of a colouring in $\mathcal{F}$. Suppose that the following hold:
		\begin{enumerate}
			\item $\pi(\mathcal{P}):=\lim_{n\rightarrow \infty}\ex(G_n, \mathcal{{P}})/e(G_n)$ exists;
			\item for all $\varepsilon>0$ there exist $\delta$,~$n_0>0$ such that if $n\geq n_0$ and $t$ is a $k$-colouring template for $E(G_n)$ with at most $\varepsilon \binom{G_n}{G_N}$ pairs~$(\phi, c)$ where $\phi:G_N\rightarrow G_n$ is an embedding and $c\in\mathcal{F}$ is a forbidden colouring of~$G_N$ such that  
			$c(e)\in t(\phi(e))$ for all $e\in E(G_n)$, then $\Ent(t)\leq (\pi(\mathcal{P})+\delta)e(G_n)$.
		\end{enumerate}	
		Then
		\[\vert \mathcal{P}_n\vert = k^{\bigl(\pi(\mathcal{P})+o(1)\bigr)e(G_n)}.\]
\end{theorem}
\begin{proof}
This is identical to the deduction of Theorem~\ref{theorem: counting result for Forb(F), F finite hereditary families} from Theorem~\ref{theorem: multi-colour container}, Proposition~\ref{proposition: entropy density} and Lemma~\ref{lemma: supersaturation}.
\end{proof}

As an illustration of the way we can mimic the proofs of Proposition~\ref{proposition: entropy density} and Lemma~\ref{lemma: supersaturation} to show conditions 1.\ and~2.\ in Theorem~\ref{theorem: counting for graph sequences} are satisfied, we include below a proof of those two results in the case where $\Graphseq$ is the sequence of hypercube graphs~$\Hypercube$.
\begin{proposition}[Goodness of hypercube graphs]\label{proposition: hypercube graphs are good}
The sequence $\Hypercube$ is good.
\end{proposition}
\begin{proof}
	We have $e(Q_n)=n2^{n-1}$ and $\binom{Q_n}{Q_N}=\binom{n}{N}2^{n-N}=\Omega(2^n n^N)$, establishing parts (i) and~(ii) of Definition~\ref{definition: good graph sequence}. For part~(iii), noting that two $N$-dimensional subcubes with at least two edges in common must meet in an $i$-dimensional subcube for some $i$: $2\leq i\leq N$, we have
	\begin{align*}
	I(N,n)&=\frac{1}{2}\binom{Q_n}{Q_N}\sum_{2\leq i\leq N} \binom{N}{i}2^{N-i}\binom{n-N}{N-i}=O\biggl(\binom{Q_n}{Q_N} n^{N-2}\biggr),
	\end{align*}
	which gives us $I(N,n)e(Q_n)/{\binom{Q_n}{Q_N}}^2=O(1/n)=o(1)$ as required. 
\end{proof}

Now we show that the sequence~$\Hypercube$ satisfies the hypotheses of Theorem~\ref{theorem: counting for graph sequences}.

\begin{proposition}[Entropy density for edge-colourings of hypercubes]\label{proposition: entropy density edge-hypercube}
If\/ $\mathcal{P}$ is an order-hereditary property of $k$-colourings of~$\Hypercube$, then the limit $\pi^{\Hypercube}(\mathcal{P}):=\lim_{n\rightarrow \infty}\ex(Q_n, \mathcal{P})/2^{n-1}n$ exists.	
\end{proposition}
\begin{proof}
Let $t$ be an extremal entropy template for $\mathcal{P}$ in $Q_{n+1}$. By averaging over all embeddings $\phi$ of $Q_n$ into $Q_{n+1}$, we have
\begin{equation*}
n\ex(Q_{n+1}, \mathcal{P}) =(n-1)\Ent(t)=\sum_{\phi} \Ent\bigl(t_{\vert \phi(Q_n)}\bigr)\leq 2(n+1) \ex(Q_n, \mathcal{P}), 
\end{equation*}
whence $\ex(Q_n, \mathcal{P})/(2^{n-1}n)$ is nonincreasing in $[0,1]$, and hence tends to a limit as $n\rightarrow \infty$.
\end{proof}
\begin{proposition}[Supersaturation for edge-colourings of hypercube]\label{proposition: supersaturation for hypercube graphs}
Let $N \in \N$ be fixed and let $\mathcal{F}$ be a nonempty collection of $k$-colourings of~$Q_N$. Set\/ $\mathcal{P}=\Forb_{\Hypercube}(\mathcal{F})$. For every $\varepsilon \in (0, 1)$, there exist constants $n_0\in\N$ and~$C_{0}>0$ such that for all $n\geq n_0$, if $t$ is a $k$-colouring template for $Q_n$ such that $\Ent(t)> (\pi(\mathcal{P})+\varepsilon)\binom{n}{2}$, then there are at least~$C_{0} \varepsilon\binom{Q_n}{Q_N}$ pairs~$(\phi, c_i)$ where $\phi: \ Q_N \rightarrow Q_n$ is an embedding and $c_i \in \mathcal{F}$ is a colouring with $c_i(e)\in t(\phi(e))$ for every $e\in E(Q_N)$.
\end{proposition}
\begin{proof}
We follow the proof of Lemma~\ref{lemma: supersaturation}, modifying it as needed to fit the hypercube setting. Instead of taking $X$ to be a random $m$-subset of $V(K_n)$, we take $X$ to be a random $M$-dimensional subcube of~$Q_n$, where $M\sim \Binom(n,p)$ (and the remaining~$n-M$ coordinates' values are chosen uniformly at random). The key difference lies in the bootstrap bound we give for $B(t)$. 
Explicitly, observe that the probability that $\vector{e}\in E(Q_n)$ is included in $X$ is exactly~$p\left(\frac{1+p}{2}\right)^{n-1}$, and that the probability that a fixed copy of~$Q_N$ is included as a subcube of~$X$ is~$p^N \left(\frac{1+p}{2}\right)^{n-N}$.

As in Lemma~\ref{lemma: supersaturation}, we use the monotonicity established in Proposition~\ref{proposition: entropy density edge-hypercube} to find a constant~$n_1$ such that for all $n\geq n_1$ we have $\ex(Q_n, \mathcal{P})\leq \left(\pi^{\Hypercube}(\mathcal{P}) +\frac{\varepsilon}{3}\right)n2^{n-1}$, and we let $p=\frac{8n_2}{n}$ for some large constant $n_2\geq n_1$. The hypercube setting is then slightly easier to work with than the graph setting of Lemma~\ref{lemma: supersaturation}, since we do not have to compute second moments, and can bound $B(t)$ directly:
\begin{align*}
p^N\left(\frac{1+p}{2}\right)^{n-N}B(t) &
= \mathbb{E} B\bigl(t_{\vert X}\bigr)\\
&\geq \frac{1}{\log_k (2)} \bigl(1-\mathbb{P}(\mathcal{A})\bigr)p\left(\frac{1+p}{2}\right)^{n-1} \left(\Ent(t) -\left(\pi_v^{\Hypercube}(\mathcal{P})+\frac{\varepsilon}{3}\right)n2^{n-1}\right).
\end{align*}
Thus if $\Ent(t)>(\pi_v^{\Hypercube}(\mathcal{P})+\varepsilon)n2^{n-1}$, \eqref{eq: chernoff bound on A} shows that
\begin{equation*}
B(t)> \frac{\varepsilon}{3\log_k(2)}\left(\frac{1+p}{2p}\right)^{N-1} n2^{n-1}> \frac{\varepsilon}{3\log_k(2)} \frac{N!}{(16n_2)^N}\binom{Q_n}{Q_N},
\end{equation*}
as required.
\end{proof}
\begin{corollary}[Counting for hypercube graph colourings]\label{corollary: edge hypercube counting}
If\/ $\mathcal{P}$ is an order-hereditary property of $k$-colourings of~$\Hypercube$, then $\vert \mathcal{P}_n\vert=k^{\left(\pi^{\Hypercube}(\mathcal{P})+o(1)\right)2^{n-1}n}$.
\end{corollary}
\begin{proof}
Propositions \ref{proposition: hypercube graphs are good}, \ref{proposition: entropy density edge-hypercube} and~\ref{proposition: supersaturation for hypercube graphs} tell us that the hypotheses of Theorem~\ref{theorem: counting for graph sequences} are satisfied; applying it yields the desired counting result.
\end{proof}

\subsection{Hypergraphs}\label{subsection: hypergraphs}
It is trivial to generalise Theorem~\ref{theorem: multi-colour container} to the $l$-graph setting, for any $l\geq 1$: instead of $k$-colouring the edges of the complete $2$-graph $K_n$, we could consider $k$-colourings of~$K_n^{(l)}$, the complete $l$-graph on $n$ vertices. In the proof of Theorem~\ref{theorem: multi-colour container}, instead of setting $V(H)= E(K_n^{(2)})$, we would set $V(H)=E(K_n^{(l)})$, and proceed onwards as before, setting $r=\binom{N}{l}$ and modifying constants as needed. More generally, the proof of Theorem~\ref{theorem: general multicolour container result for graph sequences} carries over to the setting of $l$-graph sequences~$\Graphseq$ without any change. 
\begin{theorem}\label{theorem: general multicolour container result for hypergraph sequences}
	Let $l\in \N$ and let $\Graphseq$ be a good $l$-graph sequence. Let $k$,~$N \in \N$. Let $\mathcal{F}$ be a collection of forbidden $k$-colourings of $G_N$, and let\/ $\mathcal{P}$ be the order-hereditary property of $k$-colourings of~$\Graphseq$ of not containing an embedding of a colouring in $\mathcal{F}$.

	For any $\varepsilon>0$, there exist constants $C_0$,~$n_0>0$ (depending on $\varepsilon$, $k$, $N$ and~$\Graphseq$) such that for any $n\geq n_0$ there exists a collection~$\mathcal{T}$ of $k$-colouring templates for $G_n$ satisfying:
	\begin{enumerate}[(i)]
		\item $\mathcal{T}$ is a container family for\/ $\mathcal{P}_n$;
		\item for each template $t\in\mathcal{T}$, there are at most $\varepsilon \binom{G_n}{G_N}$ embeddings~$A$ of~$G_N$ such that $c\leq t_{\vert A}$ for some $c\in \mathcal{F}$;
		\item $\vert \mathcal{T}\vert\leq \exp \left( C_0(\log k){e(G_n)}^{1 -\frac{1}{2e(G_N)-1}}\right)$.
	\end{enumerate}	
\end{theorem} 
\begin{proof}
	Identical to the proof of Theorem~\ref{theorem: general multicolour container result for graph sequences}: notice that we nowhere used the fact that $\Graphseq$ was a sequence of graphs! In fact, the only difference the $l$-graph setting makes is that in our definitions of good graph sequences the condition $N\geq 3$ must be replaced by $N\geq l+1$.
\end{proof}
There are many applications in which one is interested in \emph{vertex}-colouring of $l$-graph sequences. Again, the proof of Theorem~\ref{theorem: general multicolour container result for graph sequences} carries over to this setting with almost no change except the definition of goodness.
\begin{definition}[Good hypergraph sequence for vertex colouring]\label{def: hypergraph vertex-good}
	An $l$-uniform graph sequence~$\Graphseq$ is vertex-\emph{good} if all of the following hold:
	\begin{enumerate}[(i)]
		\item $v(G_n)\rightarrow \infty$ (`the graphs in the sequence become large');
		\item for all $N \in \N$ with $N\geq l$, $\binom{G_n}{G_N}\rightarrow \infty$ (`the sequence has many embeddings of $G_N$');
		\item for all $N \in \N$ with $N\geq 3$, $v(G_n)J(N,n)\big/ \binom{G_n}{G_N}^2 \rightarrow 0$ as $n\rightarrow \infty$ (`most pairs of embeddings of~$G_N$ are almost disjoint'), where $J(N,n)$ counts the number of joint embeddings of~$G_N$ into $G_n$ with at least~$2$ vertices in common.
	\end{enumerate}	
\end{definition}

\begin{theorem}\label{theorem: general multicolour container result for vertex colouring of hypergraph sequences}
	Let $l\in \N$ and let $\Graphseq$ be a vertex-good $l$-graph sequence. Let $k$,~$N \in \N$. Let $\mathcal{F}$ be a collection of forbidden $k$-colourings of~$V(G_N)$ and let\/ $\mathcal{P}$ be the order-hereditary property of $k$-colourings of vertices of~$\Graphseq$ of not containing an embedding of a colouring in $\mathcal{F}$.

For any $\varepsilon>0$, there exist constants $C_0$,~$n_0>0$ (depending on $\varepsilon$, $k$, $N$ and~$\Graphseq$) such that for any $n\geq n_0$ there exists a collection\/ $\mathcal{T}$ of $k$-colouring templates for $G_n$ satisfying:
	\begin{enumerate}[(i)]
		\item $\mathcal{T}$ is a container family for\/ $\mathcal{P}_n$;
		\item for each template $t\in\mathcal{T}$, there are at most $\varepsilon \binom{G_n}{G_N}$ embeddings~$A$ of $G_N$ such that $c\leq t_{\vert A}$ for some $c\in \mathcal{F}$;
		\item $\vert \mathcal{T}\vert\leq \exp \left( C_0(\log k){v(G_n)}^{1 -\frac{1}{2v(G_N)-1}}\right)$.
	\end{enumerate}	
\end{theorem} 
\begin{proof}
The only change we need to make in the proof of Theorem~\ref{theorem: general multicolour container result for graph sequences} is that we take $V(G_n)\times [k]$ as the vertex set of our hypergraph~$H$ rather than $E(G_n)\times [k]$.
\end{proof}
As an immediate corollary to Theorem~\ref{theorem: general multicolour container result for vertex colouring of hypergraph sequences}, we obtain a container result for properties of $k$-colourings of the vertices of~$Q_n$ (which can be viewed as $k$-colourings of the vertices of a $2$-graph sequence).  In the next section, we shall state this theorem and use it to deduce counting results.

\subsection{Vertex-colourings of the hypercube}\label{subsection: vertex colouring hypercube}
In this subsection, we extend our results on $k$-colourings of the edges of~$K_n$ to $k$-colourings of the vertices of~$Q_n$. 



	For $m\leq n$, an \emph{order-preserving embedding} of~$Q_m$ into $Q_n$ is a map $\phi: \ Q_m \rightarrow Q_n$ such that there exists an $m$-set~$B=\{b_1, b_2, \ldots, b_m\}\subseteq [n]$ with
	\[\phi(\vector{x})_i= \left\{\begin{array}{ll}
	\phi(\vector{0})_i & \textrm{if }i\in [n]\setminus B\\
	\vector{x}_j & \textrm{if }i=b_j\in B.\end{array}
	\right.\]
	Conversely, the \emph{subcube}~$Q_n[(B, \vector{v})]$ of~$Q_n$ induced by a set~$B=\{b_1,\ldots, b_m\}\subseteq [n]$ and a vector~$\vector{v}$ is the order-preserving embedding  of~$Q_m$ into $Q_n$ defined by
		\[\phi(\vector{x})_i= \left\{\begin{array}{ll}
		\phi(\vector{v})_i & \textrm{if }i\in [n]\setminus B\\
		\vector{x}_j & \textrm{if }i=b_j\in B.\end{array}
		\right.\]
\begin{definition}[Hypercube properties]
Let $k\in \N$. A $k$-colouring \emph{vertex property} of hypercubes~$\mathcal{P}=\left(\mathcal{P}_n\right)_{n\in\N}$ is a sequence of families $\mathcal{P}_n$ of $k$-colourings $c: \ Q_n \rightarrow [k]$ of $Q_n$. 
A vertex property of hypercubes $\mathcal{P}$ is \emph{order-hereditary} if for every $m$-set~$B\subseteq [n]$, vector~$\vector{v}\in Q_n$ and colouring $c\in \mathcal{P}_n$ the colouring $c_{\vert (B, \vector{v})}$ of the subcube $Q_n[(B, \mathbf{v})]$ induced by $c$ lies in $\mathcal{P}_m$.
\end{definition}
Applying the hypergraph sequence result proved in the previous section (Theorem~\ref{theorem: general multicolour container result for vertex colouring of hypergraph sequences}), we obtain container and counting results for vertex-properties of $k$-colourings of hypercubes.
\begin{theorem}\label{theorem: vertex hypercube containers}
	Let $\mathcal{F}$ be a nonempty family of vertex $k$-colourings of~$Q_N$ and let\/ $\mathcal{P}=\Forb_{\Hypercube}(\mathcal{F})$. For any $\varepsilon>0$, there exists $n_0=n_0(\varepsilon, k, N)$ such that for any $n\geq n_0$ there exists a collection\/ $\mathcal{T}$ of vertex $k$-colouring templates for $Q_n$ satisfying:
	\begin{enumerate}[(i)]
		\item $\mathcal{T}$ is a container family for\/ $\mathcal{P}_n$;
		\item for each template $t\in\mathcal{T}$, there are at most $\varepsilon \binom{n}{N}2^{n-N}$ pairs $(\phi,c)$ where $\phi$ is an order-preserving embedding of $Q_N$ into $Q_n$ and $c$ is a forbidden colouring $c\in \mathcal{F}$ with $c(\mathbf{x})\in t(\phi(\mathbf{x}))$ for every $\mathbf{x}\in Q_N$;   
\item $\log_k \lvert \mathcal{T} \rvert \leq 2^{n \bigl(1-\frac{1}{2^{N+1}-1}+\varepsilon\bigr)}$.
	\end{enumerate}
\end{theorem}
\begin{proof}
	Apply Theorem~\ref{theorem: general multicolour container result for vertex colouring of hypergraph sequences}. All we need to check is that $\Hypercube$ is indeed a good $1$-graph sequence.  Parts (i) and~(ii) of Definition~\ref{def: hypergraph vertex-good} are obvious, while part~(iii) is a simple calculation:
	\begin{equation*}
	\frac{v(Q_n) J(N,n)}{\binom{Q_n}{Q_N}} =\frac{2^n}{2^{n-N}\binom{n}{N}} \sum_{i>1} \binom{N}{i}\binom{n-N}{N-i}2^{N-i}=O(1/n),
	\end{equation*}
as required.
\end{proof}
\begin{proposition}\label{proposition: entropy density, vertex hypercube}
	If\/ $\mathcal{P}$ is an order-hereditary $k$-colouring vertex property of hypercubes, then the limit 
	\[\pi_v^{\Hypercube}(\mathcal{P}):=\lim_{n\rightarrow \infty} \frac{\ex(Q_n, \mathcal{P})}{2^n} \]
	exists.
\end{proposition}
We call $\pi_v^{\Hypercube}(\mathcal{P})$ the \emph{entropy density} of $\mathcal{P}$.
\begin{proof}
Let $t$ be a vertex $k$-colouring template for $Q_{n+1}$ with extremal entropy relative to $\mathcal{P}$. By averaging over the $2(n+1)$ distinct $n$-dimensional induced subcubes $Q_{n+1}[(B, \vector{x})]$ with $\vector{x}\in\{\vector{0}, \vector{1}\}$, we have
\begin{align*}
	(n+1)\ex(Q_{n+1}, \mathcal{P})=(n+1)\Ent(t)
	= \sum_{\vector{x}\in\{\vector{0}, \vector{1}\}}\sum_{B\in [n+1]^{(n)}}\Ent\bigl(t_{\vert Q_{n+1} [(B,\vector{x})]}\bigr)\leq 2(n+1) \ex(Q_n, \mathcal{P}), 
	\end{align*}
	whence $\ex(Q_{n}, \mathcal{P})/2^{n}$ is non-increasing in $[0,1]$ and converges to a limit as required.
\end{proof}
\begin{lemma}[Supersaturation]\label{lemma: vertex hypercube supersaturation}
	Let $\mathcal{F}$ be a nonempty family of vertex $k$-colourings of $Q_N$ and let\/ $\mathcal{P}=\Forb_{\Hypercube}(\mathcal{F})$.  For every $\varepsilon$ with $0<\varepsilon<1$, there exist constants $n_0\in\N$ and~$C_{0}>0$ such that for any vertex $k$-colouring template~$t$ for $Q_n$ with $n\geq n_0$ and $\Ent(t)> \left(\pi_v^{\Hypercube}(\mathcal{P})+\varepsilon\right)2^n$, there are at least~$C_{0} \varepsilon\binom{n}{N}2^{n-N}$ pairs~$(\phi, c)$, where $\phi$ is an order-preserving embedding of~$Q_N$ into $Q_n$ and $c\in \mathcal{F}$ satisfies $c(\vector{x})\in t(\phi(\vector{x}))$ for every $\vector{x}\in Q_N$.
\end{lemma}
\begin{proof}
	We follow the proof of Lemma~\ref{lemma: supersaturation}, modifying it as needed to fit the vertex-hypercube setting. Instead of taking $X$ to be a random $m$-subset of $V(K_n)$, we take $X$ to be a random $M$-dimensional subcube of $Q_n$, where $M\sim \Binom(n,p)$ (and the remaining $n-M$ coordinates' values are chosen uniformly at random). The key difference lies in the bootstrap bound we give for $B(t)$. 
Explicitly, observe that the probability $\mathbf{x}\in Q_n$ is included in $X$ is exactly $\left(\frac{1+p}{2}\right)^n$, and that the probability a fixed copy of $Q_N$ is included as a subcube of $X$ is $p^N \left(\frac{1+p}{2}\right)^{n-N}$.

As in Lemma~\ref{lemma: supersaturation}, we use the monotonicity established in Proposition~\ref{proposition: entropy density, vertex hypercube} to find a constant~$n_1$ such that for all $n\geq n_1$ we have $\ex(Q_n, \mathcal{P})\leq \left(\pi^{\Hypercube}(\mathcal{P}) +\frac{\varepsilon}{3}\right)2^n$, and we let $p=\frac{8n_2}{n}$ for some large constant~$n_2\geq n_1$. The vertex hypercube setting is then slightly easier to work with than the graph setting of Lemma~\ref{lemma: supersaturation}, since we do not have to compute second moments, and can bound $B(t)$ directly:
\begin{equation*}
p^N\left(\frac{1+p}{2}\right)^{n-N}B(t) = \mathbb{E} B\bigl(t_{\vert X}\bigr) \geq \frac{1}{\log_k (2)} \bigl(1-\mathbb{P}(\mathcal{A})\bigr)\left(\frac{1+p}{2}\right)^n \left(\Ent(t) -\left(\pi_v^{\Hypercube}(\mathcal{P})+\frac{\varepsilon}{3}\right)2^n\right).
\end{equation*}
Thus if $\Ent(t)>(\pi_v^{\Hypercube}(\mathcal{P})+\varepsilon)2^n$, we have
\begin{align*}
B(t)> \frac{\varepsilon}{3\log_k(2)}\left(\frac{1+p}{2p}\right)^N 2^n> \frac{\varepsilon}{3\log_k(2)}\frac{N!}{(8n_2)^N}\binom{n}{N}2^{n-N},
\end{align*}
as required.
\end{proof}
From there, the counting result is immediate:
\begin{corollary}\label{corollary: vertex hypercube counting}
If\/ $\mathcal{P}$ is an order-hereditary property of vertex $k$-colourings of $Q_n$, then 
\[\vert \mathcal{P}_n\vert = k^{\left(\pi_v^{\Hypercube}(\mathcal{P})+o(1)\right)2^n}.\]	
\end{corollary}
\begin{proof}
	Use Theorem~\ref{theorem: vertex hypercube containers} for the $1$-graph sequence $\Hypercube$ and Lemma~\ref{lemma: vertex hypercube supersaturation} to obtain the vertex-hypercube analogue of Theorem~\ref{theorem: counting result for Forb(F), F finite hereditary families}. Then establish a vertex-hypercube analogue of Theorem~\ref{theorem: approximation of general hereditary properties}, substituting the monotonicity of hypercube entropy that was proved in Proposition~\ref{proposition: entropy density, vertex hypercube} for Proposition~\ref{proposition: entropy density}. Finally, use this to derive Corollary~\ref{corollary: vertex hypercube counting} similarly to the deduction of Corollary~\ref{corollary: speed of arbitrary hereditary properties}.
\end{proof}




	\section{Examples and applications}\label{section: examples}
	\subsection{Order-hereditary versus hereditary}
	Here we include a quick example stressing the essential difference between hereditary and order-hereditary properties. We identify graphs with $\{0,1\}$-colourings of~$K_n$ in the usual way; as the properties we shall consider are in fact monotone, templates will consist of pairs~$e\in E(K_n)$ with $t(e)=\{0,1\}$ and entropy~$1$ and pairs~$e$ with $t(e)=\{0\}$ and entropy~$0$. We can thus represent the templates simply as the graph of edges with entropy~$1$.

	Let $\mathcal{P}_1$ be the hereditary property of graphs on $[n]$ of having maximum degree~$2$, and let $\mathcal{P}_2$ be the order-hereditary property of graphs on $[n]$ of not having any triples of vertices~$i<j<k$ with $ij$,~$jk$ both being edges.

	 It is trivial to show that $\ex(n, \mathcal{P}_1)=\lfloor n/2\rfloor$, with maximal matchings being the extremal entropy templates.  By counting the number of matchings on $[n]$, it follows that 
	 \[\bigl\lvert (\mathcal{P}_1)_n\bigr\rvert  = \left(\frac{1+\sqrt{5}}{2}\right)^{n+o(n)}.\]
	 On the other hand, $\ex(n, \mathcal{P}_2)=\lfloor \frac{n^2}{2}\rfloor$. For the lower bound, consider the template $t$ whose entropy~$1$ edges are $\{(i,j): \ 1\leq i \leq \frac{n}{2}\leq j\leq n, \ i\neq j\}$. Clearly $\Ent(t)=\lfloor n^2/4\rfloor$ and we have no $i<j<k$ with $ij$,~$jk$ both being edges. For the upper bound, suppose $\Ent(t)>n^2/4$. By Mantel's theorem, there must exist a triangle~$ijk$ of edges with full entropy, which gives us a triple of vertices~$i<j<k$ with $ij$,~$jk$ both being edges. Applying Theorem~\ref{theorem: counting result for Forb(F), F finite hereditary families}, we have that
	  \[\bigl\lvert (\mathcal{P}_2)_n\bigr\rvert  = 2^{\frac{1}{4}n^{2}\bigl(1+o(1)\bigr)}.\]

	\subsection{Graphs}\label{subsection: proof of alekseev--bollobas--thomason}

In this section, we give a new, short proof of the Alekseev--Bollob\'as--Thomason theorem (Theorem~\ref{theorem: alekseevbollobasthomason}).  Our argument is similar to the proof in~\cite{BT}.  However, our entropy results mean that our proof does not require us to use any form of the Regularity Lemma, which simplifies the argument considerably.

We shall need to use the Erd\H{o}s--Stone theorem~\cite{ErdosStone46}.

\begin{theorem}[Erd\H{o}s--Stone theorem]\label{theorem: erdosstone}
Let $r \geq 2$, $m \geq 1$ and $\varepsilon > 0$.  There exists $n_o(r, m, \varepsilon)$ such that if $G$ is a graph of order~$n \geq n_0$ and
\[
e(G) \geq \biggl(1 - \dfrac{1}{r} + \varepsilon \biggr)\binom{n}{2},
\]
then $G$ contains a copy of~$K_{r+1}(m)$.
\end{theorem}

Recall the definition of~$\mathcal{H}(r, \vector{v})$ from Definition~\ref{definition: colouring number} and observe that $\mathcal{H}(r, \vector{v})$ is a symmetric hereditary class.

\begin{lemma}\label{lemma: universal}
Let\/ $\PP$ be a symmetric hereditary property of graphs, let $r \geq 2$, let $\ell \geq 1$ and let $\varepsilon > 0$.  There exists a constant~$n_0 \in \N$ depending only on $r$, $\ell$ and~$\varepsilon$ such that if $n \geq n_0$ and $t$ is a $2$-colouring template for $K_n$ with $\langle t \rangle \subseteq \PP_n$ and
\[
\Ent(t) \geq \biggl(1 - \dfrac{1}{r} + \varepsilon \biggr)\binom{n}{2},
\]
then $\mathcal{H}(r + 1, \vector{v})_{\ell} \subseteq \PP$ for some $\vector{v} \in \{0, 1\}^{r + 1}$.
\end{lemma}

\begin{proof}
By Ramsey's theorem, for each $\ell$, there exists $m$ such that any $2$-colouring of~$E(K_m)$ contains a monochromatic copy of~$K_{\ell}$.  Let $G$ be the graph with vertex set~$[n]$ and $E(G) = \{e \in E(K_n) : \ t(e) = \{0, 1\}\}$.  Our assumption on $\Ent(t)$ and the Erd\H{o}s--Stone theorem imply that if $n$ is sufficiently large, then $G$ contains a copy~$K$ of~$K_{r + 1}(m)$.

Let $t'$ denote the restriction of~$t$ to $V(K)$ and let $V_1$, \dots,~$V_{r + 1}$ denote the classes of~$V(K)$.  Now we construct a vector~$\vector{v} \in \{0, 1\}^{r + 1}$.  By choice of~$m$, for each $i$, either $\{e \in E(K[V_i]) : 0 \in t(e)\}$ or $\{e \in E(K[V_i]) : 1 \in t(e)\}$ contains a copy of~$K_{\ell}$.  If the first case holds, then we set $v_i = 0$, and otherwise we set $v_i = 1$.  In either case, we let $U_i$ denote the vertex set of the copy of~$K_{\ell}$.  Let $H \in \mathcal{H}(r + 1, \vector{v})_{\ell}$ and let $W_1$, \dots~$W_{r + 1}$ be a partition of~$V(H)$ such that for each $i$, $W_i$ is a clique if $v_i = 1$ and an independent set if $v_i = 0$.  Because $\abs{V(H)} = \ell$, we may embed each $W_i$ into $U_i \subseteq V_i$ arbitrarily.  It follows that there is a realisation~$c$ of~$t'$ such that $H$ is a subgraph of~$c(K)$.

Finally, because $\PP$ is hereditary, it follows that $\mathcal{H}(r + 1, \vector{v})_{\ell} \subseteq \langle t' \rangle \subseteq \PP$, which is what we wanted.
\end{proof}

\begin{proof}[Proof of Theorem~\ref{theorem: alekseevbollobasthomason}.]
First, by the definition of~$\chi_c(\PP)$, there exists $\vector{v} \in \{0, 1\}^r$ such that $\mathcal{H}(r, \vector{v}) \subseteq \PP$.  By considering the graphs in $\mathcal{H}(r, \vector{v})$ such that each clique or independent set has size $\lfloor n/r \rfloor$ or~$\lceil n/r \rceil$, we see that
\[
\lvert \PP_n \rvert \geq \lvert \mathcal{H}(r, \vector{v})_n \rvert \geq 2^{\bigl(1 - 1/r + o(1)\bigr)\binom{n}{2}}.
\]

Second, suppose for a contradiction that for some $\varepsilon > 0$, there exist infinitely many $n$ such that
\begin{equation}\label{eq: P_n too big}
\lvert \PP_n \rvert \geq 2^{(1 - 1/r + \varepsilon)\binom{n}{2}}.
\end{equation}
Corollary~\ref{corollary: speed of arbitrary hereditary properties} implies that there exists $n_0$ such that for all $n \geq n_0$ such that~\eqref{eq: P_n too big} holds, there exists a template~$t$ for $2$-colourings of~$K_n$ such that $\langle t \rangle \subseteq \PP_n$ and
\[
\Ent(t) \geq \biggl(1 - \dfrac{1}{r} + \dfrac{\varepsilon}{2} \biggr)\binom{n}{2}.
\]

It follows from Lemma~\ref{lemma: universal} that for each $\ell \geq 1$, there exists $\vector{v} \in \{0, 1\}^{r+1}$ such that $\mathcal{H}(r + 1, \vector{v})_{\ell} \subseteq \PP$.  In particular, there is some $\vector{v} \in \{0, 1\}^{r+1}$ such that $\mathcal{H}(r + 1, \vector{v})_{\ell} \subseteq \PP$ for infinitely many $\ell$, and thus for all $\ell$.  However, this contradicts the definition of~$\chi_c(\PP)$, and the desired result follows
\end{proof}

	\subsection{Digraphs}\label{subsection: digraph examples}

	As mentioned earlier, hereditary properties for tournaments, orgraphs and digraphs have received significant attention from the extremal combinatorics research community, see~\cite{Bollobas07}. In a recent paper, K\"uhn, Osthus, Townsend, Zhao~\cite{KuhnOsthusTownsendZhao14} determined the typical structure of certain families of oriented and directed graphs.  In doing so, they proved a container theorem and, using it, a counting theorem for $H$-free orgraphs and $H$-free digraphs, where $H$ is a fixed orgraph with at least~two edges (Theorems~3.3 and Corollary~3.4 in~\cite{KuhnOsthusTownsendZhao14}).

	They went on to observe that their results did not extend to the case where $H$ is a digraph, giving the specific example when $H=DK_3$, the double triangle $([3], [3]\times[3])$. Their approach considered the extremal weight achievable in an $H$-free digraph where double edges receive a different weight from single edges. In the case of $DK_3$, they observed that the extremal weight did not predict the correct count of $DK_3$-free digraphs, showing that their container theorem failed to generalise in its given form to the digraph case. Giving some vindication to our entropy-based approach to containers, we use our theorems to determine the speed of the digraph property~$\mathcal{P}$ of not containing any $DK_3$. 
	
	\begin{theorem}\label{theorem: extremal entropy, no double triangle}  If\/ $\PP$ is the digraph property of not containing $DK_3$, then
\[\ex(n, \mathcal{P})= (1-\log_4 3) \biggl\lfloor \frac{n^2}{4}\biggr\rfloor + \log_43 \binom{n}{2}.\]
	\end{theorem}
	\begin{proof}
	We use the correspondence between digraphs and $4$-colourings of~$K_n$ outlined above. Let $t$ be an $n$-vertex $4$-colouring template for $\mathcal{P}$ with maximal entropy. The monotonicity of $\mathcal{P}$ and the maximality of $\Ent(t)$ imply that all edges $e$ of~$K_n$ have $t(e)=[4]$ or $t(e)=[3]$. As $\mathcal{P}$ is exactly the property of having no triangle in colour $4$, Mantel's theorem tells us that at most $\lfloor \frac{n^2}{4}\rfloor$ edges can have full entropy (entropy $1$), with the rest having entropy $\log_4 3$. Thus 
	\[\Ent(t)\leq (1-\log_4 3) \biggl\lfloor \frac{n^2}{4}\biggr\rfloor + \log_43 \binom{n}{2},\]
	as required.

	For the lower bound, consider a balanced bipartition of $[n]$ as $A\sqcup B$. Let $t$ be the $n$-vertex $4$-colouring template with $t(e)=[4]$ if $e$ is an edge from $A$ to $B$, and $t(e)=[3]$ otherwise. Clearly every realisation of $t$ contains no triangle in colour $4$, and hence lies in $\mathcal{P}$, and the entropy of $t$ exactly matches the upper bound we have established above.
	\end{proof}
	\begin{corollary}
	There are $4^{\pi(\mathcal{P})\binom{n}{2}+o(n^2)}= 3^{\binom{n}{2}-\bigl\lfloor \frac{n^2}{4}\bigr\rfloor}4^{\bigl\lfloor \frac{n^2}{4}\bigr\rfloor +o(n^2)}$ digraphs on $n$ vertices not containing any $DK_3$.
	\end{corollary}
	\begin{proof}
	Theorem~\ref{theorem: extremal entropy, no double triangle} establishes $\pi(\mathcal{P})= \frac{1}{2}+\frac{\log_43}{2}$. We then apply Corollary~\ref{corollary: speed of arbitrary hereditary properties} and are immediately done.
	\end{proof}
	Furthermore, we can characterise typical graphs in $\mathcal{P}$. Let $\mathcal{S}_n$ denote the collection of $n$-vertex $4$-colouring templates~$t$ obtained by taking a balanced bipartition $A\sqcup B$ of $[n]$ and setting $t(e)=[4]$ for all edges~$e$ from $A$ to $B$ and $t(e)=[3]$ for all edges~$e$ internal to $A$ or $B$.

	The well-known stability theorem for Mantel's theorem~\cite{Erdos67,Simonovits68} immediately implies the following result:
	\begin{proposition}\label{theorem: stability for no double triangle digraph}
	Let\/ $\PP$ denote the digraph property of not containing $DK_3$.  For every $\varepsilon>0$, there exists $\delta>0$ and $n_0$ such that if $n\geq n_0$ and $t$ is an $n$-vertex $4$-colouring template satisfying
	\begin{enumerate}[(i)]
	\item $\Ent(t)\geq \left(\pi(\mathcal{P})-\delta\right)\binom{n}{2}$, and
	\item there are at most $\delta n^3$ triples of vertices $\{a, b, c\}$ which give rise to a monochromatic triangle in colour~$4$ in some realisation of~$t$,
	\end{enumerate} 
	then $\rho(\mathcal{S}_n, t) \leq \varepsilon \binom{n}{2}$.\qed
	\end{proposition}
	\begin{corollary}
	Let\/ $\PP$ denote the digraph property of not containing $DK_3$.
For every $\varepsilon >0$ there exists $n_0>0$ such that for all $n\geq n_0$, all but $\varepsilon \vert \mathcal{P}_n\vert$ colourings in\/ $\mathcal{P}$ are within edit distance~$\varepsilon \binom{n}{2}$ of a realisation from\/ $\mathcal{S}_n$.

Equivalently, for all but an $\varepsilon$-proportion of $DK_3$-free digraphs~$D$, there exists a digraph~$H$ that is obtained by taking a balanced bipartition of the vertex set~$A\sqcup B=[n]$, setting double edges between $A$ and $B$ and letting $A$ and $B$ be quasirandom tournaments and a subdigraph~$H'$ of~$H$ that satisfies $\rho(D, H') \leq \varepsilon\binom{n}{2}$.\qed
	\end{corollary}

\begin{remark}\label{remark: digraphs in general}
Given a graph~$F$, let $DF$ be the digraph obtained by replacing each edge of~$F$ with a directed edge in each direction.  It is easy to see that all of the results of this section extend to the class of $DF$-free digraphs.
\end{remark}
	
	


	

	\subsection{Multigraphs}
	We consider multigraphs, viewed as weightings of the edges of~$K_n$ by non-negative integers. For example, let $\mathcal{P}$ be the property of multigraphs that no triple of vertices supports more than $4$ edges (counting multiplicities). Clearly no edge of such a multigraph can have weight more than $4$.

	We may view an $n$-vertex multigraph~$G$ in which all edge multiplicities are at most~$d$ as $(d + 1)$-colourings of~$E(K_n)$, with each edge coloured by its multiplicity. In this way, the problem of counting such multigraphs is placed in our framework of counting $k$-colourings.  (Mubayi and Terry~\cite{MubayiTerry16a,MubayiTerry16b} study the number and structure of multigraphs in which no $s$ vertices support more than~$q$ edges for a large class of pairs~$(s, q)$.)

	As always, we begin by first proving an extremal result, with the counting result following immediately from our extremal result by Corollary~\ref{corollary: speed of arbitrary hereditary properties}. We note that similar extremal problems on multigraphs have been considered before by F\"uredi and K\"undgen~\cite{FurediKundgen02}. However the crucial difference is that as far as counting results are concerned, we need to determine the asymptotically extremal entropy, rather than the asymptotically extremal total number of edges as was the goal in~\cite{FurediKundgen02}. Indeed, in our problem, there exist configurations which are extremal with respect to the number of edges but \emph{not} with respect to entropy. 
	\begin{example}
	Consider a balanced bipartition $V_1\sqcup V_2$ of $[n]$ and let $G_1$ be the multigraph assigning weight~$2$ to every edge from $V_1$ to $V_2$ and weight~$0$ to every other edge. Let also $t_1$ be the associated template, assigning colour list~$\{0,1,2\}$ to every edge from $V_1$ to $V_2$ and colour list~$\{0\}$ to every other edge.

Clearly $G_1\in \langle t_1 \rangle\subseteq \mathcal{P}$.  The total edge weight of~$G_1$ is $\lfloor \frac{n^2}{2}\rfloor$, and the entropy of~$t_1$ is $\log_5(3) \left\lfloor \frac{n^2}{4}\right\rfloor$. 
	\end{example}

It is not hard to show that the total edge weight of~$G_1$ is extremal:
	\begin{theorem}\label{theorem: extremal size (3,5) problem}
	If $G$ is a multigraph in\/ $\mathcal{P}_n$, for some $n\geq 3$, then $e(G)\leq \lfloor \frac{n^2}{2}\rfloor$.
	\end{theorem}
	\begin{proof}
	By induction on $n$. The base cases $n=3$,~$4$ are easily checked by hand. For $n\geq 4$, consider a multigraph~$G\in\mathcal{P}_{n+1}$ on $n+1$ vertices with $e(G)\geq \lfloor \frac{(n+1)^2}{2}\rfloor$. We claim we must have in fact equality. By the inductive hypothesis it is enough to show that we can find a pair of vertices adjacent to a total of at most $\lfloor \frac{(n+1)^2}{2}\rfloor-\lfloor \frac{(n-1)^2}{2}\rfloor=2n$ edges. Suppose $G$ contains an edge~$u_1u_2$ with weight~$3$. Then for every other vertex~$v$, the pairs $u_1v$,~$u_2v$ can have combined weight at most~$1$, whence $u_1u_2$ is adjacent to at most $4+(n-2)<2n$ edges. We may thus assume all edges in $G$ have weight at most~$2$, and, since $\lfloor \frac{(n+1)^2}{2}\rfloor >\binom{n+1}{2}$ there must be some edge $u_1u_2$ with weight~$2$. Then for every other vertex $v$, the pairs $u_1v$, $u_2v$ can have combined weight at most~$2$, whence $u_1u_2$ is adjacent to at most $2 +2(n-1)=2n$ edges. Thus $e(G)=\lfloor \frac{(n+1)^2}{2}\rfloor$, as claimed.
	\end{proof}
The total edge weight of $G_1$ is thus maximal; however, the entropy of the associated template~$t_1$ is not. Indeed we can construct a different edge-extremal construction with strictly larger entropy.
	 	\begin{example} Let $M$ be a maximal matching in $[n]$ and let $G_2$ be the multigraph assigning weight~$2$ to every edge in $M$ and weight~$1$ to every other edge. Let also $t_2$ be the associated template, assigning colour list~$\{0,1,2\}$ to every edge of~$M$ and colour list~$\{0,1\}$ to every other edge.

	 	As before, we have $G_2\in \langle t_2 \rangle\subseteq \mathcal{P}$ and $e(G_2)=\lfloor \frac{n^2}{2}\rfloor$. However
	 	\[\Ent(t_2)=\log_5(2) \binom{n}{2}+\log_5\left(\frac{3}{2}\right)\left\lfloor \frac{n}{2}\right\rfloor=\log_5(\sqrt{2}) n^2+o(n^2) > \log_5(3^{1/4})n^2\geq \Ent(t_1).\]
	 	\end{example}
	 	It is straightforward to show that $t_2$ is indeed an entropy-extremal template for $\mathcal{P}$:
	\begin{theorem}\label{theorem: extremal entropy (3,5) problem}
For all $n\geq 3$, $\ex(n, \mathcal{P})=\log_5(2) \binom{n}{2}+\log_5\left(\frac{3}{2}\right)\left\lfloor \frac{n}{2}\right\rfloor$.
	\end{theorem} 
			\begin{proof}
			This is a proof by induction on $n$ again, very similar to the proof of Theorem~\ref{theorem: extremal size (3,5) problem}. The base cases $n=3$,~$4$ are again easily checked by hand. For $n\geq 4$, consider a template~$t$ for a $5$-colouring of~$E(K_{n+1})$ with $\langle t \rangle\subseteq \mathcal{P}_{n+1}$, with colours from $\{0,1,2,3,4\}$ corresponding to edge weights. Suppose $\Ent(t)\geq \log_5(2) \binom{n+1}{2}+\log_5\left(\frac{3}{2}\right)\left\lfloor \frac{n+1}{2}\right\rfloor$. We claim that we must in fact have equality. By the inductive hypothesis it is enough to show that we can find a pair of vertices~$u_1u_2$ such that the sum of the entropies of the edges incident to $u_1$ or~$u_2$ is at most~$2(n-1)\log_5(2)+ \log_5(3)$. By monotonicity of the property~$\mathcal{P}$, we may assume that for every edge~$e$ if $i<j$ and $j\in t(e)$ then $i\in t(e)$. Thus the possible entropies for a single edge are $0$ (weight~zero), $\log_5 (2)$ (weight $0$ or~$1$),  $\log_5(3)$ (weight $0$, $1$ or~$2$), and so on.

			Suppose $G$ contains an edge~$u_1u_2$ with entropy at least~$\log_5(4)$. Then $3\in t(u_1u_2)$, and thus for every other vertex~$v$, the combined weight of $u_1v$,~$u_2v$ in any realization of~$t$ must be at most~one, so that $\log_5 \vert t(u_1v)\vert +\log_5 \vert t(u_2v)\vert \leq \log_5(2)$. Thus the total entropy of the edges incident to $u_1$ or~$u_2$ is at most~$\log_5(5)+(n-1)\log_5(2)<2(n-1)\log_5(2)$. We may thus assume that every edge~$u_1u_2$ has entropy at most~$\log_5(3)$ in $t$, and, given the bound we are trying to prove, that there is some edge with entropy exactly~$\log_5(3)$. Then $2\in t(u_1u_2)$, and for every other vertex~$v$ the pairs $u_1v$,~$u_2v$ can have combined weight at most~$2$ in every realisation of $t$. In particular,
\[\log_5 \vert t(u_1v)\vert +\log_5 \vert t(u_2v)\vert \leq \max\bigl\{ \log_5(3)+\log_5(1), \log_5(2)+\log_5(2)\bigr\}=2\log_5(2).\]
Thus the total entropy of the edges incident to $u_1$ or~$u_2$ is at most~$2(n-1)\log_5(2)+\log_5(3)$, as required, and
\[\Ent(t) \leq \log_5(2) \binom{n+1}{2}+\log_5\left(\frac{3}{2}\right)\left\lfloor \frac{n+1}{2}\right\rfloor. \qedhere\]
			\end{proof}
We may thereby deduce a counting result for $\mathcal{P}$: 	
	\begin{corollary}
		There are $2^{\binom{n}{2}+o(n^2)}$ multigraphs on $[n]$ for which no triple of vertices supports more than~$4$ edges (counting multiplicities).
	\end{corollary}
	\begin{proof}
		Immediate from Theorem~\ref{theorem: extremal entropy (3,5) problem} and Corollary~\ref{corollary: speed of arbitrary hereditary properties}.
	\end{proof}
\begin{remark}
With only a little more work, it can be shown that $t_2$ and its isomorphic copies constitute a strong stability template for $\mathcal{P}$, and that typical members of~$\mathcal{P}$ are close to realisations of~$t_2$ --- and thus far from realisations of~$t_1$, despite the fact that $t_1$ was constructed from an edge-extremal graph.
\end{remark}

	\subsection{\texorpdfstring{$3$}{3}-coloured graphs}
	Let $\mathcal{P}$ denote the set of $3$-coloured graphs with no rainbow triangle, where a triangle is called \emph{rainbow} if it has an edge in each of the three colours~$\{1,2,3\}$. We use our multicolour container results to count the number of graphs in $\mathcal{P}$ and to characterise typical elements of~$\mathcal{P}$. This is related to the multicolour Erd{\H o}s--Rothschild problem~\cite{Erdos74}, which has received significant attention, see e.g.\ Alon, Balogh, Keevash and Sudakov's proof of a conjecture of Erd{\H o}s and Rothschild in~\cite{AlonBaloghKeevashSudakov04}, as well as the recent work of Benevides, Hoppen and Sampaio~\cite{BenevidesHoppenSampaio16}, Pikhurko, Staden and Yilma~\cite{PikhurkoStadenYilma16} and Hoppen, Lefmann and Odermann~\cite{HoppenLefmannOdermann}.

	\begin{theorem}[Extremal entropy]~\label{theorem: extremal entropy no rainbow k3}
Let\/ $\mathcal{P}$ denote the set of $3$-coloured graphs with no rainbow triangle.  For all $n\geq 3$, 
	\[\ex(n, \mathcal{P})=(\log_3 2)\binom{n}{2}.\]
	Furthermore, the unique extremal templates~$t$ are obtained by choosing a pair of colours~$\{c_1,c_2\}$ from $\{1,2,3\}$ and setting $t(e)=\{c_1,c_2\}$ for every $e\in E(K_n)$.
	\end{theorem}
	\begin{proof}
	Our theorem shall follow from the following observation and a straightforward averaging argument.
	\begin{observation}\label{observation: rainbow K3}
	Suppose $\langle t \rangle\subseteq \mathcal{P}$ and $e=\{v_1,v_2\}$ is some edge of~$K_n$. Then rainbow $K_3$-freeness implies the following:
	\begin{enumerate}[(i)]
	\item if $\vert t(e)\vert =3$, then for all $x\in V(K_n)\setminus e$ and $i \in \{1, 2\}$, we have $\vert t(xv_i)\vert =1$;
	\item if $\vert t(e)\vert =\vert t(f)\vert= 2$ and $t(e)\neq t(f)$, then $e\cap f = \emptyset$;
	\item if $\vert t(e)\vert =2$ and $c$ is the colour missing from $t(e)$, then for every $x\in V(K_N)\setminus e$, either $t(xv_1)=t(xv_2)=\{c\}$ or $c$ is missing from both $t(xv_1)$ and $t(xv_2)$. \qedhere
	\end{enumerate}
	\end{observation}
	In particular, for any $3$-set $A\subseteq [n]$, we have $\Ent(t\vert_A)\leq 3\log_32$, with equality attained if and only if all~three edges of $A$ are assigned the same pair of colours~$\{c_1,c_2\}$ by $t$.

	Now, suppose $t$ is a template with $\Ent(t)\geq \log_3 2 \binom{n}{2}$. The average entropy of $t\vert_A$ over all $3$-sets~$A\subseteq [n]$ is:
	\[\frac{1}{\binom{n}{3}} \sum_A \Ent\bigl(t_{\vert_A}\bigr) =\frac{1}{\binom{n}{3}} (n-2)\Ent(t)\geq 3\log_3 2.\]
	Our previous bound on the entropy inside triangles then tells us that we must have equality, and that $t$ must have entropy $3\log_32$ inside \emph{every} $3$-set $A$. In particular all edges $e$ must have $\vert t(e)\vert=2$. Finally by (ii) in Observation~\ref{observation: rainbow K3}, we must have $t(e)=\{c_1,c_2\}$ for some pair of colours $\{c_1,c_2\}$ and \emph{all} edges $e\in E(K_n)$. This concludes the proof of the theorem.
	\end{proof}
	\begin{corollary}[Counting]
	For all $\varepsilon>0$, there exists $n_0\in \N$ such that for all $n\geq n_0$,
	\[3\cdot3^{(\log_3 2)\binom{n}{2}} -3 \leq \vert \mathcal{P}_n\vert \leq 3^{(\log_3 2)\binom{n}{2}+\varepsilon \binom{n}{2}}. \]
	\end{corollary}
	\begin{proof}
	The lower bound equals the number of colourings of~$E(K_n)$ such that each edge receives one of a prescribed pair of colours.  For the upper bound, Theorem~\ref{theorem: extremal entropy no rainbow k3} gives $\pi(\mathcal{P})=\log_32$, and the claimed result then follows from Corollary~\ref{corollary: speed of arbitrary hereditary properties}.
	\end{proof}

We note that the stronger bound $\abs{\PP_n} \leq 3^{(\log_3 2)\binom{n}{2}+ O(n \log n)}$ was proved in~\cite{BenevidesHoppenSampaio16}.

	With a bit more case analysis, we can obtain the following stability result (see the Appendix for a proof).  Recall the definition of a strong stability family (Definition~\ref{definition: strong stability family}).
	\begin{theorem}[Stability]\label{theorem: stability no rainbow K3}
	The sequence of templates~$t$ such that there exist a pair of colours~$\{c_1, c_2\}$ such that $t(e) = \{c_1, c_2\}$ for all $e \in E(K_n)$ is a strong stability family for\/ $\mathcal{P}$.  That is, for all $\varepsilon>0$, there exist $\delta=\delta(\varepsilon)>0$ and $n_0=n_0(\delta)\in \N$ such that the following holds: if $t$ is a $3$-colouring template on $n\geq n_0$ vertices satisfying 
	\begin{enumerate}[(i)]
		\item $\Ent(t)\geq (\log_3 2  -\delta) \binom{n}{2}$, and
		\item there at most $\delta \binom{n}{3}$ rainbow triangles in $K_n$ which can be realised from $t$,
	\end{enumerate}
	then there exists a pair of colours $\{c_1,c_2\}$ such that $t(e)=\{c_1,c_2\}$ for all but at most $\varepsilon \binom{n}{2}$ edges of~$K_n$. 
	\end{theorem}

	It follows from Theorem~\ref{theorem: stability no rainbow K3} that the family of template sequences $\mathcal{S}=\{\mathbf{t}_A : \ A\in [3]^{(2)}\}$, where $(t_A)_n$ is the colouring template for $K_n$ that assigns the colour pair~$\{c_1, c_2\}$ to every edge, is a strong stability family for $\mathcal{P}$.
	\begin{corollary}[Typical colourings]
	Almost all $3$-coloured graphs with no rainbow triangle are almost $2$-coloured: for every $\varepsilon>0$ there exists $n_0$ such that for all $n\geq n_0$ at most $3^{\varepsilon \binom{n}{2}}$ rainbow $K_3$-free $3$-colourings of~$K_n$ have at least $\varepsilon \binom{n}{2}$ edges in each of the colours $\{1,2,3\}$.
	\end{corollary}
	\begin{proof}
	Instant from Theorems \ref{theorem: strong stability and containers} and~\ref{theorem: stability no rainbow K3}.
	\end{proof}

	We note that there are (many) examples of rainbow $K_3$-free $3$-coloured graphs in which all three colours are used.  Indeed, consider a balanced bipartition $[n]=A\sqcup B$. Colour the edges from $A$ to $B$ Red, and then arbitrarily colour the edges internal to $A$ Red or Blue and the edges internal to $B$ Red or Green.  The resulting $3$-colouring has no rainbow $K_3$, and by randomly colouring the edges inside $A$ and~$B$ we can in fact ensure that all three colours are used on at least $(1+o(1))\frac{n^2}{16}$ edges.

\subsection{Hypercubes}
Let $\mathcal{P}$ be the vertex-hypercube property of not containing a copy of the square~$Q_2$---that is of not containing four distinct subsets of~$[n]$ of the form $A$, $A\sqcup\{i\}$, $A\sqcup\{j\}$ and~$A\sqcup \{i,j\}$. 
We have $\pi^{\Hypercube}_v(\mathcal{P})\geq \frac{2}{3}$, as may be seen for example by removing every third layer of~$Q_n$, i.e.\ taking as our construction the family of all $\vector{x}\in Q_n$ with $\sum_i \vector{x}_i\not\cong 0 \mod 3$, which clearly contains no induced square. Kostochka~\cite{Kostochka76} and, later and independently, Johnson and Entringer~\cite{JohnsonEntringer89} showed that this lower bound is tight: 
\[\pi^{\Hypercube}(\mathcal{P})=\frac{2}{3}.\]
By Corollary~\ref{corollary: vertex hypercube counting} this immediately implies the following counting result:
\begin{corollary}
There are $\vert \mathcal{P}_n\vert = 2^{\left(\frac{2}{3}+o(1)\right)2^n}$ induced $Q_2$-free subgraphs of~$Q_n$.
\end{corollary}

In a different direction, let $\mathcal{Q}$ be the property of hypercube subgraphs -- i.e.\ of $2$-colourings of~$\Hypercube$ -- of not containing a copy of the square~$Q_2$. A long-standing conjecture of Erd{\H o}s~\cite{Erdos84} states that the edge-Tur\'an density (entropy density relative to $\Hypercube$) of this property is~$\pi^{\Hypercube}(\mathcal{Q})=1/2$. The lower bound is obtained by deleting all edges between layer~$2i$ and layer~$2i+1$ for $0\leq i \leq \lfloor n/2\rfloor$. The best upper bound to date is~$0.603\ldots$ from applications of flag algebras due to Baber~\cite{Baber12} and Balogh, Hu, Lidick{\'y} and Liu~\cite{BaloghHuLidickyLiu14}. Again, by Corollary~\ref{corollary: edge hypercube counting} we have the following:
\begin{corollary}
If Erd{\H o}s's conjecture on $\pi^{\Hypercube}(\mathcal{Q})$ is true, then there are
\[\vert \mathcal{Q}_n\vert =2^{\left(\frac{1}{2}+o(1)\right)2^{n-1}n}=2^{\left(\frac{n}{4}+o(n)\right)2^n}\]
$Q_2$-free subgraphs of $Q_n$.
\end{corollary}

\subsection{A non-example: sparse graph sequences} \label{subsection: nonexample containers fail for paths}
Let $\Path=(P_n)_{n\in\N}$ be the sequence of paths on $[n]$ introduced in Section~\ref{subsection: other host graphs}. An easy calculation reveals that $\Path$ fails to satisfy the `goodness' condition introduced in Definition~\ref{definition: good graph sequence}, and is therefore not covered by Theorem~\ref{theorem: general multicolour container result for graph sequences}. As mentioned earlier, there is a good reason for this: the conclusion Theorem~\ref{theorem: general multicolour container result for graph sequences} does not hold for $\Path$ (or, more generally, for tree-like graph sequences).

Let $\mathcal{P}$ be the order-hereditary property of $3$-colourings of $\Path$ of not having two consecutive edges in the same colour. It is easy to see that $\vert \mathcal{P}_n\vert= 3\cdot 2^{n-2}=3^{n\log_3(2) -O(1)}$. On the other hand, the extremal entropy of $\mathcal{P}_n$ is only about $ n \log_3\sqrt{2}$.
\begin{theorem}
For any $n\geq 3$, $\ex(P_n, \mathcal{P})=\lceil (n-1)/2\rceil \log_32$.
\end{theorem}
\begin{proof}
If $f$ and $f'$ are consecutive edges and $t$ is a $3$-colouring template with $\langle t \rangle\subseteq \mathcal{P}$ then $t(f)\cap t(f')=\emptyset$, from which it follows that $\log_3(\vert t(f)\vert) + \log_3(\vert t(f')\vert)\leq \log_32$. Further there can be no edge~$f$ with $t(f)=[3]$, since otherwise we would have a realisation of $t$ with two consecutive edges of the same colour. Partitioning the path~$P_n$ into disjoint pairs of consecutive edges and at most~one single edge, we get $\Ent(t)\leq \lceil(n-1)/2\rceil \log_3(2)$ as desired. For the lower bound, consider the template~$t$ defined by setting $t(\{2i+1, 2i+2\})=[2]$ and $t(\{2i, 2i+1\})=\{3\}$ for $0\leq i\leq \lfloor (n-1)/2\rfloor$. This has the correct entropy and all of its realisations clearly lie in $\mathcal{P}$. 
\end{proof} 
Now, $\binom{P_n}{P_3}=n-3$, and it is easy to see that we have supersaturation of sorts for $\mathcal{P}$: if $t$ is a template with $\Ent(t)\geq n \log_3(\sqrt{2})+\varepsilon n$, there are at least~$\Omega(\varepsilon n)=\Omega(\varepsilon\binom{P_n}{P_3})$ pairs of consecutive edges which can be made monochromatic in some realisation of~$t$. In particular, templates having $o(n)$ such pairs must have entropy at most~$\log_3(\sqrt{2}) n +o(n)$. A collection of $3^{o(n)}$ such templates can thus cover at most~$2^{n/2}3^{o(n)}=o(2^n)= o(\vert \mathcal{P}_n\vert)$ colourings---in particular, it cannot form a container family for $\mathcal{P}_n$. This shows that the analogue of Theorem~\ref{theorem: general multicolour container result for graph sequences} does not hold for the graph sequence $\Path$, and that the `goodness' condition in the statement of that theorem is necessary, as we claimed.



\section{A cut metric for \texorpdfstring{$k$}{k}-decorated graphons}\label{section: cut metric for graphons}

In this and the following sections, we turn our attention to limits of sequences of $k$-coloured graphs.  As mentioned earlier, a thorough treatment of the theory of graph limits is given in the monograph~\cite{LovaszBook}.

\subsection{Notation and definitions}\label{se:DecoratedDefs}

Recall that, given a set~$K$, a \emph{$K$-decorated graph} with vertex set~$[n]$ is a labelling of~$E(K_n)$ with elements of~$K$.  In particular, $[k]$-decorated graphs correspond to $k$-coloured graphs.
In~\cite{LS:decorated}, Lov\'asz and Szegedy extended ideas and results from graph limit theory, including homomorphism densities and convergence, to $K$-decorated graphs, where $K$ is any second-countable compact Hausdorff space.  Before we state and prove our results, we shall give an overview of important definitions and notation from~\cite{LS:decorated}  (see also~\cite[Chapter~17]{LovaszBook}).  As suggested above, we primarily consider the case $K = [k]$.  Fortunately, in this case, most of the necessary definitions are straightforward variations on the corresponding definitions for ordinary graph limits.  That said, some of our results extend to $K$-decorated graphs for arbitrary~$K$.  In order to describe these results, we shall also give some definitions in full generality.

First, we discuss homomorphism densities.  As observed in~\cite{LS:decorated}, in general, there is no natural way to define the homomorphism density of one $K$-decorated graph into another.  Instead, it makes sense to define homomorphisms from $C[K]$-decorated graphs into $K$-decorated graphs, as follows.
Suppose that $G$ is a $C[K]$-decorated graph and that $H$ is a $K$-decorated graph such that $\abs{V(H)} \geq \abs{V(G)}$.  For every map $\phi : V(G) \to V(H)$, define
\begin{equation}\label{eq:homphi}
\hom_{\phi}(G, H) = \prod_{1 \leq i < j \leq \abs{V(H)}} G_{ij}\bigl(H_{\phi(i)\phi(j)}\bigr).
\end{equation}
Let
\[
\hom(G, H) = \sum_{\phi : V(G) \to V(H)} \hom_{\phi}(G, H).
\]
The \emph{homomorphism density} of $G$ into $H$ is
\begin{equation}\label{eq:homdensity}
t(G, H) := \dfrac{\hom(G, H)}{\abs{V(H)}^{\abs{V(G)}}}.
\end{equation}

Observe that $C([k]) \cong \R^k$.  The definition~\eqref{eq:homphi} only makes sense when the edges of~$G$ are labelled with functions, but in the case $K = [k]$, it will usually be sufficient to think of the labels as vectors.

Now we discuss convergence of sequences of $K$-decorated graphs.  Much as for ordinary graphs, we say that a sequence of $K$-decorated graphs~$\sequence{G}$ is is \emph{convergent} if and only if for every $C[K]$-decorated graph~$F$, the sequence of homomorphism densities~$\{t(F, G_n)\}$ converges.

Finally, we consider representations of the limit of a convergent sequence.  If $\sequence{G}$ is a convergent sequence of simple graphs, then its limit~$\Gamma$ can be represented by a \emph{graphon}, a symmetric, measurable function~$W : \oi^2 \to \oi$.  An alternative way of thinking about a graphon~$W$ is as a function that maps each $(x, y) \in \oi^2$ to the probability distribution on $\{0, 1\}$ that assigns mass~$W(x, y)$ to~$1$ and mass~$1 - W(x, y)$ to~$0$.
This perspective yields the appropriate generalization to the case of $K$-decorated graphs.  

\begin{definition}[Decorated graphons and kernels]\label{definition: decorated graphons}
Given~$K$, let $\BorelMeasures(K)$ denote the set of finite Borel measures on $K$ and let $\measures(K) \subseteq \BorelMeasures(K)$ denote the set of Borel probability measures on $K$.
A \emph{$K$-decorated graphon}, or \emph{$K$-graphon}, is a symmetric, measurable function~$W : \oi^2 \to \measures(K)$.
A \emph{$K$-kernel} is a symmetric, measurable function~$W : \oi^2 \to \BorelMeasures(K)$.
\end{definition}
It is shown in~\cite{LS:decorated} that the limit of a sequence of $K$-decorated graphs can be represented by a $K$-graphon.

Similarly to ordinary graphons, this representation is defined in terms of homomorphism densities of $C[K]$-decorated graphs into $K$-graphons, which we now define.  Given a $K$-graphon~$W$, $f \in C[K]$ and $(x, y) \in \oi^2$, we define
\begin{equation}\label{eq:WfDef}
W_f(x, y) = \int_K f \, dW(x,y).
\end{equation}
Let $F$ be a $C[K]$-decorated graph and, for all $i$,~$j \in V(F)$, let $F_{i j}$ denote the label of~$ij$.  If $W$ is a $K$-graphon, then
\begin{equation}\label{eq:KGraphonHomomorphismDensity}
t(F, W) = \int_{\oi^{v(F)}} \prod_{i < j} W_{F_{i j}}(x_i, x_j) \, dx_1 \cdots dx_{v(F)}.
\end{equation}
It is shown in~\cite[Theorem~2.5]{LS:decorated} that if $\sequence{G}$ is a convergent sequence of $K$-decorated graphs, then there exists a $K$-graphon~$W$ such that for all $C[K]$-decorated graphs~$F$, $t(F, G_n) \to t(F, W)$.

Different $K$-graphons may represent the same graph limit.  We say that two $K$-graphons $W$ and~$W'$ are \emph{equivalent}, and write $W \cong W'$, if $t(F, W) = t(F, W')$ for every $C[K]$-decorated graph~$F$.

For each $j \in [n]$, let $I_j = [(j - 1)/n, j/n)$.  Given a $K$-decorated graph~$G$, we define a $K$-graphon~$W_G$ by setting $W_G(x, y) = \delta_{G_{ij}}$ for all $(x, y) \in I_i \times I_j$, where $\delta$ denotes the Dirac delta measure.  It is easy to see that if $F$ is a $C[K]$-decorated graph, then
\[
t(F, W_G) = t(F, G).
\]

We write $k$-graphon for $[k]$-graphon and $k$-decorated (or $k$-coloured) graph for $[k]$-decorated graph.  The definitions of homomorphism densities into $k$-decorated graphs, the functions~$W_f$, and homomorphism densities into $k$-graphons are identical to those in \eqref{eq:homdensity}--\eqref{eq:KGraphonHomomorphismDensity}.

Frieze and Kannan~\cite{FK99} introduced a ``cut norm'' $\cutnorm{\cdot}$ that has become central to the theory of graph limits.  (For an overview of the history of the cut norm in other contexts, see~\cite[Section~4]{Jan13}.)  Given a graphon~$W$, the \emph{cut norm} of~$W$ is
\begin{equation*}\label{eq:CutNormDef}
\cutnorm{W} = \sup_{S, T \subseteq \oi} \biggl\lvert \int_{S \times T} W(x, y)\,dx dy\biggr\rvert,
\end{equation*}
where the supremum is over all pairs of measurable subsets of~$\oi$.  If $U$ and $W$ are graphons, then
\begin{equation}\label{eq:CutDistanceDef}
\dcut(U, W) = \cutnorm{U - W} = \sup_{S, T \subseteq \oi} \biggl\lvert \int_{S \times T} \bigl(U(x, y) - W(x, y)\bigr)\,dx dy\biggr\rvert.
\end{equation}

Lov\'asz and Szegedy~\cite{LS:decorated} did not consider a version of the cut distance for $K$-decorated graphs.  We introduce an appropriate generalization here.  If $G$ and $H$ are two $k$-decorated graphs with vertex set~$[n]$, we define
\begin{equation}\label{eq:dkcutGraphs}
\dkcut(G, H) = \max_{S, T \subseteq [n]} \dfrac{1}{n^2} \sum_{i=1}^k\biggl\lvert \sum_{(u, v) \in S \times T} \bigl(\mathds{1}(G_{uv} = i) - \mathds{1}(H_{uv} = i)\bigr)\biggr\rvert.
\end{equation}

Given a $k$-graphon~$W$, for all $i \in [k]$, set
\begin{equation}\label{eq:omegaiDef}
W_i(x, y) = \Prob\bigl(W(x, y) = i\bigr)
\end{equation}
and observe that each $W_i$ is a graphon.  If $W$ and $U$ are $k$-graphons, we define
\begin{equation}\label{eq:dkcutDef}
\dkcut(U, W) = \sup_{S, T \subseteq \oi} \sum_{i=1}^k \biggl\lvert \int_{S \times T} \bigl(U_i(x, y) - W_i(x, y)\bigr)\,dx dy\biggr\rvert.
\end{equation}
Given a measure-preserving transformation~$\varphi : \oi \to \oi$, we define $W^{\varphi}$ by $W^{\varphi}(x, y) = W(\varphi(x), \varphi(y))$.  The \emph{cut distance} between $U$ and $W$ is
\begin{equation*}\label{eq:deltakcutDef}
\deltakcut(U, W) = \inf_{\varphi : \oi \to \oi} \dkcut(U, W^{\varphi}),
\end{equation*}
where the infimum is taken over all measure-preserving transformations~$\varphi : \oi \to \oi$.

We sometimes write, e.g., $\dkcut(G, W)$ for $\dkcut(W_G, W)$.

\begin{remark}\label{re:DifferentCutNorms}
If $U$ and $W$ are graphons, then comparing \eqref{eq:dkcutDef} (with $k = 2$) and~\eqref{eq:CutDistanceDef} shows that $d_{\square_2}(U, W) = 2\cutnorm{U - W}$.
\end{remark}

\begin{remark}\label{re:CutNormDisjoint}
It is easy to see that if we restrict the maximum on the right-hand side of~\eqref{eq:dkcutGraphs} to pairs~$(S, T)$ with $S \cap T = \emptyset$, then the resulting quantity is at least~$\dkcut(G, H) / 4$.  The same observation holds for~\eqref{eq:dkcutDef}.
\end{remark}

Given two $k$-graphons $U$ and~$W$, let
\begin{equation}\label{eq:L1kgraphons}
d_1(U, W) = \sum_{i=1}^k  \int_{\oi^2} \bigl\lvert U_i(x, y) - W_i(x, y)\bigr\rvert\,dx dy
\end{equation}
and let
\[
\delta_1(U, W) = \inf_{\varphi : \oi \to \oi}  d_1(U, W^{\varphi}),
\]
where once again the infimum is taken over all measure-preserving transformations~$\varphi : \oi \to \oi$.
Observe that 
\begin{equation*}\label{eq:CutNormL1NormIneq}
\deltakcut(U, W) \leq \delta_1(U, W).
\end{equation*}

\begin{definition}[Step-function]\label{definition: step-function}
We say that a $k$-kernel~$W$ is a \emph{step-function} if there exists a partition~$\partition$ of~$\oi$ such that $W$ is constant on each product of cells of~$\partition$.
\end{definition}

We shall consider two ways of randomly generating decorated graphs from a $k$-kernel.  Recall that $\BorelMeasures([k])$ denotes the set of measures on $[k]$ and that $\measures([k])$ denotes the set of  probability measures on $[k]$. Given a $[k]$-kernel~$W$, we define a random $\BorelMeasures([k])$-decorated graph~$H(n, W)$ as follows.  Let $X_1$, \dots,~$X_n$ be i.i.d.\ uniformly chosen points from $\oi$.  For $1 \leq i < j \leq n$, label $H(n, W)_{ij}$ with $W(X_i, X_j)$.

Observe that $H(n, W)$ naturally defines a $k$-kernel.  Recall that, for $1 \leq i \leq n$, we set $I_i = [(i-1)/n, i/n)$.  For each $(x, y) \in I_i \times I_j$, let
\begin{equation}\label{eq:omegaHnDef}
W_{H_n}(x, y) = H(n, W)_{ij} = W(X_i, X_j).
\end{equation}
%
Note in particular that if $W$ is a $k$-graphon, then $H(n, W)$ is a $\measures([k])$-decorated graph and $W_{H_n}$ is a $k$-graphon.

If $W$ is a $k$-graphon and $G$ is a $k$-decorated graph, we shall want to compute the cut distance between $W_{H_n}$ and $G$.   Because $W_{H_n}$ and $W_G$ are each constant on products of intervals of the form~$[(j - 1)/n, j/n)$, we may restrict the supremum on the right-hand side of~\eqref{eq:dkcutDef} to pairs~$(S, T)$ where each of $S$~and~$T$ is a union of intervals of the form~$[(j - 1)/n, j/n)$.  Hence,
\begin{equation}\label{eq:HnwCutNorm}
\dkcut\bigl(G, W_{H_n}\bigr) = \max_{S, T \subseteq [n]} \dfrac{1}{n^2} \sum_{i=1}^k\biggl\lvert \sum_{(u, v) \in S \times T} \Bigl(\mathds{1}(G_{uv} = i) - \Prob\bigl(W(X_u, X_v) = i\bigr)\Bigr)\biggr\rvert.
\end{equation}
Also, by~\eqref{eq:L1kgraphons}, if $U$ and $W$ are two $k$-kernels, then
\begin{equation}\label{eq:L1DistanceH}
d_1(U_{H_n}, W_{H_n}) =\dfrac{1}{n^2} \sum_{i=1}^k \sum_{(u, v) \in [n]^2} \Bigl\lvert \bigl(U_i(X_u, X_v) - W_i(X_u, X_v)\bigr) \Bigr \rvert.
\end{equation}

The second way of generating random graphs only makes sense when $W$ is a $k$-graphon.  Because $H(n, W)$ is $\measures([k])$-decorated, we may think of it as a model of random $k$-decorated graphs.  Let $G(n, W)$ denote the random $k$-decorated graph such that the label of each edge~$ij$ is obtained by sampling from $H(n, W)_{ij}$ independently of all other edges.

Homomorphism densities provide a convenient description of the distribution of~$G(n, W)$.  Given a $k$-decorated graph~$F$ on $n$ vertices, we define an associated $\R^k$-decorated graph~$G_F$ on $n$ vertices as follows: for each pair $(i, j)$, we label $(G_F)_{ij}$ with the basis vector~$e_{F_{ij}}$.  If $W$ is a $k$-graphon, then \eqref{eq:WfDef} and~\eqref{eq:omegaiDef} imply that
\[
W_{e_{F_{ij}}}(x, y) = \Prob\bigl(W(x, y) = F_{ij}\bigr) = W_{F_{ij}}(x, y).
\]
Hence, \eqref{eq:KGraphonHomomorphismDensity} becomes
\begin{equation}\label{eq:GFdensityFinite}
t(G_F, W) =  \int_{\oi^n} \prod_{1 \leq i < j \leq n} \Prob\bigl(W(x_i, x_j) = F_{ij}\bigr)\,dx_1 \cdots dx_n.
\end{equation}
Now observe that
\[
\Prob\bigl(G(n, W) = F \mid X_1, \dots, X_n\bigr) = \prod_{1 \leq i < j \leq n} \Prob\bigl(W(X_i, X_j) = F_{ij}\bigr).
\]
Taking expectations and comparing the result with~\eqref{eq:GFdensityFinite} implies that
\begin{equation}\label{eq:wRandomDistribution}
\Prob\bigl(G(n, W) = F\bigr) = t(G_F, W).
\end{equation}

\subsection{Main results}\label{se:compactness}

Let $\WWk$ denote the space of $k$-graphons and let $\widetilde{\WWk}$ denote the quotient of~$\WWk$ obtained by identifying $U$ and $W$ whenever $\deltakcut(U, W) = 0$.

The main result of this section is the following.  

\begin{theorem}\label{th:compactnessFinite}
The space $(\widetilde{\WWk}, \deltakcut)$ is compact.
\end{theorem}

Most of the work in the proof of Theorem~\ref{th:compactnessFinite} will go toward proving that equivalent $k$-graphons have cut distance~$0$.

\begin{theorem}\label{th:deltakcut}
Let $U$ and $W$ be $k$-graphons.  If $U \cong W$, then $\deltakcut(U, W) = 0$.
\end{theorem}

\begin{remark}\label{re:WantToProveEquivalenceFinite}
The graphon case of Theorem~\ref{th:deltakcut} was originally proved in~\cite{BCLSV08} using complicated estimates of the relationship between rates of convergence in subgraph counts and the rate of convergence in cut distance.  However, as observed by Schrijver (see~\cite[Remark 11.4]{LovaszBook}), it is possible to prove the result using simpler analytic arguments, and this is the course that we shall pursue.
\end{remark}

In order to prove Theorem~\ref{th:compactnessFinite}, we shall also need two results from~\cite{LS:decorated}.  The first result says that $\WWk$ is compact with respect to the topology defined by convergence of homomorphism densities.

Let $B^{\infty} \subseteq (\R^k, \norm{\cdot}_{\infty})$ denote the unit ball.

\begin{theorem}\label{th:ExistenceFinite}
If $\sequence{W}$ is a sequence of $k$-graphons such that the sequence~$\{t(F, W_n)\}$ converges for every $B^{\infty}$-decorated graph~$F$, then there exists a $k$-graphon~$W$ such that $t(F, W_n) \to t(F, W)$ for every $B^{\infty}$-decorated graph~$F$.
\end{theorem}

The second result is a counting lemma that is essentially identical to~\cite[Lemma~3.5]{LS:decorated}.  

\begin{lemma}[Counting Lemma]\label{le:CountingFinite}
Let $U$ and $W$ be $k$-graphons and, for all $i \in [k]$, let $m_i = \max(\norm{U_i}_{\infty}, \norm{W_i}_{\infty})$.  If $F$ is a $B^{\infty}$-decorated graph on $q$ vertices, then
\begin{equation*}\label{eq:countingFinite}
\bigl\lvert t(F, U) - t(F, W) \bigr\rvert \leq 4\dbinom{q}{2} \biggl(\prod_{u, v \in V(F)} m_{F_{u, v}}\biggr) \deltakcut(U, W).
\end{equation*}
\end{lemma}

\begin{remark}\label{rk:CompactnessGeneral}
Theorem~\ref{th:compactnessFinite} in fact holds for arbitrary $K$.  In this case, the definition of the cut norm differs slightly from~\eqref{eq:dkcutDef}: given two $K$-graphons $U$ and~$W$, we set
\begin{align*}\label{eq:dKcutDef}
\dKcut(U, W) &= \sup_{f \in B^{\infty}} \sup_{S, T \subseteq \oi} \biggl\lvert \int_{S \times T} \bigl(U_f(x, y) - W_f(x, y)\bigr)\,dx dy\biggr\rvert \nonumber\\
	&= \sup_{f \in B^{\infty}} \cutnorm{U_f - W_f}.
\end{align*}

We do not see a way to generalize the proof of Theorem~\ref{th:compactnessFinite} given below to arbitrary $K$.  Instead, one can adapt the proof of compactness for graphons given in~\cite{LS:analyst}, which we now sketch.

Let $\GenSet$ be a countable dense subset of~$B^{\infty}$; note that the Stone--Weierstrass theorem guarantees that such a set exists.  Let $\sequence{W}$ be a sequence of $K$-graphons that convergences in cut distance.  For each $n$ and each $f \in \GenSet$, let $W_{n,f}$ be as in~\eqref{eq:WfDef}.  Using the Weak Regularity Lemma (Lemma~\ref{le:WeakRegGraphons} below) and the Martingale Convergence Theorem, it is possible to show that there exists a subsequence~$\{n_k\}_{k=1}^{\infty}$ such that for all $f \in \GenSet$, $W_{n_k,f}$ converges in cut distance to a graphon~$U_f$.

By a result from~\cite{LS:decorated}, the collection~$\{U_f\}_{f \in \GenSet}$ specifies a $K$-graphon~$U$.  To conclude the proof, one shows that $\deltaKcut(W_n, U) \to 0$ using a $3\varepsilon$-argument.

Finally, we do not know whether Theorem~\ref{th:deltakcut} holds for arbitrary $K$, and leave this as a n open problem.  
\end{remark}

\subsection{Proofs}\label{se:CompactnessProofs}

Now we prove our main results.

\begin{proof}[Proof of Theorem~\ref{th:deltakcut}.]
Fix $\varepsilon > 0$.  We prove the result in three steps.  Our argument is similar to Schrijver's proof of the result for graphons, which was mentioned in Remark~\ref{re:WantToProveEquivalenceFinite}.  It is outlined in~\cite[Exercise~11.27]{LovaszBook}.

Recall the definition of~$W_{H_n}$ from~\eqref{eq:omegaHnDef}.

\begin{claim}\label{cl:wHnConvergence}
Almost surely, $\delta_1(W_{H_n}, W) \to 0$ as $n \to \infty$.
\end{claim}

\begin{proof}
Here, it is convenient to allow $W$ to be a $k$-kernel.  We first consider the case when $W$ is a step-function.
Observe that if $A \subseteq \oi$ is a measurable set, then
\[
\dfrac{1}{n}\bigl\lvert\set{i : X_i \in A}\bigr\rvert \pto \mu(A),
\]
where $\mu$ denotes the Lebesgue measure.  It follows that for each $n$, there is a rearrangement of~$H(n, W)$, which we denote by $\widetilde{H(n, W)}$, such that
\begin{equation}\label{eq:ptoWRandom}
\mu\bigl(\set{W_{\widetilde{H_n}} \neq W}\bigr) \pto 0
\end{equation}
as $n \to \infty$.  Now
\[
d_1(W_{\widetilde{H_n}}, W) = \sum_{i=1}^k \iint_{\oi^2} \lvert W_{\widetilde{H_n}, i} - W_i \rvert = \sum_{i=1}^k \iint_{\oi^2} \lvert W_{\widetilde{H_n}, i} - W_i \rvert \mathds{1}(W_{\widetilde{H_n}} \neq W) \leq \mu\bigl(\set{W_{\widetilde{H_n}} \neq W}\bigr).
\]
It follows from~\eqref{eq:ptoWRandom} that $d_1(W_{\widetilde{H_n}}, W) \pto 0$ and thus that $\delta_1(W_{H_n}, W) \pto 0$.  Finally, by passing to a subsequence, we have that $\delta_1(W_{H_n}, W) \to 0$ almost surely. 

Now let $W$ be an arbitrary $k$-kernel.  For each $i \in [k]$, there exists a non-decreasing sequence of measurable step-functions~$(W_{m,i})_{m=1}^{\infty}$ such that $W_{m,i} \to W_i$ pointwise as $m \to \infty$.  Let $W_m$ be the $k$-kernel defined by the $W_{m,i}$.  (Even if $W$ is a $k$-graphon, we do not require the $W_m$ to be $k$-graphons.)  Each $W_m$ is a step-function.  By the triangle inequality, for all~$n$,
\begin{equation}\label{eq:wHnConvergenceTriangleIneq}
\delta_1(W, W_{H_n}) \leq \lim_{m \to \infty} \delta_1(W, W_m) + \inf_{m \geq 1} \delta_1(W_m, W_{m_{H_n}}) + \lim_{m \to \infty} \delta_1(W_{m_{H_n}}, W_{H_n}).
\end{equation}
Thus, the claim will follow if we can show that almost surely each term on the right-hand side of~\eqref{eq:wHnConvergenceTriangleIneq} tends to~$0$ as $n \to \infty$.

First, by definition, $W_m \to W$ pointwise as $m \to \infty$.

Second, by compactness and the fact that each $W_m$ is a step-function, there exists a subsequence such that for all~$m$, almost surely $\delta_1(W_{m_{H_n}}, W_m) \to 0$ as $n \to \infty$.

Finally, we show that for all~$n$, $d_1(W_{m_{H_n}}, W_{H_n}) \to 0$ as $m \to \infty$.  For each~$n$, we may couple $H(n, W)$ and the $H(n, W_m)$ by generating them from the same sequence of random variables $X_1$, \dots,~$X_n$.  For each $(u, v) \in [n]^2$, let $\lambda_{u,v} = W(X_u, X_v)$ and let $\lambda^{(m)}_{u,v} = W_m(X_u, X_v)$.  Thus, for each $i \in [k]$, \eqref{eq:L1DistanceH} implies that
\[
\norm{W_{m_{H_n}, i} - W_{H_n, i}}_1 = \dfrac{1}{n^2} \biggl\lvert \sum_{(u, v) \in [n]^2} \bigl(\lambda^{(m)}_{u,v}(i) - \lambda_{u,v}(i) \bigr) \biggr\rvert.
\]
It follows that
\begin{equation*}\label{eq:wHnStepFunctionL1Distance}
\lim_{m \to \infty} \norm{W_{m_{H_n}, i} - W_{H_n, i}}_1 \leq \lim_{m \to \infty} \dfrac{1}{n^2} \sum_{(u, v) \in [n]^2} \bigl\lvert \lambda^{(m)}_{u,v}(i) - \lambda_{u,v}(i) \bigr\rvert = 0.
\end{equation*}
This completes the proof.
\end{proof}

\begin{claim}\label{cl:wHnGnwConvergence}
Almost surely, $\dkcut(W_{H_n}, G(n, W)) \to 0$ as $n \to \infty$.
\end{claim}

\begin{proof}
Given disjoint sets $S$,~$T \subseteq [n]$ and $i \in [k]$, let
\[
Z^{(S, T)}_i = \bigl\lvert\bigl\{ (u, v) \in S \times T \, : \, G(n, W)_{uv} = i \bigr\}\bigr\rvert.
\]
We shall use a concentration inequality to bound the deviation of each~$Z^{(S, T)}_i$ from its mean.  The argument is very similar to the proof of~\cite[Lemma~10.11]{LovaszBook}.

Because $G(n, W)$ is obtained by sampling from $H(n, W)$, \eqref{eq:HnwCutNorm} implies that
\begin{equation}\label{eq:wHnGnwCutDistance}
\dkcut\bigl(G(n, W), W_{H_n}\bigr) = \max_{S, T \subseteq [n]} \dfrac{1}{n^2} \sum_{i=1}^k \Bigl\lvert Z^{(S, T)}_i - \mathbb{E}Z^{(S, T)}_i \Bigr\rvert.
\end{equation}
Let us call a triple~$(S, T, i)$ \emph{bad} if $\abs{Z^{(S, T)}_i - \mathbb{E}Z^{(S, T)}_i} \geq \varepsilon n^2/4k$.  By~\eqref{eq:wHnGnwCutDistance}, the union bound and Remark~\ref{re:CutNormDisjoint}, if there are no bad triples, then $\dkcut(G(n, W), W_{H_n}) < \varepsilon$.

Changing the colour of a single edge changes $Z^{(S, T)}_i$ by~$1$, so the standard bounded differences inequality implies that
\begin{equation*}\label{eq:ZSTiConcentration}
\Prob\Bigl(\Bigl\lvert Z^{(S, T)}_i - \mathbb{E}Z^{(S, T)}_i \Bigr\rvert \geq \varepsilon n^2/4k \Bigr) \leq 2\exp\biggl(-\dfrac{\varepsilon^2 n^4}{8k^2 \abs{S}\abs{T}}\biggr) \leq 2\exp\biggl(-\dfrac{\varepsilon^2 n^2}{8k^2}\biggr).
\end{equation*}
There are at most~$3^n$ pairs of disjoint non-empty sets~$(S, T)$, so the probability that there is a bad triple is at most~$\exp(-c\varepsilon^2 n^2)$ for some~$c > 0$.

It follows from the Borel--Cantelli lemma that almost surely $\dkcut(G(n, W), W_{H_n}) < \varepsilon$ for all but finitely many~$n$.  Because $\varepsilon$ is arbitrary, the claim follows.
\end{proof}

\begin{claim}\label{cl:EquivalentGnConvergence}
If $W$ and $U$ are equivalent $k$-graphons, then almost surely
\[
\dkcut\bigl(G(n, W), G(n, U)\bigr) \to 0
\]
as $n \to \infty$.
\end{claim}

\begin{proof}
Observe that~\eqref{eq:wRandomDistribution} implies that for all $n$, $G(n, W)$ and $G(n, U)$ have the same distribution.  We may couple $G(n, W)$ and $G(n, U)$, again generating them by the same sequence of random variables $X_1$, \dots,~$X_n$.

Given disjoint sets $S$,~$T \subseteq [n]$ and a colour~$i \in [k]$, let $Z^{(S, T)}_i(W)$ and $Z^{(S, T)}_i(U)$ be as in the proof of Claim~\ref{cl:wHnGnwConvergence}.  It follows that
\begin{equation}\label{eq:GnwGnuCutDistance}
\dkcut\bigl(G(n, W), G(n, U)\bigr) = \max_{S, T \subseteq [n]} \dfrac{1}{n^2} \sum_{i=1}^k \Bigl\lvert Z^{(S, T)}_i(W) - Z^{(S, T)}_i(U) \Bigr\rvert.
\end{equation}
Because $G(n, W)$ and $G(n, U)$ are both generated by the~$X_j$, for all $S$,~$T$ and~$i$ we have $\mathbb{E}Z^{(S, T)}_i(W) = \mathbb{E}Z^{(S, T)}_i(U)$.  It follows that the right-hand side of~\eqref{eq:GnwGnuCutDistance} is at most
\begin{equation*}
\max_{S, T \subseteq [n]} \dfrac{1}{n^2} \sum_{i=1}^k \Bigl(\Bigl\lvert Z^{(S, T)}_i(W) - \mathbb{E}Z^{(S, T)}_i(W)\Bigr\rvert + \Bigl\lvert  \mathbb{E}Z^{(S, T)}_i(U) - Z^{(S, T)}_i(U) \Bigr\rvert \Bigr).
\end{equation*}
By the same argument as in the proof of Claim~\ref{cl:wHnGnwConvergence}, almost surely $\dkcut(G(n, W), G(n, U)) < 2\varepsilon$ for all but finitely many~$n$.  This proves the claim.
\end{proof}

Claims \ref{cl:wHnConvergence}, \ref{cl:wHnGnwConvergence} and~\ref{cl:EquivalentGnConvergence} and the triangle inequality imply that
\[
\deltakcut(W, U) = 0,
\]
which is what we wanted.
\end{proof}

For later use, we record a corollary of Claims \ref{cl:wHnConvergence} and~\ref{cl:wHnGnwConvergence} in the proof of Theorem~\ref{th:deltakcut} above.

\begin{corollary}\label{cor:wRandomConvergenceCutNorm}
If $W$ is a $k$-graphon, then almost surely
\[
\deltakcut(G(n, W), W) \to 0
\]
as $n \to \infty$.
\end{corollary}

Now we are ready to prove our main result, Theorem~\ref{th:compactnessFinite}.

\begin{proof}[Proof of Theorem \ref{th:compactnessFinite}.]
Let~$\sequence{W}$ be a sequence of $k$-graphons that converges in $\deltakcut$.  Lemma~\ref{le:CountingFinite} implies that for every $B^{\infty}$-decorated graph~$F$, the sequence $\{t(F, W_n)\}$ converges.  Thus, by Theorem~\ref{th:ExistenceFinite}, there exists a $K$-graphon~$W$ such that for every $\R^k$-decorated graph~$F$, $t(F, W_n) \to t(F, W)$.  Finally, Theorem~\ref{th:deltakcut} implies that $\deltakcut(W_n, W) \to 0$ as $n \to \infty$, as claimed.
\end{proof}

\section{From containers to the entropy of graph limits}\label{section: containers to graphon entropy}
As mentioned earlier, Hatami, Janson and Szegedy~\cite{HJS} used results on the entropy of graph limits to prove Theorem~\ref{theorem: alekseevbollobasthomason} and to give counting and characterization results for hereditary properties of graphs.  In this and the next section, we shall prove multicoloured generalizations of their results.

\begin{definition}
We define the \emph{$k$-ary entropy} of a $k$-dimensional vector~$\vector{P}=(p_1, p_2, \ldots, p_k)$ to be 
\[h_k(\vector{P}):= \sum_{c \in[k]} -p_c\log_k p_c.\]
The \emph{entropy} of a $k$-decorated graphon $W$ is 
\[\Ent(W):= \iint_{[0,1]^2} h_k \bigl(W_k(x,y)\bigr)dA.\]
\end{definition}
Note that $0\leq \Ent(W)\leq 1$. For $k=2$ our definition of decorated graphon entropy coincides with that of Hatami, Janson and Szegedy.  Furthermore, if $W \cong W'$, then $\Ent(W) = \Ent(W')$.  Thus, if $\Gamma$ is a limit of $k$-decorated graphs and $W$ is a $k$-graphon that represents $\Gamma$, then we may set $\Ent(\Gamma) = \Ent(W)$.

Given a property~$\mathcal{P}$ of ($k$-decorated) graphs, we let $\limitsx{\mathcal{P}}$ denote its completion under the cut norm.  We also let $\maxlimits{\mathcal{P}}$ denote the set of elements of~$\limitsx{\mathcal{P}}$ of maximum entropy. Hatami, Janson and Szegedy related the speed of a hereditary property~$\mathcal{P}$ to the maximum entropy of graphons in its completion.

In this section, we prove a $k$-decorated version of the main counting result from~\cite{HJS}.

\begin{theorem}\label{theorem: HJS for decorated graphons}
If\/ $\mathcal{P}$ is an (order-)hereditary property of $k$-colourings, then
\begin{enumerate}[(i)]
\item $\pi(\mathcal{P})=\max_{W \in \limitsx{\mathcal{P}}}\Ent(W)$ and
\item We have
\[\lim_{n \to \infty} \dfrac{\log_k \abs{\PP_n}}{\binom{n}{2}} = \max_{W \in \limitsx{\PP}} \Ent(W).\]
\end{enumerate}
\end{theorem}

Note that Part~(ii) corresponds to \cite[Theorem~1.5]{HJS}.  

In this section, we use our multicolour container results from Section~\ref{section: containers} and our compactness result from Section~\ref{section: cut metric for graphons} to prove Theorem~\ref{theorem: HJS for decorated graphons}.  In Section~\ref{section: decorated generalization of HJS}, we shall use Theorem~\ref{th:compactnessFinite} and counting results for Szemer\'edi partitions to give a second proof of part~(ii) of Theorem~\ref{theorem: HJS for decorated graphons}, as well as other counting and characterization results from~\cite{HJS}.

Recall from Definition~\ref{definition: decorated graphons} that a $k$-decorated graphon~$W$ is a symmetric, measurable function that assigns a probability distribution on $[k]$ to each point in $\oi^2$. We may use this to define random templates and $k$-colourings from a $k$-decorated graphon. 


\begin{definition}[$W$-random templates and colourings]
Given a set of $n$ points $X_1$, $X_2$, \dots,~$X_n$ from $[0,1]$ and a $k$-decorated graphon~$W$, we may define a $k$-colouring template for $K_n$, $t_W[X_1, \ldots, X_n]$, by setting $t(ij)=\{c\in [k]: \ W(X_i,X_j)_c>0\}$. Further we may define a random $k$-colouring $c_W[X_1, \ldots, X_n]$ by setting $c(ij)$ to be a random colour from $[k]$ drawn according to the probability distribution given by $W(X_i, xX_j)$.

Finally, we define the \emph{$W$-random template}~$t_W(n)$ and the \emph{$W$-random colouring}~$c_W(n)$ by selecting the points $X_1$, $X_2$, \dots,~$X_n$ uniformly at random from $[0,1]$, and then taking the resulting (induced) $k$-colouring template and random $k$-colouring respectively.
\end{definition}
The definitions of $t_W(n)$~and~$c_W(n)$ are very similar to the definitions of $H(n, W)$~and~$G(n, W)$, respectively, given in Section~\ref{se:DecoratedDefs}.  The only difference is that, for each $i$ and $j$, $t(ij)$ is uniformly distributed on $\supp(W(X_i, X_j))$, which does not in general hold for the edges of~$H(n, W)$.

Our $W$-random templates and colourings give us a way of going from $k$-decorated graphons to templates for $k$-colourings of~$E(K_n)$/$k$-colourings of~$E(K_n)$. We can also go in the other direction: given a template $t$ for a $k$-colouring of~$K_n$, we may define the step-function $W_t$ by dividing $[0,1)$ into intervals $I_i=[(i-1)/n,i/n)$, $1\leq i\leq n$ and on each tile $I_i\times I_j$ setting $W_t$ to be constant and equal to $\frac{1}{\vert t(ij)\vert} (\mathbf{1}_{c\in t(ij)})_{c\in [k]}$.
By viewing a $k$-colouring $c$ of $E(K_n)$ as a (zero entropy) template, we may in the same way obtain from it a $k$-decorated graphon $W_c$. Thus we may go in a natural way from talking of a property of colourings to a property of decorated graphons.


Note that for all $k$ and~$n$, the quantity~$\binom{n}{2}\Ent(W)$ generalizes our notion of the entropy of a $k$-colouring template for $E(K_n)$ via our identification of templates with $k$-colouring graphons.  Furthermore, $\Ent(W)$ computes the expected value of the discrete $k$-ary entropy of the $W$-random colouring model~$c_W(n)$ which was  introduced above.

With these definitions and the container results of Section~\ref{section: containers} in hand, we are only a few lemmas away from the proof of Theorem~\ref{theorem: HJS for decorated graphons}.
\begin{definition}
A \emph{measurable partition} of $[0,1]$ is a partition $\mathcal{S}$ of $[0,1]$ into finitely many measurable sets $\{S_1, S_2, \ldots, S_n\}$. The \emph{conditional expectation of a $k$-decorated graphon with respect to $\mathcal{S}$}, $\mathbb{E}[W\vert \mathcal{S}]$, is the step-function on the tiling defined by $\mathcal{S}\times \mathcal{S}$, with $\mathbb{E}[W\vert \mathcal{S}]$ equal on each tile~$S_i\times S_j$ to the average value of~$W$ over that tile.
\end{definition}

\begin{lemma}\label{lemma: generalization of 3.3 in HJS}
\begin{enumerate}[(i)]
    \item For every $k$-decorated graphon $W$ and every measurable partition\/ $\mathcal{S}$ of $[0,1]$, we have $\Ent(\mathbb{E}[ W \vert \mathcal{S}]) \geq \Ent(W)$.
    \item $\Ent(\cdot)$ is lower semicontinuous on the space of all $k$-decorated graphons: if $W_m \to W$ in the cut norm, then $\limsup_{m\to \infty} \Ent(W_m) \leq \Ent(W)$.
\end{enumerate}
\end{lemma}

\begin{remark}\label{rk:EntropyNotContinuous}
It is shown in~\cite{HJS} that $\Ent(\cdot)$ is not continuous.
\end{remark}

Part~(ii) of Lemma~\ref{lemma: generalization of 3.3 in HJS} was proved for ordinary graphons in Lemma~2.1 of~\cite{ChatterjeeVaradhan11}.  Hatami, Janson and Szegedy~\cite{HJS} gave a different proof, which we follow here.

\begin{proof}[Proof of Lemma~\ref{lemma: generalization of 3.3 in HJS}]
Part (i) is immediate from Jensen's inequality and the concavity of $x\mapsto -x\log_k x$.

For part (ii), suppose $\|W_m- W\|_{\square_k} \to 0$. In particular,
\eqref{eq:dkcutDef} implies that for each $i \in [k]$ we must have convergence of the $i$-coloured graphons, i.e., $\|W^i_{m}-W^i\|_{\square_k} \to 0$. For every $j \geq 1$, let $\overline{\mathcal{S}_j}$ denote the partition of $[0,1)$ into $I_1$, \dots,~$I_j$.
By definition of the cut norm, for each $i$,~$j$ the conditional expectation~$\mathbb{E}[W^i_m \vert  \mathcal{S}_j]$ converges to $\mathbb{E}[W^{i} \vert \mathcal{S}_j]$ almost everywhere as $m \to \infty$ and thus 
$\mathbb{E}[W_m \vert  \mathcal{S}_j]\rightarrow \mathbb{E}[W \vert  \mathcal{S}_j]$
almost everywhere in $[0,1]^2$. Since our entropy function is bounded on the space of decorated graphons, we can apply the dominated convergence theorem to deduce that 
\[\lim_{m \to \infty} \Ent(\mathbb{E}[W_{m} \vert \mathcal{S}_j]) = \Ent(\mathbb{E}[W \vert \mathcal{S}_j]).\]
Further we have $\mathbb{E}[W\vert \mathcal{S}_j]\rightarrow W$ almost everywhere in $[0,1]^2$ as $j\rightarrow \infty$. By dominated convergence again, 
\[\lim_{j\rightarrow \infty} \Ent(\mathbb{E}[W \vert \mathcal{S}_j])=\Ent(W).\]
Finally applying (i) we have 
\[\limsup_{m\to\infty} \Ent(W_m) \leq \limsup_{j\rightarrow \infty} \limsup_{m\to\infty}\Ent(\mathbb{E}[W^j_{m} \vert \mathcal{S}_j]) = \limsup_{j\rightarrow\infty} \Ent(\mathbb{E}[W \vert \mathcal{S}_j])=\Ent(W). \qedhere\]
\end{proof}

Let $\mathcal{P}$ be a hereditary property of $k$-colourings. By Corollary~\ref{corollary: containers for arbitrary hereditary properties}, for every $j\in \N$, there exists $n_j$ and $C_0(j)$ such that for all $n\geq n_j$ there exists a container family $\mathcal{T}^j_n$ for $\mathcal{P}_n$ satisfying parts~(i)--(iv) of the corollary with $\varepsilon=1/j$ and $m=j$. We define a sequence of containers~$\mathcal{C}=(\mathcal{C}_n)_{n\geq n_1}$ for $(\mathcal{P}_n)_{n\geq n_1}$ by setting
\[j_{\star}(n)=\max\{j \, : \, n_j \leq n\}\]
and letting $\mathcal{C}_n=\mathcal{T}^{j_{\star}(n)}_n$. Finally, let $\limitsx{\mathcal{C}}$ denote the collection of decorated graphons obtained as limits of sequences of templates from $\mathcal{C}$. 
\begin{lemma}\label{lemma: limit of containers are in hat P}
$\limitsx{\mathcal{C}}\subseteq \limitsx{\mathcal{P}}$.
\end{lemma}
\begin{proof}
Let $(t_n)_{n\in\N}$ be a sequence of templates from $\mathcal{C}_n$ converging to a decorated graphon $W$. For any $m\in \N$, consider the random colouring $c_W(m)$ of $E(K_m)$ obtained from $W$. By a result of Lov\'asz and Szegedy (Corollary 2.6 in~\cite{LS06}),
for each $i\in[k]$, the graph of $i$-coloured edges in $c_W(m)$ converges almost surely to $W^i$ as $m\rightarrow \infty$, and thus the sequence of $k$-colourings~$(c_W(m))_{m\in\N}$ itself converges to $W$ as $m\rightarrow \infty$ almost surely. For each fixed $m$, let $X_m$ be a uniformly-random $m$-set of~$V(K_n)$ and observe that
\[\mathbb{P}(c_W(m)\notin \mathcal{P}_m) \leq \lim_{n\to\infty} \mathbb{P}\bigl((t_n)_{\mid X_m} \nsubseteq \PP_n\bigr).\]
By our construction of~$C_n$, the limit on the right-hand side is~$0$, whence almost surely $c_W(m)\in \mathcal{P}_m$ for each fixed $m$. We have thus exhibited a sequence of colourings (almost surely) in $\mathcal{P}$ which (almost surely) converges to $W$, showing that $W\in\limitsx{\mathcal{P}}$ and $\limitsx{\mathcal{C}}\subseteq \limitsx{\mathcal{P}}$ as claimed.
\end{proof}

\begin{corollary}\label{corollary: entropy of container families at most max entropy in hat P} 
$\limsup_{n\rightarrow \infty} \max_{t\in\mathcal{C}_n}\Ent(t)/\binom{n}{2} \leq \max_{W\in \limitsx{\mathcal{P}}}\Ent(W).$
\end{corollary}
\begin{proof}
Immediate from the lower semicontinuity of $\Ent(\cdot)$ and the inclusion $\limitsx{\mathcal{C}}\subseteq\limitsx{\mathcal{P}}$ we have just established. Note the maximum in the statement of the corollary exists as $\limitsx{\mathcal{P}}$ is a closed set in the compact space of decorated graphons.
\end{proof}

\begin{theorem}\label{theorem: entropy of container families is equal to max entropy of hat P}
$\lim_{n\rightarrow \infty} \max_{t\in \mathcal{C}_n}\Ent(t)/\binom{n}{2} = \max_{W\in \limitsx{\mathcal{P}}} \Ent(W)$.
\end{theorem}

\begin{proof}
Set
\[\alpha= \liminf_{n\rightarrow \infty} \max_{t\in \mathcal{C}_n}\Ent(t)/\binom{n}{2}\]
and
\[\beta= \max_{W \in \limitsx{\mathcal{P}}} \Ent(W).\]
By Lemma~\ref{lemma: limit of containers are in hat P}, $\alpha \leq \beta$. We show that in fact we have equality.

Let $W$ be an entropy maximizer in $\limitsx{\mathcal{P}}$. For every $n\in\N$, by linearity of expectation there exists a choice of $n$ points $X_1$, $X_2$, \dots,~$x_n$ from $[0,1]$ such that
\[\Ent\bigl(t_W[X_1, X_2, \ldots, X_n]\bigr)\geq \beta \binom{n}{2}.\] 
Furthermore, as $W\in\limitsx{\mathcal{P}}$, almost surely $\langle t_W[X_1, \ldots, X_n] \rangle \subseteq \mathcal{P}_n$, implying $\vert \mathcal{P}_n\vert \geq k^{\beta\binom{n}{2}}$ for every $n\in \N$. Now let $(n_i)_{i\in\N}$ be a subsequence such that $\lim_{i\rightarrow \infty}\max_{t\in \mathcal{C}_{n_i}}\Ent(t)/\binom{n_i}{2}= \alpha$. By Corollary~\ref{corollary: speed of arbitrary hereditary properties} and the construction of the sequence of container families~$(\mathcal{C}_n)_{n\in\N}$, we have 
\[k^{\beta\binom{n_i}{2}}\leq \vert \mathcal{P}_{n_i}\vert \leq k^{\bigl(\alpha +o(1)\bigr)\binom{n_i}{2}},\]
whence $\alpha\geq  \beta$, as desired. Together with Corollary~\ref{corollary: entropy of container families at most max entropy in hat P}, this shows 
\[\lim_{n\rightarrow \infty}\max_{t \in \mathcal{C}_n} \Ent(t)/\binom{n}{2}=\max_{W\in \limitsx{\mathcal{P}}} \Ent(W). \qedhere \] 
\end{proof}
We arrive at last at the proof of Theorem~\ref{theorem: HJS for decorated graphons}.
\begin{proof}[Proof of Theorem~\ref{theorem: HJS for decorated graphons}]
By Theorem~\ref{theorem: entropy of container families is equal to max entropy of hat P}, we have $\max_{W \in \limitsx{\mathcal{P}}}\Ent(W)=\lim_{n\rightarrow \infty} \max_{t\in \mathcal{C}_n} \Ent(t)/\binom{n}{2}$. The monotonicity established in Proposition~\ref{proposition: entropy density} implies  $\pi(\mathcal{P})\binom{n}{2}\leq \ex(n, \mathcal{P})$ for every $n$. Further, by construction of~$(\mathcal{C}_n)_{n\in\N}$, and in particular properties (i) and (ii) in Corollary~\ref{corollary: containers for arbitrary hereditary properties}, for $n\geq n_1$, we have 
\[\ex(n, \mathcal{P})\leq \max_{t\in \mathcal{C}_n} \Ent(t) \leq \biggl(\pi(\mathcal{P})+\frac{1}{j_{\star}(n)}\biggr)\binom{n}{2}.\]
As $j_{\star}(n)\rightarrow \infty$, this implies 
\[\max_{W \in \limitsx{\mathcal{P}}}\Ent(W)=\lim_{n\rightarrow \infty} \max_{t\in \mathcal{C}_n} \Ent(t)/\binom{n}{2}=\pi(\mathcal{P}),\]
showing the extremal entropy density and the maximum decorated graphon entropy of an (order-) hereditary property are the same, which is part (i) of our theorem.

The counting result, part (ii), is immediate from part (i) and Corollary~\ref{corollary: speed of arbitrary hereditary properties}.
\end{proof}

\section{Entropy of \texorpdfstring{$k$}{k}-decorated graphons}\label{section: decorated generalization of HJS}


As mentioned earlier, in this section, we shall use the graph limit results of the previous sections to generalise  results from~\cite{HJS} to the case of $k$-decorated graphs.


The first theorem, a counting result for arbitrary properties of $k$-decorated graphs, corresponds to~\cite[Theorem~1.1]{HJS}.

\begin{theorem}\label{th:SpeedEntropyIneq}
If $\QQ$ is a class of $k$-coloured graphs, then
\begin{equation}\label{eq:SpeedEntropyIneq}
\limsup_{n \to \infty} \dfrac{\log_k \abs{\QQ_n}}{\binom{n}{2}} \leq \max_{\Gamma \in \limitsx{\QQ}} \Ent(\Gamma).
\end{equation}
\end{theorem}

As mentioned earlier, we shall also give a second proof of Theorem~\ref{theorem: HJS for decorated graphons}(ii) using graph limits.

%

Finally, whenever $\QQ$ is such that equality holds in~\eqref{eq:SpeedEntropyIneq}, almost every graph in $\QQ$ is close to an element of~$\maxlimits{\QQ}$.  Let $\graphs^k$ denote the class of unlabelled $k$-coloured graphs and let $\labeled^k_n$ denote the class of labelled $k$-coloured graphs with vertex set~$[n]$.  This result corresponds to~\cite[Theorem~1.6]{HJS}.

\begin{theorem}\label{th:maxlimits}
Suppose that $\QQ$ is a class of $k$-coloured graphs such that equality holds in~\eqref{eq:SpeedEntropyIneq}.
\begin{enumerate}[(i)]
\item If $G_n \in \graphs^k_n$ is a uniformly random unlabelled element of~$\QQ_n$, then $\deltaKcut(G_n, \maxlimits{\QQ}) \pto 0$.
\item If $G_n \in \labeled^k_n$ is a uniformly random labelled element of~$\QQ^L_n$, then $\deltaKcut(G_n, \maxlimits{\QQ}) \pto 0$.
\end{enumerate}
\end{theorem}

Theorem~\ref{th:maxlimits} follows from Theorem~\ref{th:SpeedEntropyIneq}.  We do not give the proof, which is essentially identical to the proof of~\cite[Theorem~1.6]{HJS}.

\begin{remark}\label{rk:maxlimitsForHereditaryProperties}
Note that Theorem~\ref{theorem: HJS for decorated graphons}(ii) implies that Theorem~\ref{th:maxlimits} holds whenever $\QQ$ is a hereditary class.
\end{remark}

\subsection{A Weak Regularity Lemma for \texorpdfstring{$k$}{k}-graphons and other preliminary results}\label{se:entropyPrelims}

Here we prove assorted results about $k$-coloured graphs that are known to hold for ordinary graphons, and thus are used without proof in~\cite{HJS}.

First, Lov\'asz and Szegedy~\cite{LS:analyst} proved an analytic version of the Weak Regularity Lemma of Frieze and Kannan~\cite{FK99}.  

\begin{lemma}\label{le:WeakRegGraphons}
Let $\varepsilon > 0$.  There exists $m = 2^{O(1/\varepsilon^2)}$ such that if $W$ is a graphon, then there exists a partition\/ $\mathcal{S} = \{S_1, \dots, S_m\}$ and a step-function~$U$ that is constant on the cells of\/ $\mathcal{S}^2$ such that
\[
\cutnorm{W - U} < \varepsilon.
\]
Moreover, at the cost of increasing the implied constant in the bound on $m$, we may assume
that the elements of\/ $\mathcal{S}$ have equal measure and
that $U = \mathbb{E}[W\vert \mathcal{S}]$.
\end{lemma}

Lemma~\ref{le:WeakRegGraphons} easily implies a Weak Regularity Lemma for $k$-graphons.

\begin{lemma}\label{le:WeakRegEffective}
Let $k \in \N$.  For every $\varepsilon > 0$, there exists $m = m(\varepsilon, k)$, with
\[
m = 2^{O(1/\varepsilon^2)},
\]
such that if $W$ is a $k$-graphon, then there exists a partition\/ $\mathcal{S} = \{S_1, \dots, S_m\}$ and a step-function~$U$ that is constant on the cells of\/ $\mathcal{S}^2$ such that
\[
\dkcut(W, U) < \varepsilon.
\]
Moreover, at the cost of increasing the implied constant in the bound on $m$, we may assume
that the elements of\/ $\mathcal{S}$ have equal measure and
that $U = \mathbb{E}[W\vert \mathcal{S}]$.
\end{lemma}

\begin{proof}
Recall from~\eqref{eq:omegaiDef} that for each $i \in [k]$, $W_i$ is the graphon given by $W_i(x, y) = \Prob\bigl(W(x, y) = i \bigr)$.  By Lemma~\ref{le:WeakRegGraphons}, for each $i$, there is a partition~$\SSS_i$ of~$\oi$ into at most $2^{O(1/\varepsilon^2)}$ sets of equal measure such that
\begin{equation*}
\cutnorm{W_i - (W_i)_{\SSS_i}} < \dfrac{\varepsilon}{k}.
\end{equation*}
Let $M \leq 2^{c/\varepsilon^2}$ denote the maximum number of steps in any $\SSS_i$.  Let $\partition$ be a common refinement of the partitions defined by the $\SSS_i$ and observe that we may assume that
\[
\abs{\partition} \leq M k = 2^{O(1/\varepsilon^2)}.
\]
Furthermore, it is easy to see that we still have
\begin{equation}\label{eq:coordinateWeakApprox}
\cutnorm{W_i - \mathbb{E}[W_i\vert \partition]} < \dfrac{\varepsilon}{k}.
\end{equation}
By construction, for all~$(x, y)$, we have $\sum_{i=1}^k \mathbb{E}[W_i\vert \partition](x, y) = 1$, so the $\mathbb{E}[W_i\vert \partition]$ define a $k$-graphon on at most~$M k$ steps.  (In fact, this $k$-graphon is none other than~$\mathbb{E}[W\vert \partition]$.)  Finally, by~\eqref{eq:coordinateWeakApprox}, we have
\[
\dkcut(W, \mathbb{E}[W\vert \partition]) < \varepsilon,
\]
as claimed.
\end{proof}

The next result concerns convergence of $W$-random graphs.  It is an immediate consequence of Corollary~\ref{cor:wRandomConvergenceCutNorm} and Theorem~\ref{th:deltakcut}.  (The corresponding statement for graphons was proved in~\cite{BCLSV08,LS06}.)

\begin{theorem}\label{th:wRandomConvergence}
If $W$ is a $k$-graphon and $G(n, W)$ is a sequence of $W$-random graphs, then
\[
G(n, W) \to W
\]
almost surely.
\end{theorem}

Let $\PP$ be a hereditary property of graphs and let $W \in \limitsx{\PP}$.   It is not hard to show that the $W$-random graph~$G(n, W) \in \PP$ almost surely (see, e.g., \cite[Theorem~3.2]{Jan16}).  Here, we do the same for $K$-decorated graphs. 

Let $F$ be a $K$-decorated graph on $n$ vertices.  Recall the definition of the $C[K]$-decorated graph~$G_F$ from Section~\ref{se:DecoratedDefs}.

\begin{lemma}\label{le:wRandomInHeredProp}
Let $\QQ$ be a hereditary property of $K$-decorated graphs.  If $W$ is a $K$-graphon, then exactly~one of the following holds:
\begin{itemize}
\item We have $W \in \limitsx{\QQ}$ and for all $n \geq 1$, $G(n, W) \in \QQ$ almost surely.
\item We have $W \notin \limitsx{\QQ}$ and $\Prob(G(n, W) \in \QQ) \to 0$ as $n \to \infty$.
\end{itemize}
\end{lemma}

\begin{proof}
Let $W \in \limitsx{\QQ}$ and let $(H_n) \subseteq \QQ$ be a sequence of $K$-decorated graphs such that $H_n \to W$.  Suppose that $F \notin \QQ$.  Because $\QQ$ is hereditary, it follows from \eqref{eq:homphi} and~\eqref{eq:homdensity} that $t(G_F, H_n) = 0$.  Then, by hypothesis, $t(G_F, W) = \lim_{n \to \infty} t(G_F, H_n) = 0$.

So, if $t(G_F, W) > 0$ then $F \in \QQ$.  By~\eqref{eq:wRandomDistribution}, if $\Prob(G(n, W) = F) >0$, then $t(G_F, W) > 0$ and so $F \in \QQ$.  It follows that for each $n$, $G(n, W) \in \QQ$ almost surely, which means that almost surely the statement holds for all $n$.

Now suppose that $W \notin \limitsx{\QQ}$.  Let $\Gamma$ be the graph limit that $W$ represents.  Because $W\notin \limitsx{\QQ}= \overline{\QQ} \cap \limitsx{\graphs^K}$, we have $\Gamma \notin \overline{\QQ}$, that is, there exists an open neighbourhood~$X$ of~	$\Gamma$ in $\overline{\graphs^K}$ such that $X \cap \QQ = \emptyset$.  By Theorem~\ref{th:wRandomConvergence}, $G(n, \Gamma) \to \Gamma$ almost surely, and hence in probability.  It follows that $\Prob(G(n, \Gamma) \in \QQ) \to 0$ as $n \to \infty$, as claimed.
\end{proof}

\begin{remark}
Observe that in the second part of the proof of Lemma~\ref{le:wRandomInHeredProp}, we did not use the assumption that $\QQ$ was hereditary.  If $\family$ is an arbitrary family of $K$-decorated graphs, let $H(\family)$ denote the union of $\family$ and all induced subgraphs of elements of~$\family$.  Clearly, $H(\family)$ is a hereditary class.  The proof of Lemma~\ref{le:wRandomInHeredProp} shows that for any family~$\family$ and any $K$-graphon~$W$, we have $W \in \limitsx{\family}$ if and only if for all $n \geq 1$, we have $G(n, W) \in H(\family)$ almost surely.
\end{remark}

%
%

\begin{lemma}\label{le:randomWeightedCutDistance}
If $H$ is a\/ $\measures([k])$-decorated graph of order~$n$ and $G(H)$ is obtained by sampling from $H$, then with probability at least~$1 - e^{-n}$, we have
\[
\dkcut(G(H), H) \leq \dfrac{10}{\sqrt{n}}.
\]
\end{lemma}

We omit the proof, which is similar to the proof of Claim~\ref{cl:wHnGnwConvergence} in the proof of Theorem~\ref{th:deltakcut}.


The next result is similar to~\cite[Lemma~9.29]{LovaszBook}.  We omit the proof.

\begin{lemma}\label{le:relabel}
If $G_1$ and $G_2$ are unlabelled $k$-decorated graphs on $n$ vertices, there exists a relabelling~$\widetilde{G_1}$ of~$G_1$ such that
\[
\dkcut(\widetilde{G_1}, G_2) \leq \deltakcut(G_1, G_2) + \dfrac{17}{\sqrt{\log_2 n}}.
\]
\end{lemma}

If $G$ is a $k$-decorated graph of order~$n$ and $\partition$ is a partition of~$V(G)$, the averaged graph~$G_{\partition}$ is a $\measures([k])$-decorated graph on $n$ vertices obtained by averaging $G$ over each product of cells of~$\partition$.


We also need an important fact about entropy and $W$-random graphs.  The result for graphons is given in~\cite{Aldous,Jan13}.  The proof given in~\cite[Appendix~D]{Jan13} for graphons goes through essentially unchanged in the general case, so we omit it.

\begin{theorem}\label{th:wRandomEntropy}
If $W$ is a $k$-graphon, then
\[
\lim_{n \to \infty} \dfrac{\Ent \bigl(G(n, W)\bigr)}{\binom{n}{2}} = \Ent(W).
\]
\end{theorem}

Given a $k$-graphon~$W$, let $\average{W}_n$ denote the $\measures([k])$-random graph obtained by averaging $W$ over squares of sidelength~$1/n$.  In other words,
\begin{equation}\label{eq:averageDef}
\average{W}_{n_{ij}} = n^2 \iint_{I_i \times I_j} W(x, y)\,dxdy.
\end{equation}

\subsection{Lemmas}\label{se:entropyLemmas}

Here we state and prove several lemmas about the number of graphs on $n$ vertices that are close in cut distance to a given $k$-graphon.  This section corresponds to Section~4 of~\cite{HJS}.  As a rule, we shall only sketch the proofs and will mostly highlight the places where our arguments differ from those in~\cite{HJS}.

Let $\QQ$ be any class of $k$-decorated graphs.  For any $n \geq 1$, we have
\begin{equation}\label{eq:labelledUnlabelled}
\abs{\QQ_n} \leq \abs{\QQ_n^L} \leq n! \abs{\QQ_n}.
\end{equation}

Given an integer~$n$, $\delta > 0$, and a $k$-graphon~$W$, we define
\begin{equation*}\label{eq:numberdKcut}
\NdKcut{W} = \bigl\lvert\{G \in \labeled^k_n \, : \, \dKcut(G, W) \leq \delta\}\bigr\rvert
\end{equation*}
and
\begin{equation*}\label{eq:numberdeltaKcut}
\NdeltaKcut{W} = \bigl\lvert\{G \in \labeled^k_n \, : \, \deltaKcut(G, W) \leq \delta\}\bigr\rvert.
\end{equation*}
It is trivial that 
\begin{equation}\label{eq:Ninequality}
\NdKcut{W} \leq \NdeltaKcut{W}.
\end{equation}

First, we relate $\NdKcut[\cdot]{W}$ and $\NdeltaKcut[\cdot]{W}$.  To do this, we show that if $k$-decorated graph~$G$ is close to a $k$-graphon~$W$ in terms of $\Kcutnorm{\cdot}$, then some rearrangement of~$G$ is close to~$W$ in terms of $\deltaKcut$.

\begin{lemma}\label{le:dKcutRearrangement}
Let $W$ be a $k$-graphon.  If $G \in \labeled^k_n$, then there is a rearrangement~$\widetilde{G}$ of~$G$ such that
\[
\dKcut(\widetilde{G}, W) \leq \deltaKcut(G, W) + 2\dKcut(W, \average{W}_n) + \dfrac{18}{\sqrt{\log_2 n}}.
\]
\end{lemma}

The proof is very similar to the corresponding argument in~\cite[Lemma~4.1]{HJS}.  We include the proof despite the similarities because it uses some of the results from Section~\ref{se:entropyPrelims}.

\begin{proof} We observe that if $n < 2^{20}$, then the result is trivial, because we always have $\dKcut(\widetilde{G}, W) \leq 1$.

Recall the definition of~$\average{W}_n$ from~\eqref{eq:averageDef} and let $G(\average{W}_n)$ be the random $k$-decorated graph obtained by sampling from $\average{W}_n$.  By Lemma~\ref{le:randomWeightedCutDistance}, with probability tending to~$1$, we have
\[
\dKcut(G(\average{W}_n, \average{W}_n) \leq \dfrac{10}{\sqrt{n}},
\]
so we may choose a realization~$G'$ of~$G(\average{W}_n)$ such that
\begin{equation}\label{eq:randomAverageDistance}
\dKcut(G', \average{W}_n) \leq \dfrac{10}{\sqrt{n}}.
\end{equation}
It follows that
\begin{equation}\label{eq:dGGprime}
\deltaKcut(G, G') \leq \deltaKcut(G, W) + \dKcut(W, \average{W}_n) + \dKcut(\average{W}_n, G') \leq \deltaKcut(G, W) + \dKcut(W, \average{W}_n) + \dfrac{10}{\sqrt{n}}.
\end{equation}
By Lemma~\ref{le:relabel}, we may permute $V(G)$ to obtain a graph~$\widetilde{G}$ such that
\begin{equation}\label{eq:permuteDistance}
\dKcut(\widetilde{G}, G') \leq \deltaKcut(G, G') + \dfrac{17}{\sqrt{\log_2 n}}.
\end{equation}
Now, by \eqref{eq:randomAverageDistance}--\eqref{eq:permuteDistance},
\begin{align*}
\dKcut(\widetilde{G}, W) &\leq \dKcut(\widetilde{G}, G') + \dKcut(G', \average{W}_n) + \dKcut(\average{W}_n, W) \\
&\leq \deltaKcut(G, G') + \dfrac{17}{\sqrt{\log_2 n}} + \dfrac{10}{\sqrt{n}} + \dKcut(\average{W}_n, W)\\
&\leq \deltaKcut(G, W) + 2\dKcut(\average{W}_n, W) + \dfrac{17}{\sqrt{\log_2 n}} + \dfrac{20}{\sqrt{n}}.
\end{align*}
The claimed result follows for $n \geq 2^{20}$; as observed above, for smaller~$n$, it is trivial.
\end{proof}

Lemma~\ref{le:dKcutRearrangement} allows us to estimate $\NdeltaKcut{W}$ in terms of~$\NdKcut{W}$.  We omit the proof.

\begin{lemma}\label{le:graphCountIneq}
If $n \geq 1$, $\delta > 0$, and $W$ is a $k$-graphon, then
\[
\NdeltaKcut{W} \leq n! \NdKcut[\delta + \eta_n]{W},
\]
where $\eta_n := 18/\sqrt{\log_2 n} + 2\dKcut(\average{W}_n, W) \to 0$ as $n \to \infty$.
\end{lemma}

Now we estimate both $\NdKcut{W}$ and $\NdeltaKcut{W}$ in terms of~$\Ent(W)$.

\begin{lemma}\label{le:NdKcutEntropy}
Let $W$ be a $k$-graphon.  For any $\delta > 0$,
\[
\liminf_{n \to \infty} \dfrac{\log_k \NdKcut{W}}{\binom{n}{2}} \geq \Ent(W).
\]
\end{lemma}

\begin{proof}[Proof sketch.]
Consider the $W$-random graph~$G(n, W)$ as an element of~$\labeled^k_n$.  By Theorem~\ref{th:wRandomEntropy},
\begin{equation}\label{eq:WRandomEntropy}
\lim_{n \to \infty} \dfrac{\Ent \bigl(G(n, W)\bigr)}{\binom{n}{2}} = \Ent(W).
\end{equation}
Recall from Theorem~\ref{th:wRandomConvergence} that $G(n, W) \to W$ almost surely, and hence in probability.  It follows from Theorem~\ref{th:deltakcut} that $\deltakcut(G(n, W), W) \pto 0$.

The observations above and standard entropy calculations similar to those in~\cite{HJS} imply that
\[
\Ent \bigl(G(n, W)\bigr) \leq \log_k \NdeltaKcut{W} + o(n^2).
\]
By Lemma~\ref{le:graphCountIneq}, there is a sequence~$\eta_n \to 0$ such that
\[
\Ent \bigl(G(n, W)\bigr) \leq \log_k \NdKcut[\delta + \eta_n]{W} + o(n^2).
\]
The claimed result then follows from~\eqref{eq:WRandomEntropy} with $\delta/2$ in place of~$\delta$.
\end{proof}

Recall that the  function~$H_k$ is uniformly continuous on its domain and let $f$ be the modulus of continuity for $H_k$.  In other words, $f$ is a non-decreasing function such that if $X$ and $Y$ are two probability distributions with support in $[k]$ (which we may view as vectors in $\R^k$), we have
\begin{equation}\label{eq:modulus}
\norm{X - Y} < \delta \quad \Longrightarrow \quad \bigl\lvert H_k(X) - H_k(Y) \bigr\rvert < f(\delta).
\end{equation}

\begin{lemma}\label{le:NdKcutUpperBound}
Let $W$ be a $k$-graphon, let $n \geq m \geq 1$, and let $\delta > 0$.  If\/ $\SSS$ is an equipartition of~$[n]$ into $m$ sets, then
\[
\dfrac{\log_k \NdKcut{W}}{n^2} \leq \dfrac{1}{2}\Ent(\mathbb{E}[W\vert \SSS]) + \dfrac{1}{2}f(4m^2 \delta) + 2m^2 k \dfrac{\log_k n}{n^2},
\]
where $f$ is the function in~\eqref{eq:modulus}.
\end{lemma}

\begin{proof}[Proof sketch.]
Let $\SSS = \{V_1, \dots, V_m\}$ and, for each $i$, let $n_i = \abs{V_i}$.  Furthermore, let $\{I_1, \dots, I_m\}$ denote the corresponding partition of $(0, 1]$ obtained by setting $I_i = \bigcup_{s \in V_i} ((s - 1)/n, s/n]$ for each $i$.

Let $G \in \labeled^k_n$ be such that $\dKcut(W_G, W) < \delta$.  For each $(i, j) \in [m]^2$ and each $t \in [k]$, let $e_t(V_i, V_j)$ denote the number of edges of colour~$t$ with one endpoint in $V_i$ and the other in $V_j$.  (When $i = j$, we count edges twice.)  Also, for each $(i, j) \in [m]^2$, let $W_{ij}$ denote the average of~$W$ over $I_i \times I_j$.

Now we count the number of possible choices of~$G$.  First, similar calculations to those in the proof of~\cite[Lemma~4.5]{HJS} imply that for all $t \in [k]$, we have
\begin{equation}\label{eq:partitionEdgeDensity}
\biggl\lvert \dfrac{e_t(V_i, V_j)}{n_i n_j} - W_{ij}(t)\biggr\rvert \leq \dfrac{\delta n^2}{n_i n_j} \leq \delta \biggl(\dfrac{n}{\lfloor n/m \rfloor}\biggr)^{\!2} \leq 4m^2\delta.
\end{equation}
Fix vectors $(e_1(V_i, V_j), \dots, e_k(V_i, V_j))$ that satisfy~\eqref{eq:partitionEdgeDensity} for all $i$ and $j$, and let $N_1$ denote the number of graphs on $[n]$ with these $(e_t(V_i, V_j))_{t \in [k]}$.  Similar calculations to those in~\cite{HJS} show that
\[
\log_k N_1 \leq \dfrac{1}{2} \sum_{i, j = 1}^m n_i n_j H_k\biggl(\dfrac{e_1(V_i, V_j)}{n_i n_j}, \dots, \dfrac{e_k(V_i, V_j)}{n_i n_j}\biggr).
\]
It follows from \eqref{eq:partitionEdgeDensity} and~\eqref{eq:modulus} that
\[
\log_k N_1 \leq \dfrac{1}{2} \sum_{i, j = 1}^m n_i n_j \bigl(H_k\bigl(W_{ij}(1), \dots, W_{ij}(k)\bigr) + f(4m^2 \delta)\bigr).
\]
Dividing both sides by $n^2$ shows that
\[
n^{-2} \log_k N_1 \leq \dfrac{1}{2} \Ent(\mathbb{E}[W\vert \SSS]) + \dfrac{1}{2}f(4m^2 \delta).
\]
Each of the $e_t(V_i, V_j)$ may be chosen in at most~$n^2$ ways, for a total of at most~$n^{2m^2 k}$.  It follows that $\log_k \NdKcut{W} \leq n^{2m^2 k} N_1$ and hence that
\[
n^{-2} \log_k \NdKcut{W} \leq \dfrac{1}{2}\Ent(\mathbb{E}[W\vert \SSS]) + \dfrac{1}{2}f(4m^2 \delta) + 2m^2 k \dfrac{\log_k n}{n^2},
\]
as claimed.
\end{proof}

\begin{lemma}\label{le:NdeltaKcutUpperBound}
Let $W$ be a $k$-graphon.  For any $m \geq 1$, any $\delta > 0$, and any equipartition\/ $\SSS$ of~$\oi$ into $m$ sets,
\[
\limsup_{n \to \infty} \dfrac{\log_k \NdeltaKcut{W}}{\binom{n}{2}} \leq \Ent(\mathbb{E}[W\vert \SSS]) + f(4m^2 \delta),
\]
where $f$ is the function from~\eqref{eq:modulus}.  Hence,
\[
\lim_{\delta \to 0}\limsup_{n \to \infty} \dfrac{\log_k \NdeltaKcut{W}}{\binom{n}{2}} \leq \Ent(\mathbb{E}[W\vert \SSS]).
\]
\end{lemma}

We omit the proof of Lemma~\ref{le:NdeltaKcutUpperBound}, which follows from Lemmas \ref{le:graphCountIneq} and~\ref{le:NdKcutUpperBound} much as in~\cite{HJS}.

The next result, which is the main lemma of this section, says that if $A$ is a closed set of graph limits, then the number of graphons that are close to~$A$ in cut distance is determined (up to the asymptotic value of the logarithm) by the maximum entropy of any graphon in $A$.

\begin{lemma}\label{le:NdeltaKcutSpeed}
Let $A \subseteq \limitsx{\graphs^k}$ be a closed set of limits of $k$-decorated graphs, let $\delta > 0$, and let
\[
\NdeltaKcut{A} = \bigl\lvert\{G \in \labeled^k_n \, : \, \deltaKcut(G, A) \leq \delta\}\bigr\rvert.
\]
We have
\begin{equation}\label{eq:NdeltaKcutSpeed}
\lim_{\delta \to 0}\liminf_{n \to \infty} \dfrac{\log_k \NdeltaKcut{A}}{\binom{n}{2}} = \lim_{\delta \to 0}\limsup_{n \to \infty} \dfrac{\log_k \NdeltaKcut{A}}{\binom{n}{2}} = \max_{W \in A} \Ent(W).
\end{equation}
\end{lemma}

\begin{proof}

First, by Theorem~\ref{th:compactnessFinite}, $A$ is compact.  This and Lemma~\ref{lemma: generalization of 3.3 in HJS} imply that the maximum on the right-hand side of~\eqref{eq:NdeltaKcutSpeed} exists.

To see that the right-hand side of~\eqref{eq:NdeltaKcutSpeed} is a lower bound, note that~\eqref{eq:Ninequality} implies that for any $U \in A$,
\[
\NdKcut{A} \geq \NdKcut{U} \geq \NdeltaKcut{U}.
\]
The claim then follows from Lemma~\ref{le:NdKcutEntropy}.

Now we show that the right-hand side of~\eqref{eq:NdeltaKcutSpeed} is an upper bound.  By compactness, there exist $t \in \N$ and $\{W_1, \dots, W_t\} \subseteq A$ such that $\bigcup_{i=1}^t B_{\delta}(W_i) = A$.  It follows that
\begin{equation}\label{eq:deltaNet}
\NdeltaKcut{A} \leq \sum_{i=1}^t \NdeltaKcut[2\delta]{W_i}.
\end{equation}

Lemma~\ref{le:WeakRegEffective} implies that for every $m \geq 1$ and for every $i \in [t]$, there  exists an equipartition~$\SSS_i$ of~$\oi$ into $m$ parts such that 
\begin{equation}\label{eq:deltaNetStepfunction}
\dKcut(W_i, W_{i_{\SSS_i}}) < \dfrac{4}{\sqrt{\log_2 m}}.
\end{equation}
Hence, \eqref{eq:deltaNet} and Lemma~\ref{le:NdeltaKcutUpperBound} imply that
\begin{align}\label{eq:deltaNetMax}
\limsup_{n \to \infty} \dfrac{\log_k \NdeltaKcut{A}}{\binom{n}{2}} &\leq \max_{i \leq t} \limsup_{n \to \infty} \dfrac{\log_k \NdeltaKcut[2\delta]{W_i}}{\binom{n}{2}} \nonumber\\
&\leq  \max_{i \leq t} \limsup_{n \to \infty} \Ent(W_{i_{\SSS_i}}) + f(8m^2 \delta).
\end{align}
For each $m \geq 1$, let $\delta = 2^{-m}$, let $i(m)$ be the index that maximises the right-hand side of~\eqref{eq:deltaNetMax}, let $W'_m = W_{i(m)}$, and let $\SSS_m = \SSS_{i(m)}$.  By compactness, we may choose a convergent subsequence of the $W'_m$; let $W' \in A$ be such that $W'_m \to W'$.  It follows that
\begin{equation}\label{eq:newNdeltaKcutUB}
\limsup_{n \to \infty} \dfrac{\log_k \NdeltaKcut[2^{-m}]{A}}{\binom{n}{2}} \leq \Ent(W_{\SSS_m}) + f(8m^2 2^{-m}).
\end{equation}
By~\eqref{eq:deltaNetStepfunction}, $W_{\SSS_m} \to W'$ as $m \to \infty$.  Lemma~\ref{lemma: generalization of 3.3 in HJS} thus implies that 
\begin{equation*}\label{eq:partitionSemicontinuous}
\limsup_{m \to \infty} \Ent(W_{\SSS_m}) \leq \Ent(W').
\end{equation*}
Observe that $\NdeltaKcut{A}$ is increasing in $\delta$.  It follows from~\eqref{eq:newNdeltaKcutUB} that we have
\[
 \lim_{\delta \to 0}\limsup_{n \to \infty} \dfrac{\log_k \NdeltaKcut{A}}{\binom{n}{2}} =  \lim_{m \to \infty}\limsup_{n \to \infty} \dfrac{\log_k \NdeltaKcut[2^{-m}]{A}}{\binom{n}{2}} \leq \Ent(W') \leq \max_{W \in A} \Ent(W),
\]
as claimed.
\end{proof}

\subsection{Proofs of main results}\label{se:entropyProofs}

Now we prove Theorems \ref{th:SpeedEntropyIneq} and~\ref{theorem: HJS for decorated graphons}(ii).  As mentioned earlier, we omit the proof of Theorem~\ref{th:maxlimits}.

\begin{proof}[Proof of Theorem~\ref{th:SpeedEntropyIneq}.]
Let $\delta > 0$.  Observe that, by compactness, for all $n$~sufficiently large, if $G \in \QQ_n$, then $\deltaKcut(G, \limitsx{\QQ}) < \delta$.  It follows that $\abs{\QQ_n} \leq \abs{\QQ^L_n} \leq \NdeltaKcut{\limitsx{\QQ}}$.  The claimed result then follows from Lemma~\ref{le:NdeltaKcutSpeed}.
\end{proof}

\begin{proof}[Proof of Theorem~\ref{theorem: HJS for decorated graphons}(ii).]
By Theorem~\ref{th:SpeedEntropyIneq}, it is enough to show that if $\QQ$ is hereditary, then
\[
\liminf_{n \to \infty} \dfrac{\log_k \abs{\QQ_n}}{\binom{n}{2}} \geq \max_{\Gamma \in \limitsx{\QQ}} \Ent(\Gamma).
\]
Let $W \in \limitsx{\QQ}$.  Because $\QQ$ is hereditary, Lemma~\ref{le:wRandomInHeredProp} implies that $G(n, W) \in \QQ_n^L$ almost surely.  If we consider $G(n, W)$ as a random variable taking values in $\QQ_n^L$, then standard properties of entropy imply that
\[
\Ent\bigl(G(n, W)\bigr) \leq \log_k\abs{\QQ_n^L}.
\]
Theorem~\ref{th:wRandomEntropy} and~\eqref{eq:labelledUnlabelled} then imply that
\[
\Ent(W) \leq \liminf_{n \to \infty} \dfrac{\Ent\bigl(G(n, W)\bigr)}{\binom{n}{2}} \leq \liminf_{n \to \infty} \dfrac{\log_k \abs{\QQ_n}}{\binom{n}{2}},
\]
which is what we wanted to show.
\end{proof}

%

	\section{Concluding remarks}\label{section: concluding remarks}
\subsection{Entropy maximisation in the multicolour setting}
In the $2$-colour setting, the rough structure of entropy maximisers for hereditary properties is well-understood, via the choice number~$\chi_c$: given a hereditary property~$\PP$ with $\chi_c(\PP) = r$, partition the vertex sets into $r$ equal parts and define a template by giving the $r$-partite edges full entropy (i.e.\ free choice of their colour) and the other edges zero entropy (i.e.\ fix their colour, although different edges in the same part may have different colours). 
In particular, Theorem~\ref{theorem: alekseevbollobasthomason} implies that the set of possible entropy densities for hereditary properties is $\{0, 1/2, 2/3, 3/4, \cdots\}\cup\{1\}$.

By contrast, it is less clear what the set of possible values of entropy densities or the possible rough structure of entropy maximisers should be in the $k$-coloured setting for $k\geq 3$.  We are only aware of one partial result in this area: Alekseev and Sorochan~\cite{AlekseevSorochan:colored} showed that if $\PP$ is a symmetric hereditary property of $k$-coloured graphs, then either $\pi(\PP) = 0$ or $\pi(\PP) \geq (1/2)\log_k 2$.  Moreover, it is clear that the possible structures of entropy maximizers are much more varied than in the case $k  = 2$.  For example, we saw in the proof of Theorem~\ref{theorem: extremal entropy no rainbow k3} that if $\PP$ is the property of $3$-coloured graphs of containing no rainbow triangle, then $\Ent(\PP) = \log_3 2$ and that the unique maximum-entropy templates are those where the palette of each edge is a fixed pair of colours.  On the other hand, suppose that $\PP$ is the property of containing no monochromatic triangle in either colour~$2$ or colour~$3$.  The same argument as in the proof of Theorem~\ref{theorem: extremal entropy, no double triangle} shows that $\pi(\PP) = 1/2$ and that every unique maximum-entropy template~$t$ corresponds to a balanced partition of~$[n]$ where $t(e) = \{1\}$ for every edge within a partition class and $t(e) = [3]$ for every edge between partition classes.

\begin{problem}\label{problem: possible entropy density in multicolour setting}
	Let $k\in \N$ with $k\geq 3$. Determine the set of possible entropy densities of hereditary properties of $k$-colourings of~$K_n$ and the rough structure of entropy maximisers.
\end{problem}


\subsection{Contrasts between the graph limit and the container approaches}	
In this paper we have explored two approaches to counting, characterization and transference results  for hereditary properties: the hypergraph containers approach from Section~\ref{section: containers} and the entropy of graph limits approach from Section~\ref{section: decorated generalization of HJS}. We would be remiss to conclude the paper without comparing these two approaches. 

A clear advantage of the theory of containers is that it often gives more precise information: in particular, container theory has been successfully used to give good bounds on the size of `sparse' graph properties, such as $C_4$-free graphs or $K_{s,t}$-free graphs, see~\cite{BaloghSamotij11,BaloghWagner16,KleitmanWinston82,MorrisSaxton16}.  By comparison, the theory of (dense) graph limits sheds little light on such properties: any sequence of graphs with $o(n^2)$ edges tends to the constant graphon~$0$, and the results of~\cite{HJS} only show that any sparse graph property has size~$2^{o(n^2)}$.

We should note here that developing an appropriate limit theory for sequences of sparse graphs (to be precise, sequences of graphs with $o(n^2)$ edges but with average degree tending to infinity) constitutes an active field of research, but it is still nascent in several respects.  In particular, several different notions of convergence have been proposed (see~\cite{BollobasRiordan09,BorgsChayesCohnZhao14a,BorgsChayesCohnZhao14b,KKLS14}), and it is not yet clear which, if any, of these is the `right' one.  That said, it would be interesting to determine whether the results of any of the papers just cited would make it possible to extract information about the number and typical structure of graphs in a sparse property from the space of its appropriately-defined limits.

Secondly, the container approach adapts (relatively) straightforwardly to hypergraphs.  In contrast, extending the theory of graph limits to hypergraphs has proven to be very difficult, although some important steps have been taken~\cite{ElekSzegedy12,Zhao15}.  In particular, there is as yet no satisfactory generalization of the cut distance to hypergraphons.  (For further discussion, see~\cite[Section~23.3]{LovaszBook}.)  In view of these two points, it is natural to ask whether containers could pave the way for new developments in graph limit theory. In any case, this leads us to conclude this paper with a renewed expression of admiration for the power, applicability and elegance of the container theories of Balogh--Morris--Samotij and Saxton--Thomason.

\section*{Acknowledgements}
Victor Falgas-Ravry is grateful for an AMS-Simons award which allowed him to invite Andrew Uzzell to visit him and the two other authors at Vanderbilt University in November~2015, when the last stages of this research were carried out.  The authors are grateful to Caroline Terry for helpful remarks about~\cite{Terry16} and to Daniel Toundykov for Russian language assistance.

	\bibliographystyle{plain}
	\bibliography{containersgraphonbiblio}

\begin{thebibliography}{10}

\bibitem{Aldous}
D.~J. Aldous.
\newblock Exchangeability and related topics.
\newblock In {\em \'{E}cole d'\'et\'e de probabilit\'es de {S}aint-{F}lour,
  {XIII}---1983}, volume 1117 of {\em Lecture Notes in Math.}, pages 1--198.
  Springer, Berlin, 1985.

\bibitem{Alekseev93}
V.~E. Alekseev.
\newblock On the entropy values of hereditary classes of graphs.
\newblock {\em Discrete Math. Appl.}, 3(2):191--200, 1993.

\bibitem{AlekseevSorochan:colored}
V.~E. Alekseev and S.~V. Sorochan.
\newblock On the entropy of hereditary classes of colored graphs.
\newblock {\em Diskret. Mat.}, 12(2):99--102, 2000.

\bibitem{ABBM11}
N.~Alon, J.~Balogh, B.~Bollob{\'a}s, and R.~Morris.
\newblock The structure of almost all graphs in a hereditary property.
\newblock {\em J. Combin. Theory Ser. B}, 101(2):85--110, 2011.

\bibitem{AlonBaloghKeevashSudakov04}
N.~Alon, J.~Balogh, P.~Keevash, and B.~Sudakov.
\newblock The number of edge colorings with no monochromatic cliques.
\newblock {\em J. London Math. Soc.}, 70(2):273--288, 2004.

\bibitem{Austin08}
T.~Austin.
\newblock On exchangeable random variables and the statistics of large graphs
  and hypergraphs.
\newblock {\em Probab. Surv.}, 5:80--145, 2008.

\bibitem{Baber12}
R.~Baber.
\newblock Tur{\'a}n densities of hypercubes.
\newblock Preprint, \url{http://arxiv.org/abs/1201.3587}, 2012.

\bibitem{BaloghBollobasMorris06}
J.~Balogh, B.~Bollob{\'a}s, and R.~Morris.
\newblock Hereditary properties of ordered graphs.
\newblock In {\em Topics in discrete mathematics}, pages 179--213. Springer,
  2006.

\bibitem{BaloghBollobasMorris07b}
J.~Balogh, B.~Bollob{\'a}s, and R.~Morris.
\newblock Hereditary properties of combinatorial structures: posets and
  oriented graphs.
\newblock {\em J. Graph Theory}, 56(4):311--332, 2007.

\bibitem{BaloghBollobasMorris07}
J.~Balogh, B.~Bollob{\'a}s, and R.~Morris.
\newblock Hereditary properties of tournaments.
\newblock {\em Electronic J. Combin.}, 14(3):R60, 2007.

\bibitem{BBS04}
J.~Balogh, B.~Bollob{\'a}s, and M.~Simonovits.
\newblock The number of graphs without forbidden subgraphs.
\newblock {\em J. Combin. Theory Ser. B}, 91(1):1--24, 2004.

\bibitem{BBS09}
J.~Balogh, B.~Bollob{\'a}s, and M.~Simonovits.
\newblock The typical structure of graphs without given excluded subgraphs.
\newblock {\em Random Structures Algorithms}, 34(3):305--318, 2009.

\bibitem{BaloghHuLidickyLiu14}
J.~Balogh, P.~Hu, B.~Lidick{\'y}, and H.~Liu.
\newblock Upper bounds on the size of 4-and 6-cycle-free subgraphs of the
  hypercube.
\newblock {\em European J. Combin.}, 35:75--85, 2014.

\bibitem{BaloghMorrisSamotij15}
J.~Balogh, R.~Morris, and W.~Samotij.
\newblock Independent sets in hypergraphs.
\newblock {\em J. Amer. Math. Soc.}, 28(3):669--709, 2015.

\bibitem{BaloghSamotij11}
J.~Balogh and W.~Samotij.
\newblock The number of {$K\sb {s,t}$}-free graphs.
\newblock {\em J. Lond. Math. Soc. (2)}, 83(2):368--388, 2011.

\bibitem{BaloghWagner16}
J.~Balogh and A.~Zs. Wagner.
\newblock Further applications of the container method.
\newblock In A.~Beveridge, R.~J. Griggs, L.~Hogben, G.~Musiker, and P.~Tetali,
  editors, {\em Recent Trends in Combinatorics}, volume 159 of {\em The IMA
  Volumes in Mathematics and its Applications}, pages 191--213. Springer
  International Publishing, 2016.

\bibitem{BenevidesHoppenSampaio16}
F.~S. Benevides, C.~Hoppen, and R.~M. Sampaio.
\newblock Edge-colorings of graphs avoiding complete graphs with a prescribed
  coloring.
\newblock Preprint, \url{arxiv.org/abs/1605.08013}, 2016.

\bibitem{Bollobas07}
B.~Bollob{\'a}s.
\newblock Hereditary and monotone properties of combinatorial structures.
\newblock In A.~Hilton and J.~Talbot, editors, {\em Surveys in combinatorics
  2007}, volume 346 of {\em London Math. Soc. Lecture Note Ser.}, pages 1--39.
  Cambridge Univ. Press, Cambridge, 2007.

\bibitem{BollobasRiordan09}
B.~Bollob{\'a}s and O.~Riordan.
\newblock Metrics for sparse graphs.
\newblock In S.~Huczynska, J.~D. Mitchell, and C.~M. Roney-Dougal, editors,
  {\em Surveys in combinatorics 2009}, volume 365 of {\em London Math. Soc.
  Lecture Note Ser.}, pages 211--287. Cambridge Univ. Press, Cambridge, 2009.

\bibitem{BT}
B.~Bollob{\'a}s and A.~Thomason.
\newblock Hereditary and monotone properties of graphs.
\newblock In R.L. Graham and J.~Ne\v{s}et\v{r}il, editors, {\em The
  {M}athematics of {P}aul {E}rd{\"o}s {II}}, volume~14 of {\em Algorithms
  Combin.}, pages 70--78. Springer Berlin Heidelberg, 1997.

\bibitem{BorgsChayesCohnZhao14a}
C.~Borgs, J.~T. Chayes, H.~Cohn, and Y.~Zhao.
\newblock An {$L^p$} theory of sparse graph convergence {I}: limits, sparse
  random graph models, and power law distributions.
\newblock Preprint, \url{http://arxiv.org/abs/1401.2906}, 2014.

\bibitem{BorgsChayesCohnZhao14b}
C.~Borgs, J.~T. Chayes, H.~Cohn, and Y.~Zhao.
\newblock An {$L^p$} theory of sparse graph convergence {II}: {LD} convergence,
  quotients, and right convergence.
\newblock Preprint, \url{http://arxiv.org/abs/1408.0744}, 2014.

\bibitem{BCLSV08}
C.~Borgs, J.~T. Chayes, L.~Lov{\'a}sz, V.~T. S{\'o}s, and K.~Vesztergombi.
\newblock Convergent sequences of dense graphs {I}: Subgraph frequencies,
  metric properties and testing.
\newblock {\em Adv. Math.}, 219(6):1801--1851, 2008.

\bibitem{ChatterjeeVaradhan11}
S.~Chatterjee and S.~R.~S. Varadhan.
\newblock The large deviation principle for the {E}rd{\H o}s-{R}\'enyi random
  graph.
\newblock {\em European J. Combin.}, 32(7):1000--1017, 2011.

\bibitem{Conlon14}
D.~Conlon.
\newblock Combinatorial theorems relative to a random set.
\newblock {\em Proceedings of the International Congress of Mathematicians,
  Seoul 2014}, to appear.

\bibitem{ConlonGowers10}
D.~Conlon and W.~T. Gowers.
\newblock Combinatorial theorems in sparse random sets.
\newblock {\em Ann. Math.}
\newblock to appear.

\bibitem{ConlonGowersSamotijSchacht14}
D.~Conlon, W.~T. Gowers, W.~Samotij, and M.~Schacht.
\newblock On the {K{\L}R} conjecture in random graphs.
\newblock {\em Israel J. Math.}, 203(1):535--580, 2014.

\bibitem{DotsonNagle09}
R.~Dotson and B.~Nagle.
\newblock Hereditary properties of hypergraphs.
\newblock {\em J. Combin. Theory Ser. B}, 99(2):460--473, 2009.

\bibitem{ElekSzegedy12}
G.~Elek and B.~Szegedy.
\newblock A measure-theoretic approach to the theory of dense hypergraphs.
\newblock {\em Adv. Math.}, 231(3--4):1731--1772, 2012.

\bibitem{Erdos67}
P.~Erd{\H{o}}s.
\newblock Some recent results on extremal problems in graph theory. {R}esults.
\newblock In {\em Theory of {G}raphs ({I}nternat. {S}ympos., {R}ome, 1966)},
  pages 117--123 (English); 124--130 (French). Gordon and Breach, New York,
  1967.

\bibitem{Erdos74}
P.~Erd{\H o}s.
\newblock Some new applications of probability methods to combinatorial
  analysis and graph theory.
\newblock In {\em Proc. of the Fifth Southeastern Conference on Combinatorics,
  Graph Theory and Computing (Florida Atlantic Univ., Boca Raton, Fla., 1974},
  pages 39--51. University of Calgary, Department of Mathematics, Statistics
  and Computing Science, 1974.

\bibitem{Erdos84}
P.~Erd{\H o}s.
\newblock On some problems in graph theory, combinatorial analysis and
  combinatorial number theory.
\newblock {\em Graph Theory and Combinatorics (Cambridge, 1983), Academic
  Press, London}, pages 1--17, 1984.

\bibitem{ErdosFranklRodl86}
P.~Erd{\H o}s, P.~Frankl, and V.~R{\"o}dl.
\newblock The asymptotic number of graphs not containing a fixed subgraph and a
  problem for hypergraphs having no exponent.
\newblock {\em Graphs Combin.}, 2(1):113--121, 1986.

\bibitem{ErdosKleitmanRothschild76}
P.~Erd{\H{o}}s, D.~J. Kleitman, and B.~L. Rothschild.
\newblock Asymptotic enumeration of {$K_n$}-free graphs.
\newblock In {\em Internat. Colloq. Combin.} Atti Convegni Lincei (Rome), 1976.

\bibitem{ErdosStone46}
P.~Erd{\"o}s and A.~H. Stone.
\newblock On the structure of linear graphs.
\newblock {\em Bull. Amer. Math. Soc.}, 52:1087--1091, 1946.

\bibitem{FriedgutRodlSchacht10}
E.~Friedgut, V.~R{\"o}dl, and M.~Schacht.
\newblock Ramsey properties of random discrete structures.
\newblock {\em Random Structures Algorithms}, 37(4):407--436, 2010.

\bibitem{FK99}
A.~Frieze and R.~Kannan.
\newblock Quick approximation to matrices and applications.
\newblock {\em Combinatorica}, 19(2):175--220, 1999.

\bibitem{FurediKundgen02}
Z.~F{\"u}redi and A.~K{\"u}ndgen.
\newblock Tur{\'a}n problems for integer-weighted graphs.
\newblock {\em J. Graph Theory}, 40(4):195--225, 2002.

\bibitem{Galvin14}
D.~Galvin.
\newblock Three tutorial lectures on entropy and counting.
\newblock Preprint, \url{arxiv.org/abs/1406.7872}, 2014.

\bibitem{GreenTao08}
B.~Green and T.~Tao.
\newblock The primes contain arbitrarily long arithmetic progressions.
\newblock {\em Ann. Math. (2)}, 167(2):481--547, 2008.

\bibitem{HJS}
H.~Hatami, S.~Janson, and B.~Szegedy.
\newblock Graph properties, graph limits, and entropy.
\newblock Preprint, \url{http://arxiv.org/abs/1312.5626}, 2013.

\bibitem{HoppenLefmannOdermann}
C.~Hoppen, H.~Lefmann, and K.~Odermann.
\newblock A rainbow {E}rd{\H{o}}s-{R}othschild problem.
\newblock {\em Electron. Notes Discrete Math.}, 49:473--480, 2015.

\bibitem{Ishigami07}
Y.~Ishigami.
\newblock The number of hypergraphs and colored hypergraphs with hereditary
  properties.
\newblock Preprint, \url{arxiv.org/abs/0712.0425}, 2007.

\bibitem{Jan13}
S.~Janson.
\newblock Graphons, cut norm and distance, rearrangements and coupling.
\newblock {\em New York J. Math. Monographs}, 24:1--76, 2013.

\bibitem{Jan16}
S.~Janson.
\newblock Graph limits and hereditary properties.
\newblock {\em European J. Combin.}, 52(part B):321--337, 2016.

\bibitem{JohnsonEntringer89}
K.~A. Johnson and R.~Entringer.
\newblock Largest induced subgraphs of the {$n$}-cube that contain no
  {$4$}-cycles.
\newblock {\em J. Combinatorial Theory Ser. B}, 46(3):346--355, 1989.

\bibitem{KleitmanWinston82}
D.~J. Kleitman and K.~J. Winston.
\newblock On the number of graphs without {$4$}-cycles.
\newblock {\em Discrete Math.}, 41(2):167--172, 1982.

\bibitem{KLR}
Y.~Kohayakawa, T.~{\L}uczak, and V.~R{\"o}dl.
\newblock On {$K\sp 4$}-free subgraphs of random graphs.
\newblock {\em Combinatorica}, 17(2):173--213, 1997.

\bibitem{KohayakawaNagleRodl03}
Y.~Kohayakawa, B.~Nagle, and V.~R{\"o}dl.
\newblock Hereditary properties of triple systems.
\newblock {\em Combin. Probab. Comput.}, 12(2):155--189, 2003.

\bibitem{Kostochka76}
E.~A. Kostochka.
\newblock Piercing the edges of the n-dimensional unit cube.
\newblock {\em Diskret. Analiz Vyp}, 28(223):55--64, 1976.

\bibitem{KuhnOsthusTownsendZhao14}
D.~K{\"u}hn, D.~Osthus, T.~Townsend, and Y.~Zhao.
\newblock On the structure of oriented graphs and digraphs with forbidden
  tournaments or cycles.
\newblock Preprint \url{arxiv.org/abs/1404.6178}, 2014.

\bibitem{KKLS14}
D.~Kunszenti-Kov\'acs, L.~Lov\'asz, and B.~Szegedy.
\newblock Multigraph limits, unbounded kernels, and {B}anach space decorated
  graphs.
\newblock Preprint, \url{http://arxiv.org/abs/1406.7846}, 2014.

\bibitem{LovaszBook}
L.~Lov{\'a}sz.
\newblock {\em Large Networks and Graph Limits}, volume~60 of {\em Amer. Math.
  Soc. Colloq. Publ.}
\newblock Amer. Math. Soc., Providence, RI, 2012.

\bibitem{LS06}
L.~Lov{\'a}sz and B.~Szegedy.
\newblock Limits of dense graph sequences.
\newblock {\em J. Combin. Theory Ser. B}, 96(6):933--957, 2006.

\bibitem{LS:analyst}
L.~Lov{\'a}sz and B.~Szegedy.
\newblock Szemer\'edi's lemma for the analyst.
\newblock {\em Geom. Funct. Anal.}, 17(1):252--270, 2007.

\bibitem{LS:decorated}
L.~Lov{\'a}sz and B.~Szegedy.
\newblock Limits of compact decorated graphs.
\newblock Preprint, \url{http://arxiv.org/abs/1010.5155}, 2010.

\bibitem{MorrisSaxton16}
R.~Morris and D.~Saxton.
\newblock The number of {$C_{2\ell}$}-free graphs.
\newblock {\em Adv. Math.}, 298:534--580, 2016.

\bibitem{MubayiTerry16a}
D.~Mubayi and C.~Terry.
\newblock An extremal graph problem with a transcendental solution.
\newblock Preprint, \url{http://arxiv.org/abs/1607.07742}, 2016.

\bibitem{MubayiTerry16b}
D.~Mubayi and C.~Terry.
\newblock Structure and enumeration for some families of multigraphs.
\newblock Preprint, 2016.

\bibitem{NagleRodlSchacht06}
B.~Nagle, V.~R{\"o}dl, and M.~Schacht.
\newblock Extremal hypergraph problems and the regularity method.
\newblock {\em Topics in discrete mathematics}, pages 247--278, 2006.

\bibitem{PikhurkoStadenYilma16}
O.~Pikhurko, K.~Staden, and Z.~Yilma.
\newblock The {E}rd{\H o}s--{R}othschild problem on edge-colourings with
  forbidden monochromatic cliques.
\newblock Preprint, \url{arxiv.org/abs/1605.05074}, 2016.

\bibitem{PS92}
H.~J. Pr{\"o}mel and A.~Steger.
\newblock Excluding induced subgraphs. {III}. {A} general asymptotic.
\newblock {\em Random Structures Algorithms}, 3(1):19--31, 1992.

\bibitem{Razborov07}
A.~A. Razborov.
\newblock Flag algebras.
\newblock {\em J. Symbolic Logic}, 72(4):1239--1282, 2007.

\bibitem{Razborov13}
A.~A. Razborov.
\newblock What is{$\ldots$}a flag algebra?
\newblock {\em Notices Amer. Math. Soc.}, 60(10):1324--1327, 2013.

\bibitem{Sapozhenko87}
A.~A. Sapozhenko.
\newblock On the number of connected subsets with given cardinality of the
  boundary in bipartite graphs.
\newblock {\em Metody Diskret. Analiz.}, (45):42--70, 96, 1987.

\bibitem{Sapozhenko01}
A.~A. Sapozhenko.
\newblock On the number of independent sets in extenders.
\newblock {\em Diskret. Mat.}, 13(1):56--62, 2001.

\bibitem{SaxtonThomason15}
D.~Saxton and A.~Thomason.
\newblock Hypergraph containers.
\newblock {\em Invent. Math.}, 201(3):925--992, 2015.

\bibitem{SaxtonThomason16}
D.~Saxton and A.~Thomason.
\newblock Simple containers for simple hypergraphs.
\newblock {\em Combin. Probab. Comput.}, 25(3):448--459, 2016.

\bibitem{Schacht2009}
M.~Schacht.
\newblock Extremal results for random discrete structures.
\newblock {\em Ann. Math.}, 184(2):331--363, 2016.

\bibitem{Shannon48}
C.~Shannon.
\newblock A mathematical theory of communication.
\newblock {\em Bell System Tech. J.}, 27:379--423, 1948.

\bibitem{Simonovits68}
M.~Simonovits.
\newblock A method for solving extremal problems in graph theory, stability
  problems.
\newblock In {\em Theory of {G}raphs ({P}roc. {C}olloq., {T}ihany, 1966)},
  pages 279--319. Academic Press, New York, 1968.

\bibitem{Terry16}
C.~Terry.
\newblock Structure and enumeration theorems for hereditary properties in
  finite relational languages.
\newblock Preprint, \url{http://arxiv.org/abs/1607.04902}, 2016.

\bibitem{Zhao15}
Y.~Zhao.
\newblock Hypergraph limits: a regularity approach.
\newblock {\em Random Structures Algorithms}, 47(2):205--226, 2015.

\end{thebibliography}
	\appendix
	\section{Appendix}
	\begin{proof}[Proof of Theorem~\ref{theorem: stability no rainbow K3}]
		Fix $\varepsilon>0$, and let $t$ be an $n$-vertex $3$-colouring template satisfying properties (i) and~(ii) of the statement for some $\delta>0$ and $n_0\in \N$ to be specified later, and $n\geq n_0$. Our proof is a (lengthy) exercise in stability analysis --- essentially, we shall prove an approximate version of Observation~\ref{observation: rainbow K3}, and then run through the proof of Theorem~\ref{theorem: extremal entropy no rainbow k3} replacing each `for all pairs' by a `for almost all pairs'.

		Let $C_0=C_0(\mathcal{P})$ be the constant whose existence is asserted by our supersaturation result, Lemma~\ref{lemma: supersaturation}:  for all $\eta>0$ there exists $n_1(\eta, \mathcal{P})$ such that for all $n\geq n_1$, if $t$ is a $3$-colouring template on $n$ vertices which can realise at most $\eta\binom{n}{3}$ rainbow triangles, then $\Ent(t)\leq \left(\pi(\mathcal{P})+ C_0\eta\right)\binom{n}{2}$.

		Let $e_3'$ be the number of edges $e=\{u,v\}\in E(K_n)$ for which there are at least~$\delta n$ vertices~$x\in V(K_n)\setminus\{u,v\}$ for which $\vert t(\{x,u\})\vert+ \vert t(\{x,v\})\vert >2$. For each such edge~$e$ and each such vertex~$x$, there is at least~one rainbow triangle which can be realised inside $e\cup\{x\}$. Each such triangle is counted at most~$3$ times, so that in total we must have at least $\frac{e_3'\delta n}{3}<\delta \binom{n}{3}$ rainbow triangles, and in particular we must have $e_3'<\frac{\delta}{2}n^2$.

		Now let $e_3''$ denote the number of edges $e=\{u,v\}$ for which there are at most $\delta n$ vertices $x$ with $\vert t(\{x,u\})\vert + \vert t(\{x,v\})\vert >2$. We shall choose $\delta$~sufficiently small to ensure that (a) $(1-200C_0\delta)^2>2/3$ and (b) $\delta <\frac{2\pi(\mathcal{P}) -1/50}{200(C_0 +2)}$ (we can certainly do that since the value of the constant~$C_0$ does not depend on $\delta$).

		Suppose $n>3n_1(2\delta)$. We claim that $e_3''<200 (C_0 +1)\delta n^2$. Indeed suppose not. Then we can find a set $E_3''$ of at least $200 (C_0 +1)\delta n^2/2n = 100 (C_0 +1)\delta n:=cn$ vertex-disjoint edges $e=\{u,v\}$ with $t(e)=3$ and $\vert t(\{x,u\})\vert+\vert t(\{x,v\})\vert=2$ for all but at most $\delta n$ vertices $x$. Remove from $K_n$ the pairs of vertices $e=\{u,v\}$ from $E_3''$ one by one. This leaves us with a graph on $n'= n -2c n$ vertices, which by (a) and our assumption on $n$ is strictly greater than~$n_1(2\delta)$.

		Let $t'$ denote the subtemplate of $t$ induced by the remaining vertices. Clearly $t'$ can realise at most $\delta\binom{n}{3}$ rainbow triangles, which by (a) is at most $2\delta \binom{n'}{3}$. Now Lemma~\ref{lemma: supersaturation} and the fact that $n'>n_1(2\delta)$ implies that 
		\begin{align}\label{equation: upper bound on entropy of t'}
		\Ent({t'})&\leq \left(\pi(\mathcal{P})+ C_02\delta\right)\binom{n'}{2}\leq \pi(\mathcal{P})\binom{n}{2}+\frac{C_02\delta}{2}n^2 -2c(1-c)\pi(\mathcal{P})n^2.
		\end{align} 
		On the other hand, each of the edges $e$ from $E_3''$ we removed decreased the entropy by at most $\delta n $, so we have the following lower bound on $\Ent(t')$:
		\begin{align}\label{equation: lower bound on entropy of t'}
		\Ent(t')&\geq \Ent(t)-c\delta n^2\geq\pi(\mathcal{P})\binom{n}{2}-(c\delta +\delta)\frac{n^2}{2}.
		\end{align}
		Bringing the two bounds (\ref{equation: upper bound on entropy of t'}) and (\ref{equation: lower bound on entropy of t'}) together and cancelling terms as appropriate, we get
		\[ -\frac{c\delta}{2} -\frac{\delta}{2}\leq \frac{C_02\delta}{2}-2c(1-c)\pi(\mathcal{P}).\]
		Rearranging yields
		\[ c\left(2(1-c)\pi(\mathcal{P})-\frac{\delta}{2}\right)\leq \frac{\delta}{2}(1+2C_0).\]
		Since $c=100(C_0+1)$, this contradicts our assumption (b) on $\delta $. It follows that $e_3''<200(C_0+1)\delta n^2$, as claimed. Thus in total, there are at most~$e_3'+e_3''=(\delta/2 +200(C_0+1)\delta)n^2:=C_2 \delta n^2$ edges~$e$ with $\vert t(e)\vert=3$.

		We now move on to bounding the number $e_1$ of edges $e$ with $\vert t(e)\vert =1$. We have
		\[\bigl(\pi(\mathcal{P})-\delta\bigr)\binom{n}{2}\leq \Ent(t)\leq \pi(\mathcal{P})\biggl(\binom{n}{2} -e_1\biggr) + e_3' + e_3'',\]
		which together with our bound on $e_3'+e_3''$ implies that
		\[e_1< \frac{1}{\pi(\mathcal{P})} \left(\frac{1}{2}+C_2\right)\delta n^2.\]
		In particular, all but at most 
		\[\left(\frac{1}{\pi(\mathcal{P})} \left(\frac{1}{2}+C_2\right) +C_2\right)\delta n^2:=C_3\delta n^2\] edges $e$ have $\vert t(e)\vert=2$.

		Finally we turn to the edges assigned two colours by $t$. For each pair of colours~$A\in [3]^{(2)}$, let $V_A$ denote the collection of vertices incident to at least~$\delta^{1/3} n$ edges in colour~$A$. For any $A\neq B$, each vertex in $V_A\cap V_B$ gives rise to at least~$\delta^{2/3} n^2$ distinct rainbow triangles, whence 
		\[\frac{\vert V_A\cap V_B\vert}{3}\leq \delta \binom{n}{3},\]
		implying $\vert A\cap V_B\vert \leq  \delta^{1/3} n/2$. Suppose we had $\vert V_A\vert$ and $\vert V_B\vert$ both greater than~$(\sqrt{C_3}+3)\delta^{1/3}n$ for some colour pairs~$A\neq B$, and let $C$ denote the third colour pair from $[3]$. Then all but at most~$\delta^{1/3} n$ vertices in $A$ are incident to at most~$2\delta^{1/3} n$ edges whose $t$-colour assignment is $B$ or~$C$. In particular such vertices~$a$ must be incident to at least $\vert V_B\vert -\vert V_B\cap V_A\vert -2\delta^{1/3} n$ edges~$ab$ with $b\in V_B\setminus V_A$ and $t(a,b)\notin \{A,B,C\}$. This gives at least
		\[ \left(\vert V_A\vert - \delta^{1/3} n\right) \left(\vert V_B\vert -\vert V_B\cap V_A\vert -2\delta^{1/3} n\right)\geq (\sqrt{C_3}+1)\delta^{1/3}n \sqrt{C_3} \delta^{1/3}n>C_3\delta n^2\]
		edges~$e$ with $\vert t(e)\vert\neq 2$, a contradiction. It follows that there is at most~one colour pair, say $A$, with $\vert V_A\vert\geq (\sqrt{C_3}+3)\delta^{1/3}n$. Let $B$,~$C$ denote the two other colour pairs, and $e_B$,~$e_C$ the number of edges~$e$ with $t(e)=B$ and $t(e)=C$ respectively. By the definition of $V_B$, we have 
		\[e_B\leq \vert V_B\vert n/2 + (n-\vert V_B\vert )\delta^{1/3}n/2< \frac{(\sqrt{C_3}+4)\delta^{1/3}}{2}n^2,\]
		and similarly $e_C \leq (\sqrt{C_3}+4)\delta^{1/3}n^2/2$. We have thus shown that all but at most~$(C_3\delta +(\sqrt{C_3}+4)\delta^{1/3})n^2$ edges $e\in E(K_n)$ have $t(e)\neq A$. Picking $\delta=\delta(\varepsilon)$~sufficiently small (and $n_0\geq 3n_1(2\delta)$), this is less than $\varepsilon \binom{n}{2}$, proving the theorem.
	\end{proof}
\end{document}